\documentclass[12pt, a4, fleqn]{article}
\usepackage{geometry}                
\usepackage{graphicx, amsthm, amsfonts, amsmath}
\usepackage{amssymb}
\usepackage{arydshln, ulem}
\usepackage{epstopdf, enumerate, tikz, subfigure}
\usepackage{pdflscape, multirow, tablefootnote}
\usepackage{threeparttable}
\usepackage{booktabs, color}
\theoremstyle{theorem}
\newtheorem{theorem}{Theorem}[section]
\newtheorem{proposition}[theorem]{Proposition}
\newtheorem{lemma}[theorem]{Lemma}
\theoremstyle{definition}

\newtheorem{remark}[theorem]{Remark}
\newcommand\red[1]{\textcolor{red}{#1}}
\newcommand\blue[1]{\textcolor{blue}{#1}}
\newcommand\black[1]{\textcolor{black}{#1}}

\definecolor{mygreen}{RGB}{0,100,0}

\title{Direct and Indirect Competition in a Two-Prey-One-Predator System \\
- {\large Global Dynamics, Chaos, and the Paine Conjecture}}
\author{Mingchieh Huang and Tinghui Yang}
\numberwithin{equation}{section}
\begin{document}
\maketitle
\begin{abstract}
This study investigates the global dynamics of a three-dimensional Lotka–Volterra system describing two competing prey species and one predator, with direct competition between the prey explicitly incorporated. After rescaling the system, we systematically classify the parameter space under a symmetry assumption into four generic categories comprising 84 candidate cases. By rigorously analyzing the compatibility of the associated parameter conditions, we prove that 28 cases are impossible, leaving 56 meaningful subcases for detailed analysis. This classification provides a comprehensive framework for characterizing the global dynamics of the system. We analyze the existence and stability of boundary equilibria and the positive equilibrium, including the global asymptotic stability and bistability of boundary equilibria. In particular, we derive necessary and sufficient conditions for the local asymptotic stability of the positive equilibrium. Our analysis also identifies parameter regimes associated with periodic solutions and complex dynamics. Systematic numerical investigations reveal numerous parameter combinations exhibiting chaotic dynamics within the subcase IV-(iv)-B1. Furthermore, by establishing an explicit parameter correspondence between Gilpin’s model and ours, we show that the parameter values associated with Gilpin’s chaotic dynamics fall within this subcase. Our numerical results also suggest a period-doubling route to chaos.

From a biological perspective, we introduce the concept of extensive negative competition, which extends the notion of negative competition between the two prey species to account for the combined effects of the other prey and the predator. Together with trade-off competition and apparent competition, this concept provides a framework for interpreting the interactions among the three species and their consequences for community structure. Under the symmetry assumption, we identify the dominant prey species through the birth-to-consumption ratio. Our classification shows that the dynamics of the full three-dimensional system are largely consistent with those of the associated two-dimensional prey-competition subsystem, with only six notable exceptions. Through rigorous analysis of these exceptional cases, we demonstrate how the interplay among extensive negative competition, trade-off competition, and apparent competition helps explain differences between the dynamics of the full system and those of the prey-only subsystem. These findings further clarify how predator-mediated interactions can alter competitive outcomes and facilitate the persistence of competing prey species, providing mathematical support for Paine's conjecture within the framework of the model.
\end{abstract}

\section{Introduction}
In this work, motivated by \cite{Holt1977, Gilpin1979, Hsu1981}, we consider the following two-prey-one-predator models with direct competition between the prey in the Lotka-Volterra functional response,
\begin{equation}\label{theorigmodel}
\left\{
\begin{aligned}
\frac{dU_1}{ds} &= \bar r_1U_1(1-\frac{U_1}{K_1})-\bar a_{12}U_1U_2- \bar a_{13}U_1V,\\
\frac{dU_2}{ds} &= \bar r_2U_2(1-\frac{U_2}{K_2})-\bar a_{21}U_1U_2- \bar a_{23}U_2V,\\
\frac{d V}{ds} &= V(\bar a_{31}U_1+\bar a_{32}U_2-D).
\end{aligned}
\right.
\end{equation}
Here, all parameters are positive, and two prey obey the logistic growth rule with intrinsic birth rates $\bar r_i$ and environmental capacities $K_i$ for $i=1, 2$. The coefficients $\bar a_{i3}$ and  $\bar a_{3i}$ are the rate of consumption of $U_i$ and measure the contribution of the basal resource $U_i$ to the growth of consumer $V$, respectively. The direct competition between $U_1$, $U_2$ are $\bar a_{12}$, $\bar a_{21}$, and the mortality of $V$ is presented by $D$.

With the rescaling,
\begin{align*}
u_1&=U_1/K_1, & u_2&=U_2/K_2, & v&=V/D,\\
t&=Ds, & r_1&=\frac{\bar r_1}{D},& r_2&=\frac{\bar r_2}{D},\\
a_{12}&=\frac{\bar a_{12}K_2}{D},&a_{31}&=\frac{\bar a_{31}K_1}{D},&a_{13}&=\bar a_{13},\\
a_{21}&=\frac{\bar a_{21}K_1}{D},&a_{32}&=\frac{\bar a_{32}K_2}{D},&a_{23}&=\bar a_{23},
\end{align*}
we have a simpler model with 8 parameters,
\begin{equation}\label{themodel}
\left\{
\begin{aligned}
\dot u_1 &= r_1u_1(1-u_1)-a_{12}u_1u_2-a_{13}u_1v,\\
\dot u_2 &= r_2u_2(1-u_2)-a_{21}u_1u_2-a_{23}u_2v,\\
\dot v &= v(a_{31}u_1+a_{32}u_2-1),
\end{aligned}
\right.
\end{equation}
where $\dot ~$ means $\frac{d}{dt}$.

In the past decades, extensive research has examined the dynamics of the Lotka-Volterra-type predator-prey equation, which serves as a prototype for more complex functional responses. Moreover, three-species models serve as a crucial bridge between simple two-species interactions and the highly complex ecological communities found in the real world. It is well known that predation and competition, interactions between species, are the primary drivers of population regulation and evolutionary pressure \cite{Van1973, Sole2022}. Hence, it is interesting and important to investigate the global dynamics of system \eqref{themodel}. 

Some early works should be mentioned \cite{Hutchinson1961, Paine1966, Parrish1970, Cramer1972, Fujii1977, Gilpin1979, Hsu1981, Yang2025}. 
The Paine conjecture (often called the keystone species hypothesis) \cite{Paine1966} suggests that an ecosystem's species diversity is directly related to its top predators. It has been suggested \cite{Hutchinson1961, Paine1966} that, in some communities, although a particular trophic level would, in isolation, be unstable due to competition, predation from other levels can lead to a stable system. You can see that the system \eqref{themodel} is the simplest case with multiple prey for verifying the Paine conjecture. However, system \eqref{themodel} without direct competition between prey, that is $a_{12}$=$a_{21}$=0, it was shown that \cite{Yang2025} $u_1$-$u_2$-$v$ can coexist if the apparent competition between $u_1$ and $u_2$ is not extreme. That is, for each prey $u_i$, the ratio of the birth rate $r_i$ to the consumption rate by $v$ is almost the same. Moreover, the positive equilibrium is globally asymptotically stable if it exists \cite{Yang2025}. Therefore, this top-down regulation is dependent not only on predation but also heavily on competition among prey.

Numerical simulations \cite{Parrish1970, Cramer1972} of the two-prey-one-predator models \eqref{themodel}, using some parameter values, illustrated and supported the Paine conjecture. In addition, periodic solutions or complex dynamics of  \eqref{themodel} were also found in  \cite{Fujii1977, Gilpin1979}. Mathematically, Hsu \cite{Hsu1981} gave a theoretical analysis of the two-prey-one-predator models by assuming that the predator $V$ can survive in the $u_1$-$v$ and $u_2$-$v$ subcommunities. However, the parameters of numerical simulations of \eqref{themodel} exhibit interesting and complex dynamics \cite{Parrish1970, Cramer1972, Fujii1977, Gilpin1979} that fall outside the range of Hsu's results. So it is important to know what dynamics the system \eqref{themodel} will exhibit if Hsu's assumptions do not hold.

Biologically, the model \eqref{themodel} can be regarded as a unified framework that includes both the resource competition mechanism like the classical $R^*$ rule and the predator-induced interaction mechanism associated with the $P^*$ rule \cite{Holt1994}, as well as the predator-mediated apparent competition \cite{Hsu1981}. However, we refer to the direct competition between prey species as {\it negative competition} rather than resource competition. Resource competition, the interaction between two consumers competing for a limiting resource, induces competitive exclusion, and only one winner emerges. However, there could be a win-win outcome in the negative competition for prey for \eqref{themodel}. Moreover,  the $P^*$ rule describes predator-mediated competition, assuming each prey species has a $P_i^*$, the predator level at which the prey can just survive. A smaller $P_i^*$ usually indicates a greater tolerance for predators and an advantage in apparent competition. Parameters $a_{13}$ and $a_{23}$ in \eqref{themodel} represent this competition mechanism.

Building on the above distinction between negative competition and apparent competition, we further introduce the concept of {\it extensive negative competition}. Classical negative competition captures only the direct, pairwise suppression that one prey exerts on the other through the competition coefficients $a_{12}$ and $a_{21}$. In the three-species system \eqref{themodel}, however, the effective pressure felt by a focal prey species is jointly shaped by the other prey and by the shared predator: consumption by the predator directly reduces the focal prey, while the other prey species both competes for the same growth potential and channels a common predation pressure onto the focal species. Extensive negative competition therefore extends the pairwise notion of negative competition to the combined effect of the other prey and the predator on a focal prey species, which is measured through the joint combinations of the parameters $(r_{i}, a_{i3}, a_{ij}), i, j= 1,2$. This unified viewpoint allows the direct (prey-prey) and the indirect (prey-predator) pathways of suppression to be described within a single framework, and it is this concept, together with trade-off competition and apparent competition, that we use to interpret the parameter-dependent transitions of the system and the differences between the full system and its prey-only subsystem. We will see that the interaction of these three competitive mechanisms produces a wide variety of dynamic behaviors, including competitive exclusion, coexistence, bistable states, periodic oscillations, and complex dynamics. 

Alternatively, other three-species models have long been of considerable importance. In \cite{Hsu2015}, Hsu et al. discussed the intraguild predation model, including a specialist predator, an omnivorous predator, and a prey. The parameter space is divided into four categories, and the global dynamics in each case are analyzed in detail. Various results are obtained and proved, including species extinction, global stability, bistability, and periodic solutions. It is also fully discussed that, within a certain parameter range, the system can generate periodic solutions via Hopf bifurcation and enter a chaotic state via period-doubling bifurcation. The research methods, calculation methods, and key theories in this paper primarily follow the comprehensive research model established in \cite{Hsu2015}.

Our investigations, while retaining all direct and indirect competition parameters from \eqref{theorigmodel}, conduct a comprehensive study of the dynamic changes arising from fluctuations in these parameters. However, we do not analyze a single parameter, as interactions among species in an ecosystem are unlikely to be influenced by a single factor. Our investigation found that establishing the following parameter relationships better reveals the complex dynamic changes. We use the relationship between $r_1a_{23}$ and $r_2a_{13}$ to represent the {\it apparent competition}, use the relationships between ($r_1a_{23}$ and $a_{12}a_{13}$) and ($a_{21}a_{13}$ and $r_2a_{13}$) to represent the {\it negative competition}, use the relationships between ($r_1a_{23}$ and $a_{21}a_{13}$) and ($a_{12}a_{13}$ and $r_2a_{13}$) to represent the {\it trade-off competition}, and use the inequality \eqref{u1venc} to represent the {\it extensive negative competition} between prey. However, given the extensive existing research on apparent competition, we will focus our main investigation on negative competition, trade-off competition, and extensive negative competition to study the global dynamics of \eqref{themodel}. Because these negative-competition pathways are in fact mediated by the predator through the consumption terms $a_{13}$ and $a_{23}$, we consistently treat negative competition in its extensive form introduced above, so that the classification below simultaneously encodes the direct prey-prey suppression and the predator-mediated indirect effects. This makes the extensive negative competition term the natural counterpart of apparent competition in our analysis. This interpretation provides a biological mechanism underlying the observed parameter-dependent transitions between different dynamical regimes. We will not presuppose any results; rather, we will present all results in full, providing rigorous mathematical proofs and numerical analyses to verify their credibility.

We organize this work as follows. In Section 2, we recall results for classical two-species competition systems and establish the local stability of boundary equilibria of \eqref{themodel}. In Section 3, we classify all parameters into four generic categories with a total of nine subcases under the assumption $C_1=r_1a_{23}>C_2=r_2a_{13}$, which means the apparent competition of $u_1$ is superior to that of $u_2$ \cite{Holt1977, Holt1997}. Five global dynamics of \eqref{themodel} are confirmed from these 9 categories. In particular, one case needs to be considered in more detail; it is presented in Table 3. In Section 4, for the other 4 categories, we present a table for each category, as we explore and classify the complicated subcases in detail. There are 84 subcases, and the corresponding global results are discussed individually. Finally, some numerical computations and biological implications are given. In particular, we highlight how extensive negative competition, together with trade-off and apparent competition, provides a unified biological interpretation of the transitions among coexistence, exclusion and bistability, and how it underpins the mathematical support for the Paine keystone-species conjecture.
\section{Preliminary}
This section recalls some well-known results of two-dimensional subsystems of \eqref{themodel}. They can help us to investigate the local and global dynamics of \eqref{themodel} in $\mathbb{R}^3$.

\subsection{Some well known two-dimensional results of \eqref{themodel}}
It is obvious that the three boundaries,
$H_{1v}=\{(u_1, 0, v): u_1, v>0\}$, $H_{2v}=\{(0, u_2, v):u_2, v>0\}$ and $H_{12}=\{(u_1, u_2, 0): u_1, u_2>0\}$,
of the positive cone, $\mathbb{R}^3_+=\{(u_1, u_2, v):u_1, u_2, v>0\}$, which are invariant with respective to \eqref{themodel}. Solutions of \eqref{themodel} are nonnegative/positive with nonnegative/positive initial conditions. Hence, we only consider solutions of \eqref{themodel} with initial conditions in $\mathbb{R}^3_+$. In addition, we show that solutions of \eqref{themodel} are bounded.
\begin{lemma}
Solutions of \eqref{themodel} are bounded with initial conditions in $\mathbb{R}^3_+=\{(u_1, u_2, v):u_1, u_2, v>0\}$.
\end{lemma}
\begin{proof}
It is easy to see that $u_1\le1$ and $u_2\le 1$ eventually by differential inequality. Let $w=\frac{a_{31}}{a_{13}}u_1+\frac{a_{32}}{a_{23}}u_2+v$. Hence, for time large enough,  we have
\begin{align*}
\frac{dw}{dt} &=  \frac{a_{31}}{a_{13}}(r_1u_1(1-u_1))+\frac{a_{32}}{a_{23}}(r_2u_2(1-u_2))-\frac{a_{31}}{a_{13}}a_{12}u_1u_2-\frac{a_{32}}{a_{23}}a_{21}u_1u_2-v\\
& \le  \frac{2a_{31}r_1}{a_{13}}u_1+ \frac{2a_{32}r_2}{a_{23}}u_2- \frac{a_{31}r_1}{a_{13}}u_1- \frac{a_{32}r_2}{a_{23}}u_2- v\\
&\le M-Dw,
\end{align*}
where $M=\left\{ \frac{2a_{31}r_1}{a_{13}}+ \frac{2a_{32}r_2}{a_{23}}\right\}$ and $D=\min\{r_1, r_2, 1\}$. This differential inequality implies that $w$ is bounded, which in turn implies that species $v$ is also bounded. The proof is complete.
\end{proof}

\begin{lemma} Let $\phi(t)=(u_1(t), u_2(t), v(t))$ be a solution of \eqref{themodel}. If $a_{31}+a_{32}< 1$ then $\lim_{t\to\infty}v(t)=0$.
\end{lemma}
\begin{proof}
Let $\varepsilon=1-a_{31}-a_{32}>0$, and take time $t$ large enough such that $\max\{u_1(t), u_2(t)\}<1+\varepsilon$.
Then \begin{align*}
\frac{\dot v(t)}{v(t)} &=  a_{31}u_1(t)+a_{32}u_2(t)-1 \\ 
& < (a_{31}+a_{32})(1+\varepsilon)-1=-\varepsilon^2<0,
\end{align*}
 which implies that $\lim_{t\to\infty}v(t)=0$ 
\end{proof}
Therefore, it is generic asked to have the following hypothesis,
\begin{description}
\item[(H1)] $a_{31}+a_{32}>1$.
\end{description}

With initial conditions in $H_{1v}$, $H_{2v}$ or $H_{12}$, we get three two-dimensional subsystems of \eqref{themodel} which are 
\begin{equation}\label{u1vsystem}
\left\{
\begin{aligned}
\dot u_1 &= r_1u_1(1-u_1)-a_{13}u_1v,\\
\dot v &= v(a_{31}u_1-1),
\end{aligned}
\right.
\end{equation}

\begin{equation}\label{u2vsystem}
\left\{
\begin{aligned}
\dot u_2 &= r_2u_2(1-u_2)-a_{23}u_2v,\\
\dot v &= v(a_{32}u_2-1),
\end{aligned}
\right.
\end{equation}

\begin{equation}\label{u1u2system}
\left\{
\begin{aligned}
\dot u_1 &= r_1u_1(1-u_1)-a_{12}u_1u_2,\\
\dot u_2 &= r_2u_2(1-u_2)-a_{21}u_1u_2,\\
\end{aligned}
\right.
\end{equation}
respectively. There are six equilibria of \eqref{u1vsystem}-\eqref{u1u2system} which are also boundary equilibria of \eqref{themodel},
\begin{align*}
E_0&=(0, 0, 0), & E_1&=(1, 0, 0),& E_2&=(0, 1, 0), \\
E_{1v}&=(\bar u_1, 0, \bar v_1),& E_{2v}&=(0, \bar u_2, \bar v_2) ,&E_{12}&=(\tilde u_1, \tilde u_2, 0), 
\end{align*} 
 where
\begin{align}
\bar u_1&=\frac{1}{a_{31}} , & \bar v_1&=\frac{r_1(1-\bar u_1)}{a_{13}}, \nonumber\\ 
\bar u_2&=\frac{1}{a_{32}}, &\bar v_2&=\frac{r_2(1-\bar u_2)}{a_{23}}, \nonumber\\
 \tilde u_1&=\frac{(r_1r_2-r_2a_{12})}{(r_1r_2-a_{12}a_{21})}, 
 & \tilde u_2&=\frac{(r_1r_2-r_1a_{21})}{(r_1r_2-a_{12}a_{21})}.\label{u1u2tilde}
\end{align} 
Here, we recall some well-known results of two-dimensional subsystems \eqref{u1vsystem}-\eqref{u1u2system}.
\begin{theorem} \label{2dresults} The following statements are valid.\label{2dim}
\begin{enumerate}[{\rm (a)}]
\item For subsystem \eqref{u1vsystem}, if $a_{31}< 1$ then $E_1$ is globally asymptotically stable on $H_{1v}$. Otherwise, if $a_{31}>1$, then $E_{1v}$ exists and is globally asymptotically stable on $H_{1v}$.
\item Similarly, for subsystem \eqref{u2vsystem}, if $a_{32}< 1$ then $E_2$ is globally asymptotically stable on $H_{2v}$. On the contrast, if $a_{32}>1$, then $E_{2v}$ exists and is globally asymptotically stable on $H_{2v}$.
\item For the competitive subsystem \eqref{u1u2system}, there are four generic cases.
\begin{enumerate}[{\rm (i)}] 
\item If $\frac{r_1}{a_{12}}>1$ and $\frac{r_2}{a_{21}}<1$, then $E_{12}$ does not exist and $E_1$ is GAS, which means globally asymptotically stable, on $H_{12}$.
\item If $\frac{r_1}{a_{12}}<1$ and $\frac{r_2}{a_{21}}>1$, then $E_{12}$ does not exist and $E_2$ is GAS on $H_{12}$.
\item If $\frac{r_1}{a_{12}}>1$ and $\frac{r_2}{a_{21}}>1$, then $E_{12}$ exists and is GAS on $H_{12}$.
\item If $\frac{r_1}{a_{12}}<1$ and $\frac{r_2}{a_{21}}<1$, then
$E_{12}$ exists and is saddle on $H_{12}$. The bistability phenomenon happens.
\end{enumerate}
\end{enumerate}
\end{theorem}

\subsection{Local stability of boundary equilibria in $\mathbb{R}^3_+$}

By straightforward computations, we have the Jacobian matrix of \eqref{themodel}, 

\begin{equation}\label{jacobian}
\begin{aligned}
J(u_1, u_2, v)&=&
\left[
\begin{array}{ccc}
 J_{11}&-a_{12}u_1 & -a_{13}u_1 \\
-a_{21}u_2 & J_{22}  & -a_{23}u_2 \\
a_{31}v &a_{32}v & -1+a_{31}u_1+a_{32}u_2
\end{array}
\right],
\end{aligned}
\end{equation}
where $J_{11}=r_1-2r_1u_1-a_{12}u_2-a_{13}v$ and $J_{22}=r_2-2r_2u_2-a_{21}u_1-a_{23}v$. Hence we obtain the local stability of $E_0$, $E_1$ and $E_2$ by evaluating  the Jacobian matrix, 
\begin{align}
J(E_0)=
\begin{bmatrix}
r_1 & 0 & 0\\
0 & r_2 & 0\\
0 & 0 & -1
\end{bmatrix}, \quad\label{JE1}
J(E_1)=
\begin{bmatrix}
-r_1 & -a_{12} & -a_{13}\\
0 & r_2-a_{21} & 0\\
0 & 0 & a_{31}-1
\end{bmatrix}, \\
\text{ and } \quad J(E_2)=
\begin{bmatrix}
r_1-a_{12} & 0 & 0\\
-a_{21} & -r_2 & -a_{23}\\
0 & 0 & a_{32}-1
\end{bmatrix}.\label{JE2}
\end{align}

\begin{lemma} \label{stabilityE012}The following statements are valid for the system \eqref{themodel}.
\begin{enumerate}[{\rm (a)}]
\item The trivial equilibrium $E_0$ is saddle in $\mathbb{R}^3_+$ with the stable manifold $v$-axis and the unstable manifolds $u_1$-$u_2$ plane.
\item The equilibrium $E_1$ is stable in $\mathbb{R}^3_+$ if $r_2<a_{21}$ and $a_{31}<1$.
\item The equilibrium $E_2$ is stable  in $\mathbb{R}^3_+$ if $r_1<a_{12}$ and $a_{32}<1$.
\end{enumerate}
\end{lemma}
%

It is more complicated for the equilibria $E_{1v}$, $E_{2v}$ and $E_{12}$. By evaluating \eqref{jacobian} at $E_{1v}$, we have 
\begin{align}
J(E_{1v})=
\begin{bmatrix}
-r_1\bar u_1& -a_{12}\bar u_1 &  -a_{13}\bar u_1\\
0 & r_2-a_{21}\bar u_1-a_{23}\bar v_1 & 0\\
a_{31}\bar v_1 & a_{32}\bar v_1 & 0
\end{bmatrix}.\label{JE1v}
\end{align}
It is obvious that there is one eigenvalue, $r_2-a_{21}\bar u_1-a_{23}\bar v_1$, and two negative eigenvalues from the $2\times 2$ submatrix obtained by eliminating the second row and second column of \eqref{JE1v}. Hence $E_{1v}$ is stable if $r_2-a_{21}\bar u_1-a_{23}\bar v_1<0$. Similar arguments can be obtained for the equilibrium $E_{2v}$ with its Jacobian matrix, 
\begin{align}
J(E_{2v})=
\begin{bmatrix}
r_1-a_{12}\bar u_2-a_{13}\bar v_2 & 0 &  0\\
-a_{21}\bar u_2 & -r_2\bar u_2 &  -a_{23}\bar u_2\\
a_{31}\bar v_2 & a_{32}\bar v_2 & 0
\end{bmatrix}.\label{JE2v}
\end{align}
For the equilibrium $E_{12}$, we also obtain the Jacobian matrix by evaluating \eqref{jacobian} at $E_{12}$,
\begin{align}
J(E_{12})=
\begin{bmatrix}
-r_1\tilde u_1& -a_{12}\tilde u_1 &  -a_{13}\tilde u_1\\
-a_{21}\tilde u_2 & -r_2\tilde u_2 & -a_{23}\tilde u_2\\
0 & 0 & a_{31}\tilde u_1+a_{32}\tilde u_2-1
\end{bmatrix}.\label{JE12}
\end{align}
By Theorem \ref{2dim}(c), $E_{12}$ exists if $\frac{r_1}{a_{12}}>1, \frac{r_2}{a_{21}}>1$ or $\frac{r_1}{a_{12}}<1, \frac{r_2}{a_{21}}<1$. However, the equilibrium $E_{12}$ exists and is stable only if $\frac{r_1}{a_{12}}>1$, $\frac{r_2}{a_{21}}>1$ and $a_{31}\tilde u_1+a_{32}\tilde u_2<1$.
The local stability of the equilibria $E_{1v}$, $E_{2v}$ and $E_{12}$ in $\mathbb{R}^3_+$ can be stated as follows.

\begin{lemma}\label{stabilityE12v}
The following statements are valid for the system \eqref{themodel}.
\begin{enumerate}[{\rm (a)}]
\item The equilibrium $E_{1v}$ exists if $a_{31}>1$, and it is stable in $\mathbb{R}^3_+$ if 
\begin{align}\label{E1vstable}
r_2-a_{21}\bar u_1-a_{23}\bar v_1<0.
\end{align}
\item The equilibrium $E_{2v}$ exists if $a_{32}>1$, and it is stable in $\mathbb{R}^3_+$ if 
\begin{align}\label{E2vstable}
r_1-a_{12}\bar u_2-a_{13}\bar v_2<0.
\end{align}
\item The equilibrium $E_{12}$ exists if $\frac{r_1}{a_{12}}>1$ and $\frac{r_2}{a_{21}}>1$ hold, or $\frac{r_1}{a_{12}}<1$ and $\frac{r_2}{a_{21}}<1$ hold. And it is stable in $\mathbb{R}^3_+$ if $\frac{r_1}{a_{12}}>1$, $\frac{r_2}{a_{21}}>1$ and 
\begin{align}
a_{31}\tilde u_1+a_{32}\tilde u_2<1.\label{E12stable}
\end{align}
\end{enumerate}
\end{lemma}

Finally, we state, without proof, the stable and unstable manifolds of all boundary equilibria under some specific conditions in Table \ref{manifoldtable}, since this follows from straightforward computations. This table will help us to show the uniform persistence of system \eqref{themodel}.

\begin{table}[htp]
\caption{The stable and unstable manifolds of all boundary equilibria of system \eqref{themodel} under some specific conditions.}
\begin{center}
{\footnotesize
\begin{tabular}{|c|c|l|l|}
\hline
 & conditions & the stable manifold & the unstable manifold \\
\hline
$E_0$ &&the $v$-axis & the $u_1$-$u_2$ plane\\
\hline
\hline
\multirow{3}{*}{$E_1$} &$r_2>a_{12}$ and $a_{31}>1$ &the $u_1$-axis&the tangent vector pointed to interior of $\mathbb{R}^3_+$\\
                                    \cline{2-4}
                                    &$r_2<a_{12}$ and $a_{31}>1$ &the $u_1$-$u_2$ plane&the tangent vector lied on  the $u_1$-$v$ plane\\
                                    \cline{2-4}
                                    &$r_2>a_{12}$ and $a_{31}<1$ &the $u_1$-$v$ plane&the tangent vector lied on the $u_1$-$u_2$ plane\\
\hline
\hline
\multirow{3}{*}{$E_2$} &$r_1>a_{21}$ and $a_{32}>1$ &the $u_2$-axis&the tangent vector pointed to interior of $\mathbb{R}^3_+$\\
                                    \cline{2-4}
                                    &$r_1<a_{21}$ and $a_{32}>1$ &the $u_1$-$u_2$ plane&the tangent vector lied on  the $u_2$-$v$ plane\\
                                    \cline{2-4}
                                    &$r_1>a_{21}$ and $a_{32}<1$ &the $u_2$-$v$ plane&the tangent vector lied on  the $u_1$-$u_2$ plane\\
\hline
\hline
$E_{1v}$& $r_2-a_{21}\bar u_1-a_{23}\bar v_1>0$ & the $u_1$-$v$ plane&the vector pointed to interior of $\mathbb{R}^3_+$\\
\hline
\hline
$E_{2v}$& $r_1-a_{12}\bar u_2-a_{13}\bar v_2>0$ & the $u_2$-$v$ plane&the vector pointed to interior of $\mathbb{R}^3_+$\\
\hline
\hline
\multirow{4}{*}{$E_{12}$}& $r_1>a_{21}$, $r_2>a_{12}$ and &\multirow{2}{*}{the $u_1$-$u_2$ plane}&\multirow{2}{*}{the vector pointed to interior of $\mathbb{R}^3_+$}\\
     & $a_{31}\tilde u_1+a_{32}\tilde u_2>1$ &&\\
\cline{2-4}
& $r_1<a_{21}$, $r_2<a_{12}$ and &the tangent vector lied&\multirow{2}{*}{the vector pointed to interior of $\mathbb{R}^3_+$} \\
     & $a_{31}\tilde u_1+a_{32}\tilde u_2>1$& on the $u_1$-$u_2$ plane&\\
\hline
\end{tabular}
}
\end{center}
\label{manifoldtable}
\end{table}%



\section{Existence of the Positive Equilibrium, Parameter Classification, and Some Extinctions in $\mathbb{R}^3_+$}
In this section, we first show the necessary and sufficient conditions for the existence of a positive equilibrium. Then, the classification of all parameters is introduced. Finally, some extinction results are proved analytically for parameter values that do not satisfy the necessary conditions.
\subsection{Existence of the positive equilibrium and classification of parameter space}
In this subsection, we try to find the positive equilibrium $E^*=(u_1^*, u_2^*, v^*)$ of system \eqref{themodel}. It is easy to see that we should solve the following linear systems, 
\begin{equation}\
\left\{
\begin{aligned}
 r_1(1-u_1)-a_{13}v-a_{12}u_2 =0,\\
 r_2(1-u_2)-a_{23}v-a_{21}u_1 =0, \\
 a_{31}u_1+a_{32}u_2 =1, 
\end{aligned}
\right.
\end{equation}
where the positive equilibrium is exactly the intersection of the following two straight lines in the interior of the first quadrant of the $u_1$-$u_2$ plane, 
\begin{equation}\label{twolines-old}
\left\{
\begin{aligned}
&L_1: ( r_1a_{23}-a_{21}a_{13})u_1+ (a_{12}a_{23}-r_2a_{13})u_2  =  r_1a_{23}- r_2a_{13},\\
&L_2:a_{31}u_1+a_{32}u_2 =1, 
\end{aligned}
\right.
\end{equation}
with $0<u_1^*<1$, $0<u_2^*<1$ and $v^*>0$.
Using the simplified notations, $C_1=r_1a_{23}$, $C_2=r_{2}a_{13}$, $C_3=a_{12}a_{23}$, and $C_4=a_{21}a_{13}$, these two linear equations \eqref{twolines-old} can be rewritten as the form,
\begin{equation}\label{twolines}
\left\{
\begin{aligned}
&L_1: ( C_1-C_4)u_1+ (C_3-C_2)u_2  =  C_1- C_2,\\
&L_2:a_{31}u_1+a_{32}u_2 =1. 
\end{aligned}
\right.
\end{equation}
By symmetry of $u_1$ and $u_2$ in \eqref{themodel}, without loss of generality, we always assume that 
\begin{description}
\item[(H2)] $C_1>C_2$.
\end{description}
Biologically, the assumption {\bf (H2)} means that the species $u_1$ is superior to $u_2$ in the ratio of birth rate over consuming rate $r_1/a_{13}>r_2/a_{23}$, that is, $u_1$ wins the apparent competition with $u_2$ \cite{Holt1977, Holt2017}. For the other generic case, $C_1=r_1a_{23}<C_2=r_2a_{13}$, we can obtain similar results by interchanging the roles of $u_1$ and $u_2$.

To show the existence of positive solutions of \eqref{twolines} and motivated by the well-known two competitive species $u_1$-$u_2$ models \eqref{u1u2system}, we consider the four generic categories: 
\begin{enumerate}[(I)]
\item $C_1<C_4$ \ and \ $C_3<C_2$,
\item $C_1>C_4$ \ and \ $C_3<C_2$,
\item $C_1<C_4$ \ and \ $C_3>C_2$,
\item $C_1>C_4$ \ and \ $C_3>C_2$.
\end{enumerate}
It is easy to see that the necessary conditions to have a positive equilibrium of system \eqref{themodel} is that lines $L_1$ and $L_2$ intersect in the interior of the first quadrant. Hence, we have further subcases depending on the intersection of lines $L_1$ and $L_2$ presented in Table \ref{classification}. For example, in the Category I : $C_1<C_4$ and $C_3<C_2$, the line $L_2$ has positive intercepts on $1/a_{31}$ and $1/a_{32}$ with $x$-axis and $y$-axis, respectively. However, the line $L_1$ only has negative intercepts on $\frac{C_1-C_2}{C_1-C_4}$ and $\frac{C_1-C_2}{C_3-C_2}$ with $x$-axis and $y$-axis, respectively. Hence, there is no positive equilibrium in Category I.
Other cases can be argued in similar ways, so we omit it. Therefore, we conclude that the positive equilibrium $E^*$ does not exist in the left column of Table \ref{classification}.

\begin{table}[htp]
\caption{Classification of parameters with assumption (A2) and the type of intercepts by lines $L_1$ and $L_2$.}
\begin{center}
\begin{tabular}{|c||cc||cc|}
\hline
\hline
\multicolumn{5}{|c|}{$C_1>C_2$}\\
\hline
\hline
\parbox{2.0cm}{\footnotesize{Category I\\ $C_1<C_4$ \\ $C_3<C_2$}} &
\parbox{2.1cm}{$\scriptstyle\frac{C_1-C_2}{C_1-C_4}<0<\frac{1}{a_{31}}$ \\ 
$\scriptstyle\frac{C_1-C_2}{C_3-C_2}<0<\frac{1}{a_{32}}$} &
\parbox{3.5cm}{
\begin{tikzpicture}[scale=0.4,domain=0:4]
    \draw[->] (-3, 0) -- (3.0, 0) node[right] {$u_1$};
    \draw[->] (0, -3.0) -- (0, 3.0) node[above] {$u_2$};
    \draw[-] (0.0, 1.8) -- (2.5, 0.0) 
                      node[below] {$\frac{1}{a_{31}}$}
                      node[pos=0.5, right]{$\scriptstyle L_2$};
    \draw[-] (2.5, 0.0) -- (0.0, 1.8) 
                      node[left] {$\frac{1}{a_{32}}$};
    \draw[-] (0, -2.5) -- (-2.7, 0.0) 
                      node[above] {$\scriptstyle\frac{C_1- C_2}{C_1-C_4}$}
                      node[pos=0.6, below]{$\scriptstyle L_1$};
    \draw[-] (-2.7, 0.0) -- (0, -2.5) 
                      node[right] {$\scriptstyle \frac{C_1- C_2}{C_3-C_2}$};
\end{tikzpicture}
} & & \\
\hline
\hline
\hline
\parbox{2.0cm}{\footnotesize{Category II \\ $C_1>C_4$ \\ $C_3<C_2$}} &
\parbox{2.1cm}{(i) \\ $\scriptstyle\frac{C_1-C_2}{C_1-C_4}>\frac{1}{a_{31}}$ \\ 
$\scriptstyle\frac{C_1-C_2}{C_3-C_2}<0<\frac{1}{a_{32}}$} &
\parbox{3.2cm}{
\begin{tikzpicture}[scale=0.4,domain=0:4]
    \draw[->] (-2, 0) -- (3.2, 0) node[right] {$u_1$};
    \draw[->] (0, -3.0) -- (0,3. 0) node[above] {$u_2$};
    \draw[-] (0.0, 1.8) -- (1.0, 0.0) 
                      node[below] {$\frac{1}{a_{31}}$}
                      node[pos=0.2, right]{$\scriptstyle L_2$};
    \draw[-] (1.0, 0.0) -- (0.0, 1.8) 
                      node[left] {$\frac{1}{a_{32}}$};
    \draw[-] (0, -1.5) -- (2.8, 0.0) 
                      node[above] {$\scriptstyle\frac{C_1- C_2}{C_1-C_4}$}
                      node[pos=0.9, below]{$\scriptstyle L_1$};
    \draw[-] (2.8, 0.0) -- (0, -1.5) 
                      node[left] {$\scriptstyle\frac{C_1- C_2}{C_3-C_2}$};
\end{tikzpicture}
} &
\parbox{2.1cm}{(ii) \\ $\scriptstyle\frac{C_1-C_2}{C_1-C_4}<\frac{1}{a_{31}}$ \\ 
$\scriptstyle\frac{C_1-C_2}{C_3-C_2}<0<\frac{1}{a_{32}}$} &
\parbox{3.2cm}{
\begin{tikzpicture}[scale=0.4,domain=0:4]
    \draw[->] (-0.5, 0) -- (4.0, 0) node[right] {$u_1$};
    \draw[->] (0, -3.0) -- (0,3. 0) node[above] {$u_2$};
    \draw[-] (0.0, 1.8) -- (3, 0.0) 
                      node[below] {$\scriptstyle\frac{1}{a_{31}}$}
                      node[pos=0.8, above]{$\scriptstyle L_2$};
    \draw[-] (3, 0.0) -- (0.0, 1.8) 
                      node[left] {$\scriptstyle\frac{1}{a_{32}}$};
    \draw[-] (0, -1.5) -- (1.5, 1.5) 
                      node[pos=0.55, above]{$\scriptstyle L_1$};
    \node [below] at (1.3, 0 ){$\scriptstyle\frac{C_1- C_2}{C_1-C_4}$};
    \draw[-] (1.5, 1.5) -- (0, -1.5) 
                      node[left] {$\scriptstyle\frac{C_1- C_2}{C_3-C_2}$};
\end{tikzpicture}
} \\
\hline
\hline
\hline
\parbox{2.0cm}{\footnotesize{Category III \\ $C_1<C_4$ \\ $C_3>C_2$}} &
\parbox{2.1cm}{(i) \\ $\scriptstyle\frac{C_1-C_2}{C_1-C_4}<0<\frac{1}{a_{31}}$ \\ 
$\scriptstyle\frac{C_1-C_2}{C_3-C_2}>\frac{1}{a_{32}}$} &
\parbox{3.2cm}{
\begin{tikzpicture}[scale=0.4,domain=0:4]
    \draw[->] (-3, 0) -- (3.0, 0) node[right] {$u_1$};
    \draw[->] (0, -1.0) -- (0,2.5) node[above] {$u_2$};
    \draw[-] (0.0, 0.7) -- (2.0, 0.0) 
                      node[below] {$\frac{1}{a_{31}}$}
                      node[pos=0.9, above]{$\scriptstyle L_2$};
    \draw[-] (2, 0.0) -- (0.0, 0.7) 
                      node[above right] {$\frac{1}{a_{32}}$};
    \draw[-] (0, 2.0) -- (-1.5, 0.0) 
                      node[below] {$\scriptstyle\frac{C_1- C_2}{C_1-C_4}$}
                      node[pos=0.6, left]{$\scriptstyle L_1$};
    \draw[-] (-1.5, 0.0) -- (0, 2.0) 
                      node[left] {$\scriptstyle\frac{C_1- C_2}{C_3-C_2}$};
\end{tikzpicture}

} &
\parbox{2.1cm}{(ii) \\ $\scriptstyle\frac{C_1-C_2}{C_1-C_4}<0<\frac{1}{a_{31}}$ \\ 
$\scriptstyle\frac{C_1-C_2}{C_3-C_2}<\frac{1}{a_{32}}$} &
\parbox{3.2cm}{
\begin{tikzpicture}[scale=0.4,domain=0:4]
    \draw[->] (-2.5, 0) -- (3.5, 0) node[right] {$u_1$};
    \draw[->] (0, -0.5) -- (0,2. 5) node[above] {$u_2$};
    \draw[-] (0.0, 2.0) -- (2.0, 0.0) 
                      node[below] {$\scriptstyle\frac{1}{a_{31}}$}
                      node[pos=0.5, right]{$\scriptstyle L_2$};
    \draw[-] (2, 0.0) -- (0.0, 2.0) 
                      node[right] {$\scriptstyle\frac{1}{a_{32}}$};
    \draw[-] (1.0, 1.5) -- (-1.5, 0.0) 
                      node[below] {$\scriptstyle\frac{C_1- C_2}{C_1-C_4}$}
                      node[pos=0.6, below]{$\scriptstyle L_1$};
    \draw[-] (-1.5, 0.0) -- (1.0, 1.5);
    \node [left] at (0, 1.2 ){$\scriptstyle\frac{C_1- C_2}{C_3-C_2}$};
                     
\end{tikzpicture}
} \\
\hline
\hline
\hline
\multirow{2}{*}{\parbox{2.0cm}{\footnotesize{Category IV \\ $C_1>C_4$ \\ $C_3>C_2$}}} &
\parbox{2.2cm}{(i) \\ $\scriptstyle\frac{C_1-C_2}{C_1-C_4}>\frac{1}{a_{31}}$ \\ 
$\scriptstyle\frac{C_1-C_2}{C_3-C_2}>\frac{1}{a_{32}}$} &
\parbox{3.2cm}{
\begin{tikzpicture}[scale=0.4]
    \def\x{-4.0}
    \def\y{-4.0}
    \draw[->] (\x-0.5, \y+0) -- (\x+3.3, \y+0) node[right] {$u_1$};
    \draw[->] (\x+0, \y-0.5) -- (\x+0,\y+3. 0) node[above] {$u_2$};
    \draw[-] (\x+0.0, \y+1.0) -- (\x+1.5, \y+0.0) 
                      node[below] {$\scriptstyle\frac{1}{a_{31}}$}
                      node[pos=0.7, above]{$\scriptstyle L_2$};
    \draw[-] (\x+1.5, \y+0.0) -- (\x+0.0, \y+1.0) 
                      node[left] {$\scriptstyle\frac{1}{a_{32}}$};
    \draw[-] (\x+0, \y+2.5) -- (\x+3.2, \y+0.0) 
                      node[below] {$\scriptstyle\frac{C_1- C_2}{C_1-C_4}$}
                      node[pos=0.6, above]{$\scriptstyle L_1$};
    \draw[-] (\x+3.2, \y+0.0) -- (\x+0, \y+2.5) 
                      node[left] {$\scriptstyle\frac{C_1- C_2}{C_3-C_2}$};
\end{tikzpicture}
} &
\parbox{2.2cm}{(iii) \\ $\scriptstyle\frac{C_1-C_2}{C_1-C_4}>\frac{1}{a_{31}}$ \\ 
$\scriptstyle\frac{C_1-C_2}{C_3-C_2}<\frac{1}{a_{32}}$} &
\parbox{3.2cm}{
\begin{tikzpicture}[scale=0.4,domain=0:4]
    \def\x{-2.0}
    \def\y{-2.0}
    \draw[->] (\x-0.5, \y+0) -- (\x+3.0, \y+0) node[right] {$u_1$};
    \draw[->] (\x+0, \y-0.5) -- (\x+0, \y+3. 0) node[above] {$u_2$};
    \draw[-] (\x+0.0, \y+1.0) -- (\x+2.5, \y+0.0) 
              node[below] {$\scriptstyle\frac{C_1- C_2}{C_1-C_4}$}
              node[pos=0.5, right]{$\scriptstyle L_1$};
    \draw[-] (\x+2.5, \y+0.0) -- (\x+0.0, \y+1.0) 
              node[left] {$\scriptstyle\frac{C_1- C_2}{C_3-C_2}$};
    \draw[-] (\x+0, \y+2.5) -- (\x+0.8, \y+0.0) 
              node[below] {$\scriptstyle\frac{1}{a_{31}}$}
                      node[pos=0.3, right]{$\scriptstyle L_2$};
    \draw[-] (\x+0.8, \y+0.0) -- (\x+0, \y+2.5) 
              node[left] {$\scriptstyle\frac{1}{a_{32}}$};
\end{tikzpicture}
}  \\
\cline{2-5}
&
\parbox{2.2cm}{(ii) \\ $\scriptstyle\frac{C_1-C_2}{C_1-C_4}<\frac{1}{a_{31}}$ \\ 
$\scriptstyle\frac{C_1-C_2}{C_3-C_2}<\frac{1}{a_{32}}$} &
\parbox{3.2cm}{
\begin{tikzpicture}[scale=0.4]
    \def\x{-4.0}
    \def\y{-4.0}
    \draw[->] (\x-0.5, \y+0) -- (\x+3.3, \y+0) node[right] {$u_1$};
    \draw[->] (\x+0, \y-0.5) -- (\x+0,\y+3. 0) node[above] {$u_2$};
    \draw[-] (\x+0.0, \y+1.0) -- (\x+1.5, \y+0.0) 
                      node[below] {$\scriptstyle\frac{C_1- C_2}{C_1-C_4}$}
                      node[pos=0.7, above]{$\scriptstyle L_1$};
    \draw[-] (\x+1.5, \y+0.0) -- (\x+0.0, \y+1.0) 
                      node[left] {$\scriptstyle\frac{C_1- C_2}{C_3-C_2}$};
    \draw[-] (\x+0, \y+2.5) -- (\x+3.2, \y+0.0) 
                      node[below] {$\scriptstyle\frac{1}{a_{31}}$}
                      node[pos=0.6, above]{$\scriptstyle L_2$};
    \draw[-] (\x+3.2, \y+0.0) -- (\x+0, \y+2.5) 
                      node[left] {$\scriptstyle\frac{1}{a_{32}}$};
\end{tikzpicture}
}  &
\parbox{2.2cm}{(iv) \\ $\scriptstyle\frac{C_1-C_2}{C_1-C_4}<\frac{1}{a_{31}}$ \\ 
$\scriptstyle\frac{C_1-C_2}{C_3-C_2}>\frac{1}{a_{32}}$} &
\parbox{3.2cm}{
\begin{tikzpicture}[scale=0.4,domain=0:4]
    \draw[->] (-0.5, 0) -- (3.3, 0) node[right] {$u_1$};
    \draw[->] (0, -0.5) -- (0,3. 3) node[above] {$u_2$};
    \draw[-] (0.0, 0.7) -- (2.8, 0.0) 
                      node[below] {$\scriptstyle\frac{1}{a_{31}}$}
                      node[pos=0.8, above]{$\scriptstyle L_2$};
    \draw[-] (2.8, 0.0) -- (0.0, 0.7) 
                      node[left] {$\scriptstyle\frac{1}{a_{32}}$};
    \draw[-] (0, 2.0) -- (1.0, 0.0) 
                      node[below] {$\scriptstyle\frac{C_1- C_2}{C_1-C_4}$}
                      node[pos=0.2, right]{$\scriptstyle L_1$};
    \draw[-] (1.0, 0.0) -- (0, 2.0) 
                      node[left] {$\scriptstyle\frac{C_1- C_2}{C_3-C_2}$};
\end{tikzpicture}
}  \\
\hline
\hline
\end{tabular}
\end{center}
\label{classification}
\end{table}%

However, in the right column of Table \ref{classification}, we have a positive intersection in the interior of the first quadrant, which is not sufficient to guarantee the existence of the positive equilibrium $E^*$. Hence, we write down the explicit form of $u^*_1$, $u^*_2$, and $v^*$, and $E^*$ exists if
\begin{align}
u_1^*&=\frac{a_{32}(C_1-C_2)-C_3+C_2}{a_{32}(C_1-C_4)+a_{31}(C_2-C_3)}<1, \label{u1l1}\\ 
u_2^*&=\frac{C_1-C_4-a_{31}(C_1-C_2)}{a_{32}(C_1-C_4)+a_{31}(C_2-C_3)}<1, \label{u2l1}\\
v^*&=\frac{1}{a_{13}}\left(r_1(1-u_1^*)-a_{12}u_2^*\right)=\frac{1}{a_{23}}\left(r_2(1-u_2^*)-a_{21}u_1^*\right)\nonumber\\
&=\frac{1}{a_{13}a_{23}}\frac{a_{31}C_2(C_1-C_3)+a_{32}C_1(C_2-C_4)-C_1C_2+C_3C_4}{a_{32}(C_1-C_4)+a_{31}(C_2-C_3)}>0. \label{vp}
\end{align}

\begin{lemma}[The necessary and sufficient condition for the existence of $E^*$]
The positive equilibrium $E^*$ exists if and only if the parameters satify assumptions {\bf(H1)}, {\bf(H2)}, \eqref{u1l1}, \eqref{u2l1}, \eqref{vp}, and one of conditions {\rm II-(ii)}, {\rm III-(ii)}, {\rm IV-(iii)} and {\rm IV-(iv)}.
\end{lemma}
For the equilibrium $E_{12}$, it can be also represented in terms of $C_1, \cdots, C_4$,  
\begin{align} \label{u1u2tildeform}
\tilde u_1=\frac{C_2(C_1-C_3)}{C_1C_2-C_3C_4}, \quad \tilde u_2=\frac{C_1(C_2-C_4)}{C_1C_2-C_3C_4},
\end{align}
and the third condition of the local stability \eqref{E12stable} can be also rewritten in the form,
\begin{align}
\frac{a_{31}C_2(C_1-C_3)+a_{32}C_1(C_2-C_4)}{C_1C_2-C_3C_4}<1. \label{E12stability}
\end{align}
By comparing \eqref{vp} and \eqref{E12stability}, if $C_1>C_3$ and $C_2>C_4$ then the stability of $E_{12}$ excludes the existence of $E^*$, and vice versa. We will see this fact in the next two sections.

We also use the previous classification to reformulate the conditions for local stability of boundary equilibria in $\mathbb{R}^3_+$, and we present the results in Table \ref{StabilityTable}. For example, $E_{1v}$ exists if $a_{31}>1$ and is stable if $\eqref{E1vstable}$ holds, by straightforward computations, 
\begin{align*}
r_2-a_{21}\bar u_1-a_{23}\bar v_1<0 \ \text{ if and only if }\ C_1-C_2>\frac{C_1-C_4}{a_{31}}.
\end{align*}
Hence, we obtain the stability of $E_{1v}$ in Table \ref{StabilityTable}. Similarly, $E_{2v}$ exists if $a_{32}>1$ and is stable if and only if 
$C_1-C_2<(C_3-C_2)/a_{32}$ holds. Finally, we summarize all local stability conditions for the boundary equilibria in terms of $C_1$-$C_4$ in Table \ref{StabilityTable}.

\begin{table}[htp]
\caption{The conditions of existence and global/local stability of each boundary equilibrium of system \eqref{themodel} on $H_{12}$/$\mathbb{R}^3$ are reformulated in terms of $C_1$-$C_4$.}
\begin{center}
\begin{tabular}{|c|c|c|c|}
\hline
& Existence & Global Stability in $H_{12}$ & Local Stability in $\mathbb{R}^3_+$\\
\hline
\hline
$E_1$ & & $C_1>C_3$ and $C_2<C_4$ & $a_{31}<1$ and $C_2<C_4$\\
\hline
$E_2$ & & $C_1<C_3$ and $C_2>C_4$ & $a_{32}<1$ and $C_1<C_3$ \\
\hline
\multirow{2}{*}{$E_{12}$} & ($C_1<C_3$, $C_2<C_4$) or  & \multirow{2}{*}{$C_1>C_3$ and $C_2>C_4$} & \multirow{2}{*}{$C_1>C_3$, \ $C_2>C_4$}  \multirow{2}{*}{\ and\ \eqref{E12stability}}\\
&($C_1>C_3$, $C_2>C_4$)&  &\\
\hline
\multirow{3}{*}{$E_{1v}$} & \multirow{3}{*}{$a_{31}>1$} & &
$C_1<C_4$ \quad or \\ 
&&&\multirow{2}{*}{($\frac{C_1-C_2}{C_1-C_4}>\frac{1}{a_{31}}$ if $C_1>C_4$)}\\
&&&\\
\hline
\multirow{2}{*}{$E_{2v}$} & \multirow{2}{*}{$a_{32}>1$} & &\multirow{2}{*}{$C_3>C_2$ \quad and \quad $\frac{C_1-C_2}{C_3-C_2}<\frac{1}{a_{32}}$}\\
&&&\\
\hline
\hline
\end{tabular}
\end{center}
\label{StabilityTable}
\end{table}%

\subsection{Extinction Results and the Corresponding Global Dynamics in Categories I, II-(i), III-(i) and IV-(i)}

We first investigate the four categories in the left column of Table \ref{classification}, Categories I, II-(i), III-(i), and IV-(i). 


\begin{proposition} \label{EcaseI} If assumptions {\rm \bf{(A1)}}, {\rm \bf{(A2)}} and Category {\rm (I)} hold, then $\lim_{t\to\infty} u_2(t)=0$. Furthermore, if $a_{31}< 1(a_{31}>1)$, then $E_1(E_{1v})$ is globally asymptotically stable.
\end{proposition}
\begin{proof}
Consider 
\begin{align*}
a_{13}\frac{\dot u_2}{u_2}-a_{23}\frac{\dot u_1}{u_1} &=C_2(1-u_2)-C_4u_1-C_1(1-u_1)+C_3u_2\\
&=-(C_1-C_2)-(C_2-C_3)u_2-(C_4-C_1)u_1\\
&\leq-(C_1-C_2)<0,
\end{align*}
which implies that $u_2^{a_{13}}(t)\le Ke^{-(C_1-C_2)t}$ for $t\ge 0$ and some positive constant $K$ by integrating both side of the above inequality. This is clear that $\lim_{t\to\infty} u_2(t)=0$. We conduct that \eqref{themodel} asymptotically approaches to system \eqref{u1vsystem} by Markus's limiting theorem. The remainder of the global results are followed by Theorem \ref{2dim}(a). We complete the proof.
\end{proof}

\begin{proposition}\label{EcaseII-i} If assumptions {\rm \bf{(A1)}}, {\rm \bf{(A2)}} and Category {\rm II-(i)} hold, that is,  
\begin{align} 
\frac{C_1-C_2}{C_1-C_4}>\frac{1}{a_{31}} \quad\text{ and }\quad
\frac{C_1-C_2}{C_3-C_2}<0<\frac{1}{a_{32}}
\end{align}
hold, then $\lim_{t\to\infty} u_2(t)=0$. Furthermore, if $a_{31}< 1(a_{31}>1)$, then $E_1(E_{1v})$ is globally asymptotically stable.
\end{proposition}
\begin{proof}
Take a positive number $\mu$ such that 
\begin{align*}
1<\mu<a_{31}\left(\frac{C_1-C_2}{C_1-C_4}\right),
\end{align*}
 and consider
\begin{align*}
a_{13}\frac{\dot u_2}{u_2}&-a_{23}\frac{\dot u_1}{u_1}-\left(\frac{C_1-C_2}{\mu}\right)\frac{\dot v}{v}\\
&=(\frac{1}{\mu}-1)({C_1-C_2})+\left(({C_1-C_4})-\frac{a_{31}}{\mu}({C_1-C_2})\right)u_1\\
&\quad+\left(({C_3-C_2})-\frac{a_{32}}{\mu}({C_1-C_2})\right)u_2\\
&\le(\frac{1}{\mu}-1)({C_1-C_2})<0
\end{align*}
which implies that $\lim_{t\to\infty} u_2(t)=0$ by similar arguments of the previous proposition. Hence, \eqref{themodel} asymptotically approaches system \eqref{u1vsystem} by Markus's limiting theorem. The remainder of the global results are also followed by Theorem \ref{2dim}(a), and the proof is completed.
\end{proof}


\begin{proposition}\label{EcaseIII-i} If assumptions {\rm \bf{(A1)}}, {\rm \bf{(A2)}} and case {\rm III-(i)} hold, that is, 
\begin{align} 
\frac{C_1-C_2}{C_1-C_4}<0<\frac{1}{a_{31}} \quad\text{ and }\quad
\frac{C_1-C_2}{C_3-C_2}>\frac{1}{a_{32}}
\end{align}
hold, then $\lim_{t\to\infty} u_2(t)=0$. Furthermore, if $a_{31}< 1(a_{31}>1)$, then $E_1(E_{1v})$ is globally asymptotically stable.
\end{proposition}
\begin{proof}
Take a positive number $\mu$ such that 
\begin{align*}
1<\mu<a_{32}\left(\frac{C_1-C_2}{C_3-C_2}\right),
\end{align*}
 and consider
\begin{align*}
a_{13}\frac{\dot u_2}{u_2}&-a_{23}\frac{\dot u_1}{u_1}-\left(\frac{C_1-C_2}{\mu}\right)\frac{\dot v}{v}\\
&=(\frac{1}{\mu}-1)({C_1-C_2})+\left(({C_1-C_4})-\frac{a_{31}}{\mu}({C_1-C_2})\right)u_1\\
&\quad+\left(({C_3-C_2})-\frac{a_{32}}{\mu}({C_1-C_2})\right)u_2\\
&\le(\frac{1}{\mu}-1)({C_1-C_2})<0
\end{align*}
which implies that $\lim_{t\to\infty} u_2(t)=0$. Similarly, \eqref{themodel} asymptotically approaches system \eqref{u1vsystem} by Markus's limiting theorem, and the remainder of the global results are followed by Theorem \ref{2dim}(a). The proof is completed.
\end{proof}

\begin{proposition}\label{EcaseIV-i} If assumptions {\rm \bf{(A1)}}, {\rm \bf{(A2)}} and case {\rm IV-(i)} hold, that is, 
\begin{align} 
\frac{C_1-C_2}{C_1-C_4}>\frac{1}{a_{31}} \quad\text{ and }\quad
\frac{C_1-C_2}{C_3-C_2}>\frac{1}{a_{32}} 
\end{align}
hold, then $\lim_{t\to\infty} u_2(t)=0$. Furthermore, if $a_{31}< 1(a_{31}>1)$, then $E_1(E_{1v})$ is globally asymptotically stable. 
\end{proposition}
\begin{proof}
Take a positive number $\mu$ such that 
\begin{align*}
1<\mu<\min\left\{a_{31}\left(\frac{C_1-C_2}{C_1-C_4}\right), a_{32}\left(\frac{C_1-C_2}{C_3-C_2}\right)\right\},
\end{align*}
 and consider
\begin{align*}
a_{13}\frac{\dot u_2}{u_2}&-a_{23}\frac{\dot u_1}{u_1}-\left(\frac{C_1-C_2}{\mu}\right)\frac{\dot v}{v}\\
&=(\frac{1}{\mu}-1)({C_1-C_2})+\left(({C_1-C_4})-\frac{a_{31}}{\mu}({C_1-C_2})\right)u_1\\
&\quad+\left(({C_3-C_2})-\frac{a_{32}}{\mu}({C_1-C_2})\right)u_2\\
&\le(\frac{1}{\mu}-1)({C_1-C_2})<0
\end{align*}
which implies that $\lim_{t\to\infty} u_2(t)=0$. Similarly, \eqref{themodel} asymptotically approaches system \eqref{u1vsystem} by Markus's limiting theorem, and the remainder of the global results are followed by Theorem \ref{2dim}(a). The proof is completed.
\end{proof}

\begin{remark}\label{negativecompetition}
It is easy to see that the proofs of these four global results of the categories I, II-(i), III-(i), and (IV)-(i), are similar, and the biological mechanisms are almost the same, where we call this the negative competition. We say that $u_1$ wins the negative competition with $u_2$ if $C_1 > C_3$, or equivalently, $r_1 > a_{12}$. Biologically, it means that the $r$-strategy of $u_1$ can overcome the negative effect from the species $u_2$. Similarly, we say that $u_2$ wins the negative competition with $u_1$ if $C_2>C_4$, or equivalently, $r_2>a_{21}$. 
Hence, there are four generic cases for $u_1$-$u_2$ on the negative competition, win-win if $C_1>C_3$ and $C_2>C_4$, or win-loss if $C_1>C_3$ and $C_2<C_4$, and so on. 
Obviously, the negative competition is independent of the predator $v$. Furthermore, this competition can also be applied to the two-dimensional $u_1$-$u_2$ subsystem and governs its global dynamics of \eqref{u1u2system}. The results of Theorem \ref{u1u2system}(c) can be rephrased as follows: 
\begin{enumerate}[(i)]
\item if win-loss, then $u_1$ and $u_2$ cannot coexist, $u_1$ exists permenantly, and $u_2$ dies out eventually;
\item if loss-win, then then $u_1$ and $u_2$ cannot coexist, $u_2$ exists permenantly, and $u_1$ dies out eventually;
\item if win-win, then $u_1$ and $u_2$ can globally stablly coexist;
\item if loss-loss, then $u_1$ and $u_2$ can coexist unstablly and the bistability phenomenon happens. 
\end{enumerate}

We imply that the negative competition is also the most important factor for the whole system \eqref{themodel}. The cases, I, II-(i), III-(i), and IV-(i) of \eqref{themodel} we just proved are discussed here, and others are presented in the final section.

In case I, with assumptions {\bf (A2)}, $C_1<C_4$, and $C_3<C_2$, it is obvious that $ u_1$ wins the negative competition of $u_1$-$u_2$. In case II-(i), it is easy to see that $u_1$ wins the negative competition, since $C_1>C_2>C_3$. On the other hand, for $a_{31}\le1$, it implies that $u_2$ losses since $\frac{C_1-C_2}{C_1-C_4}>\frac{1}{a_{31}}\ge1$. However, for $a_{31}>1$, it is possible to have the win-win state. But we will show that $u_2$ cannot win too much by this inequality $\frac{C_1-C_2}{C_1-C_4}>\frac{1}{a_{31}}$, which we call the apparent-tradeoff competition ratio. More details will be presented in the final section after we define the tradeoff competition.
\end{remark}
\section{Global Dynamics of \eqref{themodel} for Cases II-(ii) and III-(ii)}\label{II-(ii)}
For the following two sections, Section  \ref{II-(ii)} and  Section \ref{IV}, we consider the other five cases, II-(ii), III-(ii), IV-(ii), IV-(iii), and IV-(iv). A complete classification and its corresponding global dynamics based on parameters $C_1$, $C_2$, $C_3$, $C_4$ classifications are given. These cases, except for the case IV-(ii), satisfy the necessary condition of the existence of $E^*=(u^*_1, u^*_2, v^*)$ where the straight lines $L_1$ and $L_2$ intersect in the interior of the first octant. Hence, the positive equilibrium $E^*$ does not exist in the case IV-(ii), and for other four cases, $E^*$ exists if and only if $u_1^*<1$, $u_2^*<1$, and $v^*>0$. In particular, if 
\begin{align}\label{Hsu-condition}
\min\{a_{31}, a_{32}\}>1
\end{align}
then $E^*$ exists if and only if $v^*>0$ since the $u_1$- and $u_2$-coordinates of the intersection of $L_1$ and $L_2$ are less than $1/a_{31}$ and $1/a_{32}$, respectively. It is worth noting that the assumption \eqref{Hsu-condition} was assumed in \cite{Hsu1981}. Hence, we generalize the results of \cite{Hsu1981} for all positive $a_{31}$ and $a_{32}$. The dynamics of these five categories are more complex than we expected; therefore, in Section 4.1, we first establish the conditions for local stability of the positive equilibrium. Then these five categories are investigated in Sections 4 and 5.

\subsection{Local Stability of the Positive Equilibrium}
Now, we investigate the local stability of $E^*$ when it exists. 
By evaluating the Jacobian matrix,
\begin{equation*}
\begin{aligned}
J(u^*_1, u^*_2, v^*)&=&
\left[
\begin{array}{ccc}
 -r_1u^*_1&-a_{12}u^*_1 & -a_{13}u^*_1 \\
-a_{21}u^*_2 & -r_2u^*_2  & -a_{23}u^*_2 \\
a_{31}v^* &a_{32}v^* & 0
\end{array}
\right],
\end{aligned}
\end{equation*}
we can obtain the characteristic polynomial,
\begin{align}\label{E*charpoly}
P= \lambda^{3} + & \lambda^{2} \left(r_{1} u^*_1+r_{2} u^*_2\right) + \nonumber\\
&\lambda(a_{13} a_{31} u^*_1 v^* + a_{23} a_{32} u^*_2 v^*+
 r_{1} r_{2} u^*_1 u^*_2-a_{12} a_{21} u^*_1 u^*_2)+\nonumber \\
&(r_1a_{23}a_{32}+r_2a_{13} a_{31}-a_{12} a_{23} a_{31} - a_{13} a_{32}a_{21} )u^*_1 u^*_2 v^*.
\end{align}
Using the Routh–Hurwitz criterion, $E^*$ is stable if and only if
\begin{align}
& a_{13} a_{31} u^*_1 v^* + a_{23} a_{32} u^*_2 v^*+
 r_{1} r_{2} u^*_1 u^*_2-a_{12} a_{21} u^*_1 u^*_2>0,\label{RH-1}\\
&r_1a_{23}a_{32}+r_2a_{13} a_{31}-a_{12} a_{23} a_{31} - a_{13} a_{32}a_{21}>0, \qquad {\rm and }\label{RH-2}\\
& (r_{1} u^*_1+r_{2} u^*_2)(a_{13} a_{31} u^*_1 v^* + a_{23} a_{32} u^*_2 v^*+r_{1} r_{2} u^*_1 u^*_2-a_{12} a_{21} u^*_1 u^*_2)>\nonumber \\
& \qquad (r_1 a_{23} a_{32} +r_2a_{13} a_{31} -a_{12} a_{23} a_{31} - a_{13} a_{32} a_{21})u^*_1 u^*_2 v^*.\label{RH-3}
\end{align}
Obviously, if $a_{12}=a_{21}=0$, i.e., there is no direct competition between $u_1$ and $u_2$, then $E^*$ is always stable when it exists. In this case, $E^*$ is actually GAS when it exists \cite{Yang2025}. However, if 
\begin{align}\label{Existence-1}
r_1r_2>a_{12}a_{21},
\end{align} 
or equivalently $C_1C_2>C_3C_4$, then it is easy to see that \eqref{RH-1} and \eqref{RH-3} hold. Since \eqref{RH-1} and \eqref{RH-3} can be rewritten as the forms, 
\begin{align*}
a_{13} a_{31} u^*_1 v^* + a_{23} a_{32} u^*_2 v^*+
 (r_{1} r_{2}-a_{12} a_{21}) u^*_1 u^*_2> 0,
\end{align*}
and
\begin{align*}
r_{1} (u^*_1)^2&a_{13} a_{31}v^*+ r_{2}a_{23} a_{32} (u^*_2)^2 v^*+\\
&(r_1u_1^*+r_2u_2^*)(r_{1} r_{2} -a_{12} a_{21}) u^*_1 u^*_2 +
(a_{12} a_{23} a_{31} +a_{13} a_{32} a_{21})u^*_1 u^*_2 v^*>0,
\end{align*}
respectively. For the inequality \eqref{RH-2}, we also have an alternative form in terms of $C_1, \cdots, C_4$,
\begin{align}\label{RH-2-1}
 (C_1-C_4)a_{32}+(C_2-C_3)a_{31}>0.
\end{align}

\begin{lemma}\label{RH13} If the positive equilibrium $E^*$ exists, then the following statements are valid.
\begin{enumerate}[{\rm (i)}]
\item $E^*$ is stable if and only if \eqref{RH-1}, \eqref{RH-2} and \eqref{RH-3} hold.
\item If \eqref{Existence-1}, or equivalently $C_1C_2>C_3C_4$, and \eqref{RH-2-1} are true, then $E^*$ is stable.
\end{enumerate}
\end{lemma}

\subsection{The Category II-(ii)}

In this case, to clarify the dynamics of \eqref{themodel}  with the assumption {\bf (H2)}, we use the following logical strategy by separating the following all possible generic three cases of Category II with $C_1>C_4$ and $C_3<C_2$,
\begin{description}
\item[II-(ii) 1 :] $C_1>C_4>C_2>C_3$,
\item[II-(ii) 2 :] $C_1>C_2>C_4>C_3$,
\item[II-(ii) 3 :] $C_1>C_2>C_3>C_4$.
\end{description} 
In addition, by coupling two remaining parameters, $a_{31}$ and $a_{32}$, which are considered generically in the four ways, $\max\{a_{31}, a_{32}\}<1$, $a_{31}<1<a_{32}$, $a_{32}<1<a_{31}$, and $\min\{a_{31}, a_{32}\}>1$, we have all 12 subcases in Table \ref{tableII-ii} with labeled A1-A3, $\cdots$, D1-D3. In particular, 2 of the 12 cases labeled `contradiction' cannot occur. Since, in cells C1 and D1 of Table \ref{tableII-ii}, the inequality $\frac{C_1-C_2}{C_1-C_4}<\frac{1}{a_{31}}<1$ implies $C_2>C_4$, which is contradicts to $C_1>C_4>C_2>C_3$. To shorten our statements, the term ``the parameters are classified in the cell {\rm A1} of Table \ref{tableII-ii}'' means that the assumptions $\frac{C_1-C_2}{C_1-C_4}<\frac{1}{a_{31}}$, $\frac{C_1-C_2}{C_3-C_2}<\frac{1}{a_{32}}$, $\max\{a_{31}, a_{32}\}<1$, and $C_1>C_4>C_2>C_3$ hold for parameters, and others are similar. 

\begin{table}[htp]
\begin{threeparttable}
\caption{\textcolor{black}{All 12 subcases of the category II-(ii) under inequalities $0<\frac{C_1-C_2}{C_1-C_4}<\frac{1}{a_{31}}$ and $\frac{C_1-C_2}{C_3-C_2}<0<\frac{1}{a_{32}}$. }}
\begin{center}
\begin{tabular}{|c||c|c|c|c|}
\hline
\hline
&A &B&C&D\\
& \footnotesize{$\max\{a_{31}, a_{32}\}<1$} & \footnotesize{$a_{31}<1<a_{32}$} & \footnotesize{$a_{32}<1<a_{31}$} & \footnotesize{$1<\min\{a_{31}, a_{32}\}$} \\
\hline
\hline
1\ :\ \red{\scriptsize (W, L, W, W)} & \multirow{2}{*}{\black{$E_1$ is GAS}} &\multirow{2}{*}{\black{$E_1$ is GAS}} & \multirow{2}{*}{\blue{contradiction}} &\multirow{2}{*}{\blue{contradiction}} \\
{\scriptsize $C_1>C_4>C_2>C_3$} & & & &\\
\hline
2\ :\ \red{\scriptsize (W, W, W, W)}& \multirow{2}{*}{\black{$E_{12}/E^*$}\tnote{1} }&\multirow{2}{*}{\black{$E_{12}/E^*$}\tnote{1}} & \multirow{2}{*}{\black{$E_{12}/E^*$}\tnote{1}} & \black{\multirow{2}{*}{$E^*$\tnote{2}}}\\
{\scriptsize $C_1>C_2>C_4>C_3$} &  & &  & \\
\hline
3\ :\ \red{\scriptsize (W, W, W, W)} & \multirow{2}{*}{\black{$E_{12}/E^*$}\tnote{1}}& \multirow{2}{*}{\black{$E_{12}/E^*$}\tnote{1}} & \multirow{2}{*}{\black{$E_{12}/E^*$}\tnote{1}} &\black{\multirow{2}{*}{$E^*$\tnote{2}}}\\
 {\scriptsize $C_1>C_2>C_3>C_4$} &&&&\\
\hline
\hline
\end{tabular}
\begin{tablenotes}
\item[1] In these 6 cells, the equilibrium $E_{12}$ exists, and the notation $E_{12}/E^*$ means that either $E_{12}$ is stable or $E^*$ exists. Furthermore, $E^*$ is stable if it exists.
\item[2] In these 2 cells, the equilibrium $E^*$ exists and is also stable. 
\end{tablenotes}
\end{center}
\label{tableII-ii}
\end{threeparttable}
\end{table}%
Next, for all valid parameters in Table \ref{tableII-ii}, it is easy to see that $E_{12}$ exists in the second and third rows by checking Table \ref{StabilityTable}, and we check the existence of the positive equilibrium $E^*$ by verifying \eqref{u1l1}-\eqref{vp}. Then, clarify the dynamics of the cell by cell by checking the conditions in Table \ref{StabilityTable}. Senarioly, if we find that there is only one stable boundary equilibrium, we expect it to be GAS. If we find two stable boundary equilibria, we expect bistability. Finally, if we cannot find any stable boundary equilibrium, we expect that the positive equilibrium exists and may be stable.

Let us check the existence of $E^*$ in cells of Table \ref{tableII-ii}. Since the numerators and denominators of \eqref{u1l1}-\eqref{u2l1} are positive for cases (II)-(ii) and (IV)-(iv), hence \eqref{u1l1}-\eqref{vp} can be simplified as the forms,  
\begin{align} 
& a_{32}(C_2-C_4)>(1-a_{31})(C_2-C_3), \label{u1sc}\\
& a_{31}(C_1-C_3)>(1-a_{32})(C_1-C_4), \label{u2sc}\\
& a_{31}C_2(C_1-C_3)+a_{32}C_1(C_2-C_4)-C_1C_2+C_3C_4>0,\label{vpositive}
\end{align}
respectively, by straightforward computations. It is easy to see that $E^*$ does not exist for cells A1 and B1, since \eqref{u1sc} does not hold. However, we find a critical condition to verify the existence of the positive equilibrium $E^*$ for the second and third rows.

\begin{lemma}\label{II-ii-1}
Let the parameters be classified in {\rm Table \ref{tableII-ii}}. The following statements are valid.
\begin{enumerate}[{\rm (i)}]
\item $E^*$ does not exist for parameters classified in cells {\rm A1} and {\rm B1}.
\item In cells, {\rm A2}, {\rm B2}, {\rm C2}, {\rm A3}, {\rm B3} and {\rm C3} of {\rm Table \ref{tableII-ii}}, $E^*$ exists if and only if \eqref{vpositive} holds. Moreover, if $E^*$ exists, then it is also stable.
\item In cells {\rm D2} and {\rm D3}, $E^*$ exists and is stable.
\end{enumerate}
\end{lemma}
\begin{proof}
(ii) By \eqref{vpositive}, 
\begin{align*}
a_{31}C_2(C_1-C_3)&=C_1C_2-C_3C_4-a_{32}C_1C_2+a_{32}C_1C_4\\
&>C_1C_2-C_2C_4-a_{32}C_1C_2+a_{32}C_2C_4\\
&=C_2(1-a_{32})(C_1-C_4),
\end{align*}
which implies \eqref{u2sc}. By \eqref{vpositive} again,
\begin{align*}
&a_{32}C_1(C_2-C_4)\\
&=C_1(1-a_{31})(C_2-C_3)+C_3\left(C_1-C_4-a_{31}(C_1-C_2)\right)\\
&>C_1(1-a_{31})(C_2-C_3),
\end{align*}
which implies \eqref{u1sc}, since $0<\frac{C_1-C_2}{C_1-C_4}<\frac{1}{a_{31}}$. 

It is clear that \eqref{RH-2-1} is always true, and we can further get the inequality \eqref{Existence-1} by $C_1C_2= r_1a_{23}r_2a_{13}> C_3C_4 = a_{12}a_{23}a_{21}a_{13}$. Therefore, by Lemma \ref{RH13}, $E^*$ is stable. 

(iii) Consider 
\begin{align*}
a_{31}&C_2(C_1-C_3)+a_{32}C_1(C_2-C_4)-C_1C_2+C_3C_4\\
& =a_{31}C_2(C_1-C_3)+a_{32}C_1(C_2-C_4)-C_2(C_1-C_3)-C_3(C_2-C_4)\\
&=C_2(a_{31}-1)(C_1-C_3)+(a_{32}C_1-C_3)(C_2-C_4)>0,
\end{align*}
since $\min\{a_{31}, a_{32}\}>1$. And it is stable by Lemma \ref{RH13}.
We have completed the proof.
\end{proof}

By checking Table \ref{StabilityTable}, it is easy to check that $E_1$ is the only stable equilibrium in cells A1 and B1, and we will show that $E_1$ is GAS in these two cells. For $E_{12}$, stability holds if \eqref{E12stable} is true. However, it is easy to verify that the reverse of \eqref{E12stable},
\begin{align}\label{E12unstable}
a_{31}\tilde u_1+a_{32}\tilde u_2>1,
\end{align} 
is equivalent to \eqref{vpositive}. Hence, the existence of $E^*$ excludes the stability of $E_{12}$, and vice versa.
Here, we would like to show the global dynamics of all cases of Category II-(ii). One technical lemma is given first.
\begin{lemma}\label{maxless1}
Let $\phi(t)=\big(u_1(t), u_2(t), v(t))$ be a solution of \eqref{themodel}. If $C_4>C_2$ and $\max\{a_{31}, a_{32}\}<1$, then $\lim_{t\to\infty}u_2^{a_{13}}(t)v^\xi(t)=0$ for some positive number $\xi$.
\end{lemma}
\begin{proof}
We take $\xi$ such that $C_2<\xi\le\min\{\frac{C_2}{a_{32}}, \frac{C_4}{a_{31}}\}$ and consider 
\begin{align*}
a_{13}&\frac{\dot u_2}{u_2}+\xi\frac{\dot v}{v}\\
&=C_2-\xi+(\xi a_{31}-C_4)u_1+(\xi a_{32}-C_2)u_2-a_{13}a_{23}v\\
&\le C_2-\xi<0.
\end{align*}
This inequality implies that $\lim_{t\to\infty}u_2^{a_{13}}(t)v^\xi(t)=0$.
\end{proof}

\begin{proposition}\label{table4A1}
If the parameters are classified in the cell {\rm A1} of Table \ref{tableII-ii}, then $E_1$ is GAS.
\end{proposition}
\begin{proof}
Let $\phi(t; p)=(u_1(t), u_2(t), v(t); p)$ be a solution of \eqref{themodel} with initial condition $p\in\mathbb{R}^3_+$. We have $\lim_{t\to\infty}u_2^{a_{13}}(t)v^\xi(t)=0$ for some positive constant $\xi$, since parameters satisfy the assumptions of Lemma \ref{maxless1}. Similarly, we consider two cases: the first case is that $u_2(t_k)\to 0$ as $k\to\infty$ for some sequence $\{t_k\}$. For another case, there is a positive number $\epsilon$ such that $u_2(t)\ge\epsilon$ for all $t\ge 0$.

If the first case happens, a $q\in H_{1v}$ exists, such as $q\in\omega(p)$. And $\lim_{t\to\infty}\phi(t; q)=E_1$ by Theorem \ref{2dim}(a), hence $E_1\in\omega(p)$. This implies $\lim_{t\to\infty}\phi(t; p)=E_1$, since $E_1$ is stable on $\mathbb{R}^3_+$. 

On the other hand, if $u_2(t)\ge\epsilon$ for all $t\ge 0$ for some positive number $\epsilon$. The system \eqref{themodel} approaches to $H_{12}$ asymptotically. In this situation, we also have $\lim_{t\to\infty}\phi(t; p)=E_1$ by Theorem \ref{2dim}(c)(i). We have completed the proof.
\end{proof}

\begin{proposition}
If the parameters are classified in the cell {\rm B1} of Table \ref{tableII-ii}, then $E_1$ is GAS.
\end{proposition}
\begin{proof}
Let $\phi(t)=(u_1(t), u_2(t), v(t))$ be a solution of \eqref{themodel} and $\xi$ be a positive number such that $\frac{C_2}{C_1}<\xi<\frac{C_4}{C_1}<1<\frac{C_2}{C_3}$, and consider
\begin{align*}
a_{13}&\frac{\dot u_2}{u_2}-\xi a_{23}\frac{\dot u_1}{u_1}\\
&=C_2-\xi C_1+(\xi C_1-C_4)u_1+(\xi C_3-C_2)u_2-(1-\xi)a_{13}a_{23}v\\
&\le C_2-\xi C_1<0,
\end{align*}
which implies that $\lim_{t\to\infty} u_2(t)=0$. The system \eqref{themodel} approaches to $H_{1v}$ asymptotically. Hence, the boundary equilibrium $E_1$ is GAS by Theorem \ref{2dim}(a), so we complete the proof.
\end{proof}

By Lemma \ref{II-ii-1}(iii), $E^*$ exists and is stable for the parameters in the cells D2 and D3. Numerical observations indicate that $E^*$ is GAS in these two cells, but we cannot show this analytically. However, we can show that the system \eqref{themodel} is uniformly persistent.  
\begin{proposition}\label{ThoremII-ii-1}
System \eqref{themodel} is uniformly persistent for the parameters in cells {\rm D2} and {\rm D3} of Table \ref{tableII-ii}. 
\end{proposition}
\begin{proof}
We use the result in \cite{Freedman1994} to verify the uniform persistence of \eqref{themodel}. It suffices to show that the boundary equilibria in the first octant of the solution to \eqref{themodel} are isolated and acyclic. Consider all isolated boundary equilibria, $E_0$, $E_1$, $E_2$, $E_{1v}$, $E_{2v}$, and $E_{12}$. By checking Table \ref{manifoldtable}, it is impossible that the cycle between equilibria exists for any combination of the boundary equilibria. We complete the proof.
\end{proof}

\begin{theorem}\label{ThoremII-ii}
Let parameters be classified in the cells {\rm A2}, {\rm B2},  {\rm C2}, {\rm A3}, {\rm B3} and {\rm C3} of Table \ref{tableII-ii}. The statements are valid.
\begin{enumerate}[{\rm (i)}]
\item If \eqref{E12stable} holds, then $E^*$ does not exist
and $E_{12}$ is GAS.
\item If \eqref{E12unstable} holds, then $E_{12}$ is unstable with , $E^*$ exists and system \eqref{themodel} is uniformly persistent.
\end{enumerate}
\end{theorem}

\begin{lemma}\label{E12lemma}
Let $\phi(t)=\big(u_1(t), u_2(t), v(t); p\big)$ be a solution of \eqref{themodel} with initial condition $p \in\mathbb{R}^3_+$ and parameters classified in cells  {\rm A2}, {\rm B2}, {\rm C2}, {\rm A3}, {\rm B3}, and {\rm C3} of Table \ref{tableII-ii}. If $a_{31}\tilde u_1+a_{32}\tilde u_2<1$,
then $\lim_{t\to\infty}v(t)=0$. 
\end{lemma}
\begin{proof}
First of all, it is clear that the equilibrium $E_{12}$ is stable in $\mathbb{R}^3_+$ and GAS in the subspace $H_{12}$ for parameters in cells  {\rm A2}, {\rm B2}, {\rm C2}, {\rm A3}, {\rm B3}, and {\rm C3} of Table \ref{tableII-ii}. Please refer to Table \ref{StabilityTable}.
Define two numbers,
\[\alpha \equiv \frac{r_2a_{31}-a_{21}a_{32}}{r_1r_2-a_{12}a_{21}}\quad\text{ and }\quad \beta\equiv\frac{r_1a_{32}-a_{12}a_{31}}{r_1r_2-a_{12}a_{21}},\]
and consider, by the identities $r_1(1-\tilde u_1)-a_{12}\tilde u_2=0$ and $r_2(1-\tilde u_2)-a_{21}\tilde u_1=0$,  
\begin{align}
\frac{\dot v}{v}+&\alpha\frac{\dot u_1}{u_1}+\beta\frac{\dot u_2}{u_2}\nonumber \\
&\le -1+a_{31}u_1+a_{32}u_2+\alpha\left(r_1(\tilde u_1-u_1)+a_{12}(\tilde u_2-u_2)\right)+\nonumber \\
&\qquad\beta\left(r_2(\tilde u_2-u_2)+a_{21}(\tilde u_1-u_1)\right) \nonumber\\
&=-1+a_{31}\tilde u_1+a_{32}\tilde u_2+a_{31}(u_1-\tilde u_1)+a_{32}(u_2-\tilde u_2)+\nonumber\\
& \qquad\alpha\left(r_1(\tilde u_1-u_1)+a_{12}(\tilde u_2-u_2)\right)+\beta\left(r_2(\tilde u_2-u_2)+a_{21}(\tilde u_1-u_1)\right)\nonumber\\
& =-1+a_{31}\tilde u_1+a_{32}\tilde u_2<0.\label{vvanish}
\end{align}

Generically, we have four subcases:
(i) $\alpha<0$ and $\beta<0$, 
(ii) $\alpha<0$ and $\beta>0$, 
(iii) $\alpha>0$ and $\beta<0$,
(iv) $\alpha>0$ and $\beta>0$, 
to be considered. For the subcase (i) $\alpha<0$ and $\beta<0$, by similar previous arguments, we can obtain that $\lim_{t\to\infty}v(t)=0$ easily. Hence, we only show the most difficult subcase (iv) $\alpha>0$ and $\beta>0$, and the other two subcases (ii) and (iii) are similar and omitted.

By the inequality \eqref{vvanish}, we have $\lim_{t\to\infty}v(t)u_1^\alpha(t)u_2^\beta(t)=0$. Logically, we consider two cases: the first case is that $v(t_k)$ approaches zero as $k\to\infty$ for some infinite sequence $\{t_k\}$. For another case, there is a positive number $\epsilon$ such that $v(t)\ge\epsilon$ for all $t\ge 0$. If $v(t_k)$ approaches zero as $k\to\infty$ for some infinite sequence $\{t_k\}$, then there is a $q\in H_{12}$ such that $q\in\omega(p)$, the omega limit set of $p$. Consequently, the equliubrium $E_{12}$ also belong to $\omega(p)$ because of $\lim_{t\to\infty}\phi(t; q)=E_{12}$ and invariance of $\omega(p)$. This implies $\lim_{t\to\infty}\phi(t; p)=E_{12}$ since $E_{12}$ is stable in $\mathbb{R}^3_+$, and we can conclude that $\lim_{t\to\infty}v(t)=0$. 

On the other hand, let $v(t)\ge\epsilon$ for all $t\ge 0$ for some positive number $\epsilon$, and we show this cannot happen by contradiction. If 
$v(t)\ge\epsilon$ for all $t\ge 0$ for some positive number $\epsilon$, then we have $\lim_{t\to\infty}u_1^\alpha(t)u_2^\beta(t)=0$. 

We assume that $a_{31}<1$ where this inequality is true for cells A2, A3, B2, and B3. Similarly, we also have two cases to be considered; there is an infinite sequence $\{s_k\}$ such that $u_2(s_k)\to 0$ as $s_k\to\infty$, or there is a positive number $\delta$ such that $u_2(t)\ge\delta$ for all $t\ge0$. If there is an infinite sequence $\{s_k\}$ such that $u_2(s_k)\to 0$ as $s_k\to\infty$, then there is a $q_1\in H_{1v}$ such that $q_1\in\omega(p)$ and $\phi(t; q_1)\to E_1$ as $t\to \infty$ by  Theorem \ref{2dresults}. But this implies $E_1\in\omega(p)$ which contradicts to $v(t)\ge \epsilon>0$ for all $t\ge0$. On the other hand, if there is a positive number $\delta$ such that $u_2(t)\ge\delta$ for all $t\ge0$, then we have $\lim_{t\to\infty}u_1(t)=0$. Then the system \eqref{themodel} asymptotically approaches the system \eqref{u2vsystem} by the Markov Limit Theorem. On the subspace $H_{2v}$, if $a_{32}<1$ or $a_{32}>1$, then $\phi(t; p)$ approaches $E_2$ or $E_{2v}$, respectively. This implies $p$ is on the stable manifold of $E_2$ or $E_{2v}$, which is a contradiction.

 For cells C2, C3, we have $a_{32}<1$, and the arguments are simialr just considering $u_1$ instead of $u_2$. So we omit and complete the proof.
\end{proof}

\begin{proof}[Proof of Theorem \ref{ThoremII-ii}] ~\\
(i) Since $a_{31}\tilde u_1+a_{32}\tilde u_2<1$ is equvalent to the reverse of \eqref{vpositive}, $E^*$ does not exists. By Lemma \ref{E12lemma}, we have $\lim_{t\to\infty}v(t)=0$. Hence system \eqref{themodel} approaches to $H_{12}$ subsystem. We can conclude that $E_{12}$ is GAS by checking the conditions of Table \ref{StabilityTable}.

\noindent (ii) If $a_{31}\tilde u_1+a_{32}\tilde u_2>1$, then $E_{12}$ is unstable, $E^*$ exists and is stable. The proof of the uniform persistence for system \eqref{themodel} is similar to Proposition \ref{ThoremII-ii-1}, so we omit it. The proof is complete.
\end{proof}

\begin{remark}
In this category II-(ii), we show the global stability of $E_1$ for \eqref{themodel} with parameters in cells A1 and B1. For cells D2 and D3, the positive equilibrium exists and is stable, and we show that the system \eqref{themodel} is uniformly persistent. The GAS of $E^*$ cannot be established analytically, but we believe it holds. For the other six cells, the equilibrium $E_{12}$ exists and is stable if \eqref{E12stable} holds. Meanwhile, if $E_{12}$ is stable, then the positive equilibrium $E^*$ does not exist. Conversely, if \eqref{E12unstable} is true, then $E_{12}$ is unstable and $E^*$ exists stablely. In this case, we also believe that $E^*$ is GAS.

Biologically, all dynamics of these cells are compatible with those of $u_1$-$u_2$ subsystem \eqref{u1u2system} under negative competition. For cells A1 and B1, the win-loss state of the negative competition implies the global stability of $E_1$. For cells in the second and third rows, the win-win state implies $E_{12}$ or $E^*$ exists and is stable, where the key point is the inequality \eqref{E12unstable}. Without the top predator $v$, the constants $\tilde u_1$ and $\tilde u_2$ are the balance states of $u_1$ and $u_2$ when they negatively affect each other. However, when the top predator $v$ invades, the new balance states of $u_1$ and $u_2$ will be less than $\tilde u_1$ and $\tilde u_2$ under the negative effect of $v$. Hence the predator $v$ can get benifit from $u_1$ and $u_2$ is less than $a_{31}\tilde u_1+a_{32}\tilde u_1$. The holdiness of \eqref{E12unstable} means $v$ can survive by obtaining nutrients from $u_1$ and $u_2$.  
\end{remark}
\subsection{The Category III-(ii)}


For the category (III)-(ii), similarly, with assumption {\bf (H2)}, we further obtain finer classifications by separating the following three possible generic subcases with different values of $a_{31}$ and $a_{32}$ in Table \ref{tableIII-ii}.
\begin{description}
\item[III-(ii) 1 :] $C_3>C_4>C_1>C_2$,
\item[III-(ii) 2 :] $C_4>C_3>C_1>C_2$,
\item[III-(ii) 3 :] $C_4>C_1>C_3>C_2$,
\end{description} 


\begin{table}[htp]
\begin{threeparttable}
\caption{\textcolor{black}{All 12 subcases of the category III-(ii) under inequalities $\frac{C_1-C_2}{C_1-C_4}<0<\frac{1}{a_{31}}$ and $0<\frac{C_1-C_2}{C_3-C_2}<\frac{1}{a_{32}}$. }}
\begin{center}
\begin{tabular}{|c||c|c|c|c|}
\hline
\hline
&A &B&C&D\\
& \footnotesize{$\max\{a_{31}, a_{32}\}<1$} & \footnotesize{$a_{31}<1<a_{32}$} & \footnotesize{$a_{32}<1<a_{31}$} & \footnotesize{$1<\min\{a_{31}, a_{32}\}$} \\
\hline
\hline
1\ :\  \red{\scriptsize (L, L, L, L)}& Bistability\tnote{1} &Bistability\tnote{1}& Bistability\tnote{1} &Bistability\tnote{1} \\
{\scriptsize $C_3>C_4>C_1>C_2$} & $E_1/E_2$ & $E_{1}/E_{2v}$ & $E_{1v}/E_{2}$ &$E_{1v}/E_{2v}$\\
\hline
2\ :\ \red{\scriptsize (L, L, L, L)}& Bistability\tnote{1} &Bistability\tnote{1}& Bistability\tnote{1}& Bistability\tnote{1}\\
{\scriptsize $C_4>C_3>C_1>C_2$} & $E_1/E_2$ & $E_{1}/E_{2v}$ & $E_{1v}/E_{2}$ &$E_{1v}/E_{2v}$\\
\hline
3\ :\  \red{\scriptsize (W, L, L, L)}&\multirow{2}{*}{\black{$E_1$ is GAS}}&\multirow{2}{*}{\blue{contradiction}} & \multirow{2}{*}{ \black{$E_{1v}$ is GAS} }&\multirow{2}{*}{\blue{contradiction}}\\
 {\scriptsize $C_4>C_1>C_3>C_2$} &&&&\\
\hline
\hline
\end{tabular}
\begin{tablenotes}
\item[1] Please see the details in Remark \ref{remakrIII-ii}.
\end{tablenotes}
\end{center}
\label{tableIII-ii}
\end{threeparttable}
\end{table}%

Looking the inequalities \eqref{u1l1}-\eqref{vp} for the existence of the positive equilibrium $E^*$, since the numerators and denominators of \eqref{u1l1} and \eqref{u2l1} are negative for cases (III)-(ii) and (IV)-(iii), the conditions $u^*_1<1$, $u^*_2<1$ and $v^*>0$ are equivalent to 
\begin{align} 
&a_{32}(C_4-C_2)>(1-a_{31})(C_3-C_2)\label{u1pIII-ii},\\
& a_{31}(C_3-C_1)>(1-a_{32})(C_4-C_1)\label{u2pIII-ii},\\
& a_{31}C_2(C_3-C_1)+a_{32}C_1(C_4-C_2)+C_1C_2-C_3C_4>0.\label{vpIII-ii}
\end{align}
You can see that $E^*$ does not exist for cells A3 and C3 because \eqref{u2pIII-ii} is invalid. We have the following global result for these two cells.
\begin{proposition}
If $C_4>C_1>C_3>C_2$, then $\lim_{t\to\infty}u_2(t)=0$. Moreover, if $a_{31}<1$ then $E_{1}$ is GAS. Otherwise, if $a_{31}>1$ then $E_{1v}$ exists and is GAS.
\end{proposition}
\begin{proof}
Without loss of generality, we may assume $u_2(t)\le 1$ for $t$ is large enough. Consider 
\begin{align*}
a_{13}&\frac{\dot u_2}{u_2}-a_{23}\frac{\dot u_1}{u_1}\\
&= C_2-C_1-(C_4-C_1)u_1+(C_3-C_2)u_2\\
&\le C_2-C_1+C_3-C_2=C_3-C_1<0,
\end{align*}
which implies that $\lim_{t\to\infty}u_2(t)=0$. By Markus's limiting theorem, system \eqref{themodel} approaches \eqref{u1vsystem} asymptotically. Therefore, the results of this proposition are just a consequence of  Theorem \ref{2dim}(a). We have completed the proof.
\end{proof}

For the remaining eight cells of Table \ref{tableIII-ii}, the equilibrium $E_{12}$ exists, though it is unstable in $\mathbb{R}^3_+$ by Table \ref{StabilityTable}. 
Similar to Category II-(ii), we can establish results similar to Lemma \ref{II-ii-1}.
\begin{lemma}\label{lemmaIII-ii}
In cellls, {\rm A1}, {\rm B1}, {\rm C1}, {\rm D1}, {\rm A2}, {\rm B2}, {\rm C2} and {\rm D2} of Table \ref{tableIII-ii}, the positive equilibrium $E^*$ exists if and only if \eqref{vpIII-ii} holds. Moreover, $E^*$ does not exist for cells {\rm A1} and {\rm A2}.
\end{lemma}
\begin{proof}
It is sufficient to show that \eqref{vpIII-ii} implies \eqref{u1pIII-ii} and \eqref{u2pIII-ii}. Consider 
\begin{align*}
a_{32}C_1(C_4-C_2)&>-a_{31}C_2C_3+a_{31}C_1C_2-C_1C_2+C_3C_4\\
&>-a_{31}C_1C_3+a_{31}C_1C_2-C_1C_2+C_1C_3\\
&= (1-a_{32})C_1(C_3-C_2),
\end{align*}
which implies \eqref{u1pIII-ii} holds.
Consider 
\begin{align*}
a_{31}C_2&(C_3-C_1)>-a_{32}C_1(C_4-C_2)-C_1C_2+C_3C_4\\
&=(a_{32}C_1C_2-a_{32}C_2C_4+C_2C_4-C_1C_2)\\
&\qquad +(-a_{32}C_1C_4+a_{32}C_2C_4+C_3C_4-C_2C_4)\\
&=(1-a_{32})C_2(C_4-C_1)+C_4\big(C_3-C_2-a_{32}(C_1-C_2)\big)\\
&>(1-a_{32})C_2(C_4-C_1), 
\end{align*}
because of $\frac{C_1-C_2}{C_3-C_2}<\frac{1}{a_{32}}$. This implies that \eqref{u2pIII-ii} holds. 

To show the non-existence of $E^*$ for cells A1 and A2, rewrite the reverse of \eqref{vpIII-ii} as the form,
\[ (a_{31}-1)C_2(C_3-C_1)+(a_{32}C_1-C_3)(C_4-C_2)<0.\]
This inequality is true since $a_{31}<1$ and $a_{32}C_1-C_3<0$ by checking 
\begin{align*}
\frac{C_1-C_2}{C_3-C_2}<\frac{1}{a_{32}} \ \text{ which implies } \ a_{32}C_1-C_3<(a_{32}-1)C_2<0.
\end{align*}
So we complete the proof.
\end{proof}

However, $E^*$ is unstable even if it exists for the other 6 cells of the first and second rows of Table \ref{tableIII-ii} because of $C_1<C_4$ and $C_2<C_3$ by \eqref{RH-2-1}. We find that, by checking the stability conditions of Table \ref{StabilityTable}, $E_1/E_{1v}$ and $E_2/E_{2v}$ are stable simultaneously, which depend on the values of $a_{31}$ and $a_{32}$ less than or greater than 1, respectively. Bistability is expected in these eight cells.
\begin{lemma}\label{bistabilitylemma}
Let parameters be classified in the first and second rows of Table \ref{tableIII-ii}. Then the following statements are valid for system \eqref{themodel}.
\begin{enumerate}[{\rm (i)}]
\item If \eqref{E12stable} holds, then $E_{12}$ exists and $E^*$ does not exist. Moreover, $E_{12}$ is a saddle with a two-dimensional stable manifold, where one stable eigenvector lies on $H_{12}$ and another is transversal to $H_{12}$.
\item If \eqref{E12unstable} holds, then $ E _ {12} $ exists and is a saddle with a one-dimensional stable manifold. Moreover, $E^*$ exists and is a saddle with a two-dimensional stable manifold.
\end{enumerate}
\end{lemma}
\begin{proof}
It is easy to see that $E_{12}$ exists, and by observing the Jacobian matrix of $E_{12}$ \eqref{JE12}, $E_{12}$ is a saddle with a two-dimensional stable manifold if \eqref{E12stable} holds. Or $E_{12}$ is a saddle with a one-dimensional stable manifold if \eqref{E12unstable} holds.

If \eqref{E12unstable} holds, we determine the number of negative eigenvalues of $E^*$ using the Routh-Hurwitz criterion. Let $P(\lambda)=a_3\lambda^3+a_2\lambda^2+a_1\lambda+a_0$ be the characteristic polynomial for $E^*$, and, by \eqref{E*charpoly}, we obtain \begin{align*}
a_3&=1, \\
a_2&=r_{1} u^*_1+r_{2} u^*_2,\\
a_1&=a_{13} a_{31} u^*_1 v^* + a_{23} a_{32} u^*_2 v^*+
 r_{1} r_{2} u^*_1 u^*_2-a_{12} a_{21} u^*_1 u^*_2,\\
a_0&=(r_1a_{23}a_{32}+r_2a_{13} a_{31}-a_{12} a_{23} a_{31} - a_{13} a_{32}a_{21} )u^*_1 u^*_2 v^*\\
&=\big((C_1-C_4)a_{32}+(C_2-C_3)a_{31}\big)u^*_1u^*_2v^*<0.
\end{align*}
By computing the Routh table, we obtain
\begin{align*}
b_1&=\frac{a_2a_1-a_3a_0}{a_2}=- a_{12} a_{21} r_{1} u_{1}^{2} u_{2} - a_{12} a_{21} r_{2} u_{1} u_{2}^{2} +a_{12} a_{23} a_{31} u_{1} u_{2} v +\\
&  a_{13} a_{21} a_{32} u_{1} u_{2} v +a_{13} a_{31} r_{1} u_{1}^{2} v + a_{23} a_{32} r_{2} u_{2}^{2} v + r_{1}^{2} r_{2} u_{1}^{2} u_{2} + r_{1} r_{2}^{2} u_{1} u_{2}^{2},\\
b_2&=b_3=\cdots=0,\\
c_1&=\frac{b_1a_0}{b_1}=a_0<0,\\
c_2&=c_3=\cdots=0.
 \end{align*}
Therefore, the first column of the Routh table is $a_3$, $a_2$, $b_1$, and $c_1$. We can check $c_1<0$ by Lemma \ref{IV-iiiE*unstable}. Hence, there is only one sign change whenever $b_1>0$ or $b_1<0$. Hence, the dimension of the stable manifold of $E^*$ is two. We complete the proof.
\end{proof}

Therefore, we can easily obtain the following bistability results by Lemma \ref{lemmaIII-ii} and \ref{bistabilitylemma}, omitting the proofs.
\begin{proposition}[Bistability] \label{III-iib} For parameters classified in cells, {\rm A1}, {\rm A2}, {\rm B1}, {\rm B2}, {\rm C1}, {\rm C2}, {\rm D1}, and {\rm D2}, the following statements are valid.
\begin{enumerate}[{\rm (i)}]
\item For cells {\rm A1} and {\rm A2}, $E^*$ doest not exist and $E_{12}$ is a saddle with two-dimensionan stable manifold. The absorbing regions of the stable nodes, $E_1$ and $E_2$, are separated by the stable manifold of $E_{12}$.
\item For cells {\rm B1}-{\rm D1} and {\rm B2}-{\rm D2}, 
\begin{enumerate}[{\rm (a)}]
\item If \eqref{E12stable} holds, then $E^*$ does not exist, and $E_{12}$ is a saddle with a two-dimensional stable manifold. The absorbing region of each stable node is separated by the stable manifold of $E_{12}$;
\item if \eqref{E12unstable} holds, then 
$E^*$ exists and is a saddle with a two-dimensional stable manifold. The absorbing region of each stable node is separated by the stable manifold of $E^*$.
\end{enumerate}
\end{enumerate}
\end{proposition}

\begin{remark}\label{remakrIII-ii}
In this category III-(ii), we show the global stability of $E_1$ or $E_{1v}$ for \eqref{themodel} with parameters in cell A3 or C3, respectively. For the other eight cells of Table \ref{tableIII-ii}, the equilibrium $E_{12}$ exists and is a saddle. The positive equilibrium $E^*$ exists if and only if \eqref{E12unstable} holds, but it is also a saddle when it exists. The dynamics of these eight cells are clear. If \eqref{E12stable} holds, then $E^*$ does not exist and $E_{12}$ is a saddle with a two-dimensional stable manifold that separates the absorbing regions of the two stable nodes shown on each cell. On the other hand, if \eqref{E12unstable} holds, $E_{12}$ is a saddle with a one-dimensional stable manifold, and $E^*$ exists and is a saddle with a two-dimensional stable manifold that separates the absorbing regions of the two stable nodes shown on each cell. All dynamics of these eight cells are compatible with those of \eqref{u1u2system} in the loss-loss state under negative competition. We conclude that the only difference in dynamics between \eqref{u1u2system} and \eqref{themodel} is that the bistability $E_1/E_2$ is replaced by $E_{1v}/E_2$, $E_1/E_{2v}$, or $E_{1v}/E_{2v}$, depending on the values of $a_{31}$ and $a_{32}$, as separated by the two-dimensional manifold of $E_{12}$ or $E_*$.
\end{remark}

\section{Global Dynamics of the Categories IV-(ii), IV-(iii), and IV-(iv)}\label{IV}

In Category IV, the dynamics of \eqref{themodel} is more complicated compared to the previous categories, since, with the assumption {\bf (H2)}, $C_1>C_4$ and $C_3>C_2$, we should further separate the following all possible generic five subcases,
\begin{description}
\item[IV 1 :] $C_3>C_1>C_2>C_4$,
\item[IV 2 :] $C_3>C_1>C_4>C_2$,
\item[IV 3 :] $C_1>C_3>C_2>C_4$,
\item[IV 4 :] $C_1>C_3>C_4>C_2$,
\item[IV 5 :] $C_1>C_4>C_3>C_2$.
\end{description} 
And this classification applies to all remaining cases of category IV, including IV-(ii), IV-(iii), and IV-(iv), and we will also present all possible dynamics in Tables \ref{tableIV-ii}-\ref{tableIV-iv}, respectively. 

From the previous discussions, we know that the positive equilibrium $E^*$ does not exist for the case IV-(ii), but it could exist for the cases IV-(iii) and IV-(iv), since the necessary conditions for the existence of $E^*$ hold in these two cases. We first establish a result that applies to all cases of IV categories.
\begin{lemma}\label{C1dominate}
If $C_1>\max\{C_3, C_4\}$ and $\min\{C_3, C_4\}>C_2$, then $\lim_{t\to\infty}u_{2}(t)=0$.
\end{lemma}
\begin{proof}
Take a positive number $\xi$ such that $\frac{C_2}{C_1}<\xi<\min\{\frac{C_4}{C_1}, \frac{C_2}{C_3}\}<1$, and consider
\begin{align*}
 a_{13}&\frac{\dot u_2}{u_2}-\xi a_{23}\frac{\dot u_1}{u_1}\\
&= C_2-\xi C_1-(C_4-\xi C_1)u_1-(C_2-\xi C_3)u_2-(1-\xi)a_{13}a_{23}v\\
&\le C_2-\xi C_1<0,
\end{align*}
which implies that $\lim_{t\to\infty}u_2(t)=0$. We have completed the proof.
\end{proof}
This lemma applies to the above classifications, IV 4 and IV 5; hence, $u_2(t)$ always approaches 0 in rows 4 and 5 of the classification tables of category IV in the subsections that follow.
\subsection{The Category IV-(ii)}
Similarly, by coupling $a_{31}$ and $a_{32}$, we have all 20 subcases in Table \ref{tableIV-ii}, and 10 of the 20 subcases labeled `contradiction' cannot occur. 
\begin{table}[htp]
\caption{\textcolor{black}{Dynamics of All 20 subcases of the category IV-(ii) under inequalities $0<\frac{C_1-C_2}{C_1-C_4}<\frac{1}{a_{31}}$ and $0<\frac{C_1-C_2}{C_3-C_2}<\frac{1}{a_{32}}$. }}
\begin{center}
\begin{tabular}{|c||c|c|c|c|}
\hline
\hline
&A &B&C&D\\
& {\scriptsize $\max\{a_{31}, a_{32}\}<1$} &  $a_{31}<1<a_{32}$& $a_{32}<1<a_{31}$ & {\scriptsize $\min\{a_{31}, a_{32}\}>1$}  \\
\hline
\hline
1\ :\ \red{\scriptsize (L, W, W, L)}&  \multirow{2}{*}{\black{$E_{2}$ is GAS} } &\multirow{2}{*}{\black{$E_{2v}$ is GAS} }&\multirow{2}{*}{\black{$E_{2}$ is GAS} }& \multirow{2}{*}{\black{$E_{2v}$ is GAS} }\\
{\scriptsize $C_3>C_1>C_2>C_4$} &&&&\\
\hline
2\ :\  \red{\scriptsize (L, L, W, L)}&\black{Bistability} &\black{Bistability} &\multirow{2}{*}{\blue{contradiction} } &\multirow{2}{*}{\blue{contradiction}}\\
 {\scriptsize $C_3>C_1>C_4>C_2$} &$E_1/E_2$& $E_1/E_{2v}$  &&\\
\hline
3\ :\ \red{\scriptsize (W, W, W, L)} & \multirow{2}{*}{\black{$E_{12}$ is GAS} } &\multirow{2}{*}{\blue{contradiction} }& \multirow{2}{*}{\black{$E_{12}$ is GAS}} &\multirow{2}{*}{\blue{contradiction}} \\
{\scriptsize $C_1>C_3>C_2>C_4$} &&&&\\
\hline
4\ :\ \red{\scriptsize (W, L, W, L)} &\multirow{2}{*}{\black{$E_1$ is GAS} }& \multirow{2}{*}{\blue{contradiction} } &\multirow{2}{*}{\blue{contradiction} }&  \multirow{2}{*}{\blue{contradiction}} \\
 {\scriptsize $C_1>C_3>C_4>C_2$} &&&&\\
\hline
5\ :\ \red{\scriptsize (W, L, W, L)} &\multirow{2}{*}{\black{$E_1$ is GAS} }& \multirow{2}{*}{\blue{contradiction} } &\multirow{2}{*}{\blue{contradiction} }&\multirow{2}{*}{\blue{contradiction}} \\
 {\scriptsize $C_1>C_4>C_3>C_2$} &&&&\\
\hline
\hline
\end{tabular}
\end{center}
\label{tableIV-ii}
\end{table}%
We first show one technical lemma that can be applied to all subcases of case IV-(ii).

\begin{lemma}\label{EcaseIV-ii}Let $\phi(t)=(u_1(t), u_2(t), v(t))$ be a solution of \eqref{themodel}, and assumptions {\rm \bf{(A1)}}, {\rm \bf{(A2)}} and case {\rm IV-(ii)}, 
\begin{align*} 
\frac{C_1-C_2}{C_1-C_4}<\frac{1}{a_{31}} \quad\text{ and }\quad
\frac{C_1-C_2}{C_3-C_2}<\frac{1}{a_{32}},
\end{align*}
hold. Then the following statements are valid. 
\begin{enumerate}[{\rm (i)}]
\item $\lim_{t\to\infty}u_1^{a_{23}}(t)v^\xi(t)=0$ for some positive number $\xi$.
\item Inequality \eqref{E12stable} holds.
\end{enumerate}
\end{lemma}
\begin{proof}
(i) Take a positive number $\mu$ such that 
\begin{align*}
\max\left\{a_{31}\left(\frac{C_1-C_2}{C_1-C_4}\right), a_{32}\left(\frac{C_1-C_2}{C_3-C_2}\right)\right\}<\mu<1,
\end{align*}
 and consider
\begin{align*}
a_{23}\frac{\dot u_1}{u_1}&-a_{13}\frac{\dot u_2}{u_2}+\left(\frac{C_1-C_2}{\mu}\right)\frac{\dot v}{v}\\
&=(1-\frac{1}{\mu})({C_1-C_2})+\left(\frac{a_{31}}{\mu}({C_1-C_2})-({C_1-C_4})\right)u_1\\
&\quad+\left(\frac{a_{32}}{\mu}({C_1-C_2})-({C_3-C_2})\right)u_2\\
&\le(1-\frac{1}{\mu})({C_1-C_2})<0, 
\end{align*}
which implies that $\lim_{t\to\infty}u_1^{a_{23}}(t)v^\xi(t)=0$ with $\xi=(C_1-C_2)/\mu$.

(ii) With assumptions and the forms of $\tilde u_1$, $\tilde u_2$ \eqref{u1u2tildeform}, we have
\begin{align*}
a_{31}\tilde u_1+a_{32}\tilde u_2<\frac{C_1-C_4}{C_1-C_2}\cdot\frac{C_2(C_1-C_3)}{C_1C_2-C_3C_4}+\frac{C_3-C_2}{C_1-C_2}\cdot\frac{C_1(C_2-C_4)}{C_1C_2-C_3C_4}=1.
\end{align*}
The proof is complete.
\end{proof}


\begin{remark}
The result of this lemma excludes the existence of $E^*$ again, though we know this fact. By Table \ref{StabilityTable} and this lemma, we can conclude that $E_{12}$ is stable on cells of the third row of Table \ref{tableIV-ii}, A3 and C3, and we will show further that it is GAS.
\end{remark}

\begin{theorem}\label{IV-(ii)-results} Let $\phi(t; p)=(u_1(t), u_2(t), v(t); p)$ be a solution of \eqref{themodel} for initial condition $p\in\mathbb{R}^3_+$. Then the following statements are valid under the same assumptions as in Lemma \ref{EcaseIV-ii}.
\begin{enumerate}[{\rm (i)}]
\item Let parameters be classified in the cells {\rm A4} or {\rm A5} of Table \ref{tableIV-ii}. Then $E_1$ is GAS.
\item For parameters in {\rm A1}, {\rm B1}, {\rm C1} and {\rm D1},
if $a_{32}<1$ $(a_{32}>1)$, then $E_2$ $(E_{2v})$ is GAS. 
\item For parameters in {\rm A3} and {\rm C3}, 
$E_{12}$ is GAS.
\item For parameters in {\rm A2} and {\rm B2},
$E_{12}$ exists and is a saddle with a two-dimensional stable manifold.
\end{enumerate}
\end{theorem}
\begin{proof}
(i) It is obvious that parameters satisfy the assumptions of  Lemma \ref{C1dominate}, $C_1>\max\{C_3, C_4\}$ and $\min\{C_3, C_4\}>C_2$. Consequently, we have $\lim_{t\to\infty}u_2(t)=0$ and system \eqref{themodel} approaches subsusyem $H_{1v}$ asymptotically. Therefore $E_1$ is GAS since $a_{31}<1$ by Theorem \ref{2dim}(a).

(ii) By Lemma \ref{EcaseIV-ii}(i), we have $\lim_{t\to\infty}u_1^{a_{23}}(t)v^\xi(t)$=0 for some positive number $\xi$. Based on this limit, logically, we consider two mutually exclusive cases: the first case is that $u_1(t_k)\to 0$ as $k\to\infty$ for some sequence $\{t_k\}$. For the other case, there is a positive number $\epsilon$ such that $u_1(t)\ge\epsilon$ for all $t\ge 0$.

%

Let the first case happen, then there is a point $q_1\in H_{2v}$ such that $q_1\in\omega(p)$. On $H_{2v}$, if $a_{32}<1$ then  $\lim_{t\to\infty}\phi(t;q_1)=E_2$ by Theorem \ref{2dim}(b). This also implies that $E_2\in\omega(p)$ by invariance of the omega limit set. However, $E_2$ is stable in $\mathbb{R}^3_+$ by Lemma \ref{stabilityE012}(c) and Table \ref{StabilityTable}, hence $\lim_{t\to\infty}\phi(t; p)=E_2$. Alternatively, if $a_{32}>1$  then  $\lim_{t\to\infty}\phi(t;q_1)=E_{2v}$ by Theorem \ref{2dim}(b). This also implies that, similarly, $E_{2v}\in\omega(p)$ by invariance of the omega limit set. And $E_{2v}$ is stable in $\mathbb{R}^3_+$ by Lemma \ref{stabilityE12v}(b) and Table \ref{StabilityTable}, hence $\lim_{t\to\infty}\phi(t; p)=E_{2v}$. This completes the proof of the first case.

For the other case, if $u_1(t)\ge\epsilon$ for some positive number $\epsilon$ and time $t$ is large. Consequently, this implies 
$\lim_{t\to\infty}v(t)=0$, 
and system \eqref{themodel} is asymptotically approaches to subsystem $H_{12}$ \eqref{u1u2system}. By Theorem \ref{2dim}(c)(ii), $E_2$ is GAS on $H_{12}$ which implies $E_2\in\omega(p)$. If $a_{32}<1$, then the equilibrium $E_2$ is also stable on $\mathbb{R}^3$. Hence $\lim_{t\to\infty}\phi(t; p)=E_2$. However, if $a_{32}>1$ then $E_2$ is saddle with the positive eigenvalue $a_{32}-1$ and eigenvector $[0, \frac{a_{23}}{-r_2-a_{32}+1}v_3, v_3]$  with $v_3>0$ by some easy calculations of the Jacobian matrix \eqref{JE2} of $E_2$. By the Butler-McGehee Lemma, there is a point $q_2$ of the unstable manifold of $E_2$ such that $q_2\in\omega(p)$, and this fact contradicts the vanishing of the coordinate $v$. So this inequality $a_{32}>1$ cannot happen, and we have completed the proof of part (ii).

(iii) By Lemma \ref{EcaseIV-ii}(ii), the inequality \eqref{E12stable} holds, and we can obtain that $E_{12}$ is the only stable equilibrium by checking Table \ref{StabilityTable}. By applying Lemma \ref{EcaseIV-ii}(i), the proof of this part is almost the same as part (ii). The only difference is that we use the stability of $E_{12}$ on $\mathbb{R}^3_+$ instead of $E_1$ because of $C_1>C_3$ and $C_2>C_4$. So, we omit it.

(iv) By checking Table \ref{StabilityTable} and Lemma \ref{EcaseIV-ii}, it is obvious that $E_{12}$ exists and is a saddle with a two-dimensional stable manifold for cells A2 and B2, and the proof is completed.
\end{proof}

\begin{remark}
In this subsection, we show the non-existence of $E^*$ and the dynamics of \eqref{themodel} presented in Table \ref{tableIV-ii} are completely verified analytically by Theorem \ref{IV-(ii)-results}. Since Theorem \ref{IV-(ii)-results} (i) proves the results in the cells {\rm A4} and {\rm A5}, (ii) verifies the first row, and (iii) shows the third row of Table \ref{tableIV-ii}. On the second row of Table \ref{tableIV-ii}, similar to the results of category III-(ii), it is easy to check, by Table \ref{StabilityTable}, that $E_1$, $E_2$, and $E_{2v}$ are stable simultaneously. Since \eqref{E12stable} holds, by Lemma \ref{EcaseIV-ii}(ii), $E^*$ does not exist and $E_{12}$ is a saddle with a two-dimensional stable manifold separating the absorbing regions of the boundary stable nodes. 

It is also worth noting that the dynamics of all sub-cases of this category IV-(ii) are dominated by the negative competition. Since the win-loss state, in rows 4 and 5 of Table \ref{tableIV-ii}, implies $E_1$ is GAS; the loss-win state, in row 1 of Table \ref{tableIV-ii}, implies $E_2/E_{2v}$ is GAS; the win-win state, in row 3 of Table \ref{tableIV-ii}, implies $E_{12}$ is GAS; the loss-loss state, in row 2 of Table \ref{tableIV-ii}, implies the bistability $E_1/E_2$ or $E_1/E_{2v}$. All results are consistent with those of the two-dimensional $u_1$-$u_2$ subsystem.
\end{remark}

\subsection{The Cases IV-(iii) and IV-(iv)}
In this subsection, we investigate the final two most interesting cases, IV-(iii) and IV-(iv). Using the same classification of IV-(ii), similarly, we present two tables of all possible subcases of cases IV-(iii) and IV-(iv), Table \ref{tableIV-iii} and Table \ref{tableIV-iv}.

By Lemma \ref{C1dominate}, we expect that $E_{12}$ and $E^*$ cannot exist in the classification of IV-4 and IV-5 (rows 4 and 5 of Table \ref{tableIV-iii} and \ref{tableIV-iv}). In particular, $E_{12}$ exists only for IV-2 and IV-3 (rows 2 and 3 of Table \ref{tableIV-iii} and \ref{tableIV-iv}). Moreover, $E_{12}$ is unstable for IV-2, and it is stable for IV-3 if \eqref{E12stable} holds. For the existence of $E^*$, we first establish a result to simplify the condition of verification.

\begin{lemma}\label{vpositiveenough}
For the {\rm IV-1}, {\rm IV-2} and {\rm IV-3} classifications of the cases {\rm IV-(iii)} or {\rm IV-(iv)}, $E^*$ exists if and only if \eqref{vpIII-ii} or \eqref{vpositive} holds, repectively.

\end{lemma}
\begin{proof}
Recall that $E^*$ exists if and only if $u_1^*<1$, $u^*_2<1$ and $v^*>0$ which are equivalent to \eqref{u1pIII-ii}-\eqref{vpIII-ii} or  \eqref{u1sc}-\eqref{vpositive} for cases IV-(iii) or IV-(iv), respectively.
We only show that \eqref{vpIII-ii} implies  \eqref{u1pIII-ii} and \eqref{u2pIII-ii}, and the other case is similar.

By \eqref{vpIII-ii}, we have
\begin{align*}
a_{32}C_1(C_4-C_2)&>-a_{31}C_2(C_3-C_1)-C_1C_2+C_3C_4\\
&=(1-a_{31})C_1(C_3-C_2)+C_3(a_{31}C_2-a_{31}C_1-C_1+C_4)\\
&>(1-a_{31})C_1(C_3-C_2), 
\end{align*}
since $(C_1-C_2)/(C_1-C_4)>1/a_{31}$. Hence \eqref{u1pIII-ii} holds. Similarly, 
\begin{align*}
a_{31}&C_2(C_3-C_1)>-a_{32}C_1(C_4-C_2)-C_1C_2+C_3C_4\\
&=(1-a_{32})C_2(C_4-C_1)+C_4(-a_{32}C_1+a_{32}C_2-C_2+C_3)\\
&>C_2(1-a_{32})(C_1-C_4), 
\end{align*}
since $(C_1-C_2)/(C_3-C_2)<1/a_{32}$. Hence \eqref{u2pIII-ii} also holds. The proof is completed.
\end{proof}

\subsubsection{The Case IV-(iii)}
For case IV-(iii), we rule out 9 cells labeled ``contradiction '' in the previous similar way, and please see Table \ref{tableIV-iii}. Examining the local stability of each equilibrium in Table \ref{StabilityTable} for each cell of Table \ref{tableIV-iii}, we will show some global dynamics and demonstrate bistability.

\begin{table}[htp]
\begin{threeparttable}
\caption{\textcolor{black}{All 20 subcases of the category IV-(iii) under inequalities $\frac{C_1-C_2}{C_1-C_4}>\frac{1}{a_{31}}$ and $0<\frac{C_1-C_2}{C_3-C_2}<\frac{1}{a_{32}}$. }}
\begin{center}
\begin{tabular}{|c||c|c|c|c|}
\hline
\hline
&A &B&C&D\\
& \footnotesize{$\max\{a_{31}, a_{32}\}<1$} &  \footnotesize{$a_{31}<1<a_{32}$} & \footnotesize{$a_{32}<1<a_{31}$} & \footnotesize{$\min\{a_{31}, a_{32}\}>1$}  \\
\hline
\hline
1\ :\ \red{\scriptsize (L, W, W, L)}&  \multirow{2}{*}{\blue{contradiction} } &\multirow{2}{*}{\blue{contradiction} }&\red{Bistability}&\red{Bistability}\\
{\scriptsize $C_3>C_1>C_2>C_4$} &&&$E_{1v}$/$E_2$&$E_{1v}$/$E_{2v}$\\
\hline
2\ :\ \red{\scriptsize (L, L, W, L)} &\black{Bistability} &\black{Bistability} &\black{Bistability}  &\black{Bistability}\\
 {\scriptsize $C_3>C_1>C_4>C_2$} &$E_1$/$E_2$& $E_1$/$E_{2v}$  & $E_{1v}$/$E_2$&$E_{1v}$/$E_{2v}$\\
\hline
3\ :\ \red{\scriptsize (W, W, W, L)} & \multirow{2}{*}{\blue{contradiction} } &\multirow{2}{*}{\blue{contradiction} }& \multirow{2}{*}{\red{$E_{1v}/E_{12}$}\tnote{1}}&\multirow{2}{*}{\blue{contradiction}} \\
{\scriptsize $C_1>C_3>C_2>C_4$} &&&&\\
\hline
4\ :\ \red{\scriptsize (W, L, W, L)} &\multirow{2}{*}{\black{$E_1$ is GAS} }& \multirow{2}{*}{\blue{contradiction} } &\multirow{2}{*}{\black{$E_{1v}$ is GAS} }&  \multirow{2}{*}{\blue{contradiction}} \\
 {\scriptsize $C_1>C_3>C_4>C_2$} &&&&\\
\hline
5\ :\ \red{\scriptsize (W, L, W, L)} &\multirow{2}{*}{\black{$E_1$ is GAS} }& \multirow{2}{*}{\blue{contradiction} } &\multirow{2}{*}{\black{$E_{1v}$ is GAS}}&\multirow{2}{*}{\blue{contradiction}} \\
 {\scriptsize $C_1>C_4>C_3>C_2$} &&&&\\
\hline
\hline
\end{tabular}
\begin{tablenotes}
\item[1] If \eqref{E12stable} holds, then the bistabilty of $E_{1v}/E_{12}$ happens. Otherwise, if \eqref{E12unstable} holds, then $E_{1v}$ is GAS. 
\end{tablenotes}
\end{center}
\label{tableIV-iii}
\end{threeparttable}
\end{table}%

We first show the following result, which can be applied to all subcases of IV-(iii)
\begin{lemma}\label{IV-iiiE*unstable}
In this category {\rm IV-(iii)}, $E^*$ is unstable whenever it exists.
\end{lemma}
\begin{proof}
Because  
\begin{align*}
\frac{C_1-C_2}{C_1-C_4}>\frac{1}{a_{31}} \ \text{ and }\ \frac{C_1-C_2}{C_3-C_2}<\frac{1}{a_{32}}\ \text{ imply }\ \frac{1}{a_{32}}\frac{C_1-C_2}{C_1-C_4}>\frac{1}{a_{31}} \frac{C_1-C_2}{C_3-C_2}
\end{align*}
which is equivalent to the inequality 
$ (C_1-C_4)a_{32}+(C_2-C_3)a_{31}<0$
and contradicts to \eqref{RH-2-1}. Hence, for IV-(iii), $E^*$ is unstable if it exists. The proof is completed.
\end{proof}

\begin{proposition}\label{IV-iii-r1}
Let the parameters be classified in the cells {\rm A4}, {\rm A5}, {\rm C4} and {\rm C5} of Table \ref{tableIV-iii}. If $a_{31}<1$ then $E_1$ is GAS, otherwise if $a_{31}>1$ then $E_{1v}$ is GAS.
\end{proposition}
\begin{proof}
Let $\phi(t)=(u_1(t), u_2(t), v(t))$ be a solution of \eqref{themodel}, and we have $\lim_{t\to\infty}u_2(t)=0$ since the parameter classified in the cells {\rm A4}, {\rm A5}, {\rm C4} and {\rm C5} of Table \ref{tableIV-iii} satisfy the assumptions of Lemma \ref{C1dominate}. Hence, system \eqref{themodel} approaches $H_{1v}$ subsystem asymptotically. By Theorem \ref{2dim}(a), the results can be easily established. We have completed the proof.
\end{proof}

\begin{proposition}\label{IV-iii-bistability}
For parameters in cells {\rm C1}, {\rm D1}, {\rm A2}, {\rm B2}, {\rm C2} and {\rm D2} of Table \ref{tableIV-iii}, the following statements are valid.
\begin{enumerate}[{\rm (i)}]
\item For parameters in cells {\rm C1} and  {\rm D1}, the positive equilibrium $E^*$ exists and is unstable with a two-dimensional stable manifold. For {\rm C1}, the boundary equilibria $E_{1v}$ and $E_2$ are stable, and the boundary equilibria $E_{1v}$ and $E_{2v}$ are stable for {\rm D1}.
\item For parameters in {\rm A2}, $E^*$ does not exist, $E_{12}$ is a saddle with a two-dimensional stable manifold, and the boundary equilibria $E_1$ and $E_2$ are stable.
\item For parameters in cells {\rm B2}, {\rm C2} and {\rm D2}, if \eqref{E12stable} holds, then $E^*$ does not exist and $E_{12}$ is a saddle with a two-dimensional stable manifold. Otherwise, if \eqref{E12unstable} holds, then $E^*$ exists and is a saddle with a two-dimensional stable manifold. Moreover, the boundary equilibria $E_1$/$E_{2v}$, $E_{1v}$/$E_2$ and $E_{1v}$/$E_{2v}$ are stable with respect to parameters in the cells {\rm B2}, {\rm C2} and {\rm D2}.
\end{enumerate}
\end{proposition}
\begin{proof}
For all local stabilities of the boundary equilibria, including $E_1$, $E_2$, $E_{1v}$, and $E_{2v}$, they are easy to be verified by Table \ref{StabilityTable}. So we omit them.

(i) By Lemma \ref{vpositiveenough}, we check the existence of $E^*$ by verifying \eqref{vpIII-ii}. We rewrite the equivalent form of \eqref{vpIII-ii}, 
\begin{align}
a_{31}&C_2(C_3-C_1)-a_{32}C_1(C_2-C_4)+C_1C_2-C_3C_4\nonumber \\
&= a_{31}C_2(C_3-C_1)-C_2(C_3-C_1)-a_{32}C_1(C_2-C_4)+C_3(C_2-C_4)\nonumber \\
&=(a_{31}-1)C_2(C_3-C_1)+(C_3-a_{32}C_1)(C_2-C_4)>0.\label{vpeIV-iii}
\end{align}
It is easy to see that, for parameter in C1, \eqref{vpeIV-iii} holds, hence, $E^*$ exists, and it is unstable by Lemma \ref{IV-iiiE*unstable}. Moreover, its stable manifold is two-dimensional by Lemma \ref{bistabilitylemma} (ii).

For parameters in D1, it is sufficient to show that \eqref{vpeIV-iii} is also true. If $C_3\ge a_{32}C_1$, the inequality \eqref{vpeIV-iii} is obvious true. On the other hand, we assume that $C_3< a_{32}C_1$. Then the inequality $\frac{C_1-C_2}{C_1-C_4}>\frac{1}{a_{31}}$ implies
\[(C_3-a_{32}C_1)(C_2-C_4)>(a_{31}-1)(C_3-a_{32}C_1)(C_1-C_2).\]
Consider 
\begin{align*}
(a_{31}-1)&C_2(C_3-C_1)+(C_3-a_{32}C_1)(C_2-C_4)\\
&>(a_{31}-1)C_2(C_3-C_1)+(a_{31}-1)(C_3-a_{32}C_1)(C_1-C_2)\\
&=(a_{31}-1)\big(C_2(C_3-C_1)+(C_3-a_{32}C_1)(C_1-C_2)\big)\\
&=(a_{31}-1)C_1(C_3-C_2-a_{32}C_1+a_{32}C_2)>0, 
\end{align*}
since $\frac{C_1-C_2}{C_3-C_2}<\frac{1}{a_{32}}$. Hence \eqref{vpeIV-iii} holds.

(ii) For parameters in A2, it is easy to be verified that $E^*$ doest not exist and $E_{12}$ is a saddle with two-dimensional stable manifold by \eqref{vpeIV-iii} and \eqref{JE12}, the Jacobian matrix of $E_{12}$.

(iii) If \eqref{E12stable} holds, then $E^*$ does not exist, and $E_{12}$ is a saddle with the two-dimensional stable manifold for the same reasons as part (ii). On the other hand, if \eqref{vpIII-ii} holds, then $E^*$ exists, and it is a saddle with the two-dimensional stable manifold by Lemma \ref{IV-iiiE*unstable} and \ref{bistabilitylemma}(ii).
Hence, all proofs are completed.
\end{proof}

%

\begin{proposition}\label{IV-iiiC3}
Let $\phi(t)=\big(u_1(t), u_2(t), v(t)\big)$ be a solution of \eqref{themodel}  and parameters classified in cells  {\rm C3} of Table \ref{tableIV-iii}. Then the following statements are valid.
\begin{enumerate}[{\rm (i)}]
\item If \eqref{E12stable} holds, then the bistability of $E_{1v}/E_{12}$ happens.
\item If \eqref{E12unstable} holds, then $\lim_{t\to\infty} u_2(t)=0$ and $E_{1v}$ is GAS.
\end{enumerate}
\end{proposition}
If \eqref{E12stable} holds, then it is clear that $E_{1v}$ is stable, $E_{12}$ is stable, and $E^*$ does not exist by checking Table \ref{StabilityTable}. This shows part (i). For part (ii), we first show the following technical lemma.

\begin{lemma}
Let the parameters be classified in the cell {\rm C3} of Table \ref{tableIV-iii}.
If  \eqref{E12unstable} holds, then there exist two positive numbers $\eta_0$, $\xi_0$ such that these inequalities 
\begin{align*}
&\frac{C_2}{C_1}<\eta_0<1,
& \max\left\{\frac{\eta_0 C_1-C_4}{a_{31}}, \frac{\eta_0 C_3-C_2}{a_{32}}\right\}<\xi_0<\eta_0 C_1-C_2,
\end{align*}
are valid.
\end{lemma}
\begin{proof}
Define two functions of $\eta$ on domain $(C_2/C_1, 1]$
\begin{align*}
f(\eta)= \frac{\eta C_1-C_4}{\eta C_1-C_2} \quad\text{ and } \quad g(\eta)=\frac{\eta C_3-C_2}{\eta C_1-C_2},
\end{align*}
and it is easy to see that 
\begin{align*} 
&1<f(1)<a_{31}, \quad a_{32}<g(1)<1, \quad f'(\eta)<0, \quad g'(\eta)>0, \\
& \lim_{\eta\to(\frac{C_2}{C_1})^+}f(\eta)=\infty, \ \text{ and }\lim_{\eta\to(\frac{C_2}{C_1})^+}g(\eta)=-\infty. 
\end{align*}
By the Intermediate Value Theorem of the continuous function, there are two unique positive numbers $\eta_1, \eta_2\in(C_2/C_1, 1)$ such that $f(\eta_1)=a_{31}$ and $g(\eta_2)=a_{32}$, respectively. It is easy to solve
\begin{align*}
 \eta_1=\frac{a_{31}C_2-C_4}{C_1(a_{31}-1)},  \qquad\eta_2=\frac{C_2(1-a_{32})}{C_3-a_{32}C_1}, 
\end{align*}
and verify that the inequality \eqref{E12unstable} is equivalent to
\begin{align*}
(a_{31}-1)C_1C_2+a_{32}C_1(C_2-C_4)>C_3(a_{31}C_2-C_4),
\end{align*}
which is also equivalent to $\eta_1<\eta_2$. This implies that $f(\eta_2)<a_{31}$ and $g(\eta_2)=a_{32}$. By the continuity of the functions $f$ and $g$, we take a positive number $\epsilon$ such that $\eta_0=\eta_2-\epsilon$, $f(\eta_0)<a_{31}$ and $0<g(\eta_0)<a_{32}$. More clearly, it means that 
\begin{align*}
f(\eta_0)=\frac{\eta_0 C_1-C_4}{\eta_0 C_1-C_2}<a_{31} \qquad\text{ and } \qquad g(\eta_0)=\frac{\eta_0 C_3-C_2}{\eta_0 C_1-C_2}<a_{32},
\end{align*}
which implies that 
\[ \max\left\{ \frac{\eta_0 C_1-C_4}{a_{31}}, \frac{\eta_0 C_3-C_2}{a_{32}} \right\}<\eta_0 C_1-C_2. \]
Therefore, we can choose $\xi_0$ easily, and the proof is completed.

\end{proof}

\begin{proof}[Proof of Proposition \ref{IV-iiiC3}\ {\rm (ii)}]
Let constants $\eta_0$ and $\xi_0$ be defined on the previous lemma. Consider
\begin{align*}
a_{13}\frac{\dot u_2}{u_2}-&\eta_0 a_{23}\frac{\dot u_1}{u_1}-\xi_0\frac{\dot v}{v}\\
&\le C_2(1-u_2)-C_4u_1-\eta_0C_1(1-u_1)+\eta_0C_3u_2+\xi_0-\xi_0a_{31}u_1-\xi_0a_{32}u_2\\
&\le C_2-\eta_0 C_1+\xi_0+(\eta_0 C_1-C_4-\xi_0 a_{31})u_1+(\eta_0 C_3-C_2-\xi_0 a_{32})u_2\\
&\le C_2-\eta_0 C_1+\xi_0<0,
\end{align*}
which implies $\lim_{t\to\infty}u_2(t)=0$. It is clear that $E_{1v}$ is GAS by the Markus Limiting Theorem. We have completed the proof.
\end{proof}

\begin{remark}
In this subsection, we establish all global dynamics of \eqref{themodel} with parameters in case IV-(iii). First, we show the GAS of $E_1$ or $E_{1v}$, for the fourth and fifth rows of Table \ref{tableIV-iii} by Proposition \ref{IV-iii-r1}. Secondly, we show that $E_{12}$ exists in the second and third rows, and that $E^*$ could exist in rows 1-3, depending on whether \eqref{E12stable} or \eqref{E12unstable} holds. In particular, in the second row, $E_{12}$ is a saddle with a two-dimensional stable manifold or a one-dimensional stable manifold if \eqref{E12stable} or \eqref{E12unstable} holds, respectively. In the third row, $E_{12}$ is a stable node or a saddle with a one-dimensional stable manifold if \eqref{E12stable} or \eqref{E12unstable} holds, respectively. Finally, we find two stable nodes on cells C1 and D1. However, $E_{12}$ does not exist on the first row; hence, we show that $E^*$ indeed exists and is a saddle with a two-dimensional stable manifold separating the absorbing regions of these two stable nodes.

The dynamics of all subcases of this category IV-(iii) are NOT dominated by the negative competition anymore. For the second, fourth, and fifth rows, corresponding to the loss-loss, win-loss, and win-loss states in negative competition, respectively, these subcases are compatible with those of the two-dimensional $u_1$-$u_2$ subsystem. However, for the first row with the loss-win state under negative competition, we show bistability, and for the third row, with the win-win state, we observe different dynamics. These issues will be discussed in the last section.
\end{remark}

\subsubsection{The Case IV-(iv)}
For the most interesting and complicated case IV-(iv), we also present the parameters classification in Table \ref{tableIV-iv}.


\begin{table}[htp]
\begin{threeparttable}
\caption{\textcolor{black}{All 20 subcases of the category IV-(iv) under inequalities $0<\frac{C_1-C_2}{C_1-C_4}<\frac{1}{a_{31}}$ and $\frac{C_1-C_2}{C_3-C_2}>\frac{1}{a_{32}}$. }}
\begin{center}
\begin{tabular}{|c||c|c|c|c|}
\hline
\hline
&A &B&C&D\\
& \footnotesize{$\max\{a_{31}, a_{32}\}<1$} & \footnotesize{$a_{31}<1<a_{32}$} & \footnotesize{$a_{32}<1<a_{31}$} & \footnotesize{$\min\{a_{31}, a_{32}\}>1$}  \\
\hline
\hline
1\ :\ \red{\scriptsize (L, W, W, L)}&  \multirow{2}{*}{\blue{contradiction}} &\multirow{2}{*}{\textcolor{red}{Complex\tnote{1}}}&\multirow{2}{*}{\blue{contradiction} }&\multirow{2}{*}{\red{$E^*$ is stable\tnote{2}}}\\
{\scriptsize $C_3>C_1>C_2>C_4$} &&&&\\
\hline
2\ :\ \red{\scriptsize (L, L, W, L)} &\multirow{2}{*}{\blue{contradiction}}&\multirow{2}{*}{\red{$E_1$ is GAS\tnote{3}}} &\multirow{2}{*}{\blue{contradiction}} &\multirow{2}{*}{\blue{contradiction}}\\
 {\scriptsize $C_3>C_1>C_4>C_2$} &&&&\\
 \hline
3\ :\ \red{\scriptsize (W, W, W, L)}&\black{\multirow{2}{*}{$E_{12}/E^*$\tnote{4}}}& \black{\multirow{2}{*}{$E_{12}/E^*$\tnote{4}}} &\black{\multirow{2}{*}{$E_{12}/E^*$\tnote{4}}} & \black{$E^*$ exists} \\
{\scriptsize $C_1>C_3>C_2>C_4$} &&&& and is stable\\
\hline
4\ :\ \red{\scriptsize (W, L, W, L)} &\multirow{2}{*}{$E_1$ is GAS}& \multirow{2}{*}{\black{$E_1$ is GAS}} &\multirow{2}{*}{\blue{contradiction}}&  \multirow{2}{*}{\blue{contradiction}} \\
 {\scriptsize $C_1>C_3>C_4>C_2$} &&&&\\
\hline
5\ :\ \red{\scriptsize (W, L, W, L)} &\multirow{2}{*}{$E_1$ is GAS}& \multirow{2}{*}{\black{$E_1$ is GAS}} &\multirow{2}{*}{\blue{contradiction}}&\multirow{2}{*}{\blue{contradiction}} \\
 {\scriptsize $C_1>C_4>C_3>C_2$} &&&&\\
\hline
\hline
\end{tabular}
\begin{tablenotes}
\item[1] In this cell, we numerically found the stable positive equilibrium, the periodic solution, or complex dynamics for different parameters.
\item[2] In this cell, $E^*$ is GAS by numerical observation.
\item[3] If \eqref{E12unstable} holds, then we can show that $E_1$ is GAS analytically. However, if \eqref{E12unstable} does not hold, then $E_1$ is also GAS by numerical observation.
\item[4] In these three cells, the equilibrium $E_{12}$ exists, and the notation $E_{12}/E^*$ means that either $E_{12}$ is stable if \eqref{E12stable} holds or $E^*$ exists if \eqref{E12unstable} holds. Furthermore, if $E^*$ exists, then it is stable, and $E_{12}$ is unstable.
\end{tablenotes}
\end{center}
\label{tableIV-iv}
\end{threeparttable}
\end{table}%

We first establish the following simple global dynamics in the cells {\rm A4}, {\rm A5}, {\rm B4} and {\rm B5} of Table \ref{tableIV-iv}. Since the assumptions of Lemma \ref{C1dominate} are satisfied, by similar arguments in Theorem \ref{IV-(ii)-results}(i), $E_1$ is GAS.

For the inequalities, 
$\frac{C_1-C_2}{C_1-C_4}<\frac{1}{a_{31}}$ and $\frac{C_1-C_2}{C_3-C_2}>\frac{1}{a_{32}}$ imply
\begin{align*}
\frac{C_1-C_2}{C_1-C_4}\cdot\frac{1}{a_{32}}<\frac{C_1-C_2}{C_3-C_2}\cdot\frac{1}{a_{31}}
\end{align*}
which implies \eqref{RH-2-1}. 
On the third row of the classification of Table \ref{tableIV-iv}, $C_1>C_3$ implies $r_1>a_{12}$ and $C_2>C_4$ implies $r_2>a_{21}$. Combining these two inequalities and Lemma \ref{RH13}, the inequalities $\eqref{RH-1}$ and $\eqref{RH-3}$ hold, hence we establish the following results.

\begin{lemma}\label{IV-ivrow3stability}
In the cells {\rm A3}-{\rm D3} of Table \ref{tableIV-iv}, the positive equilibrium $E^*$ is stable if it exists.
\end{lemma}

\begin{proposition}\label{IV-ivE*}
For parameters in cells {\rm B1}, {\rm D1}, {\rm A3}-{\rm D3} of  Table \ref{tableIV-iv}, the following statements are valid.
\begin{enumerate}[{\rm (i)}]
\item For parameters in cells {\rm B1}, {\rm D1} and {\rm D3}, the positive equilibrium $E^*$ exists, and it is stable in {\rm D3}. Moreover, system \eqref{themodel} is uniformly persistent.
\item For parameters in the cells {\rm A3},  {\rm B3} and {\rm C3}, if \eqref{E12stable} holds then $E_{12}$ is GAS. Otherwise, if \eqref{E12unstable} holds, then the positive equilibrium $E^*$ exists and is stable. In addition, system \eqref{themodel} is uniformly persistent.
\end{enumerate}
\end{proposition}
\begin{proof}
(i) Recall that $E^*$ exists if and only if $u^*_1<1$, $u_2^*<1$ and $v^*>0$ which are equivalent to \eqref{u1sc}, \eqref{u2sc} and \eqref{vpositive}, respectively, and it is sufficient to show that \eqref{vpositive} holds by Lemma \ref{vpositiveenough}. With the equivalent form of \eqref{vpositive},
\begin{align}\label{E2existenceIV-ivB1}
(1-a_{31})C_2(C_3-C_1)+(a_{32}C_1-C_3)(C_2-C_4)>0,
\end{align}
it is obvious that $E^*$ exists for D3, and it is also stable by Lemma \ref{IV-ivrow3stability}.

For parameters in B1, it is sufficient to show that $a_{32}C_1-C_3>0$. Since 
\[ \frac{C_1-C_2}{C_3-C_2}>\frac{1}{a_{32}} \ \text{ implies } \ a_{32}C_1-C_3>(a_{32}-1)C_2>0.\]

For parameters in D1,
\[0<\frac{C_1-C_2}{C_1-C_4}<\frac{1}{a_{31}}\ \text{ if and only if } \ (a_{31}-1)(C_1-C_2)<C_2-C_4,\]
which implies 
\[(C_2-C_4)(a_{32}C_1-C_3)>(a_{31}-1)(a_{32}C_1-C_3)(C_1-C_2).\]
Then
\begin{align*}
(1-a_{31})&C_2(C_3-C_1)+(a_{32}C_1-C_3)(C_2-C_4)\\
&>(1-a_{31})C_2(C_3-C_1)+(a_{31}-1)(a_{32}C_1-C_3)(C_1-C_2)\\
&=(a_{31}-1)\big((a_{32}C_1-C_3)(C_1-C_2)-C_2(C_3-C_1)\big)\\
&=(a_{31}-1)\big(a_{32}(C_1-C_2)-(C_3-C_2)\big)>0, 
\end{align*}
since $\frac{C_1-C_2}{C_3-C_2}>\frac{1}{a_{32}}$. 

The proof of the uniform persistence for system \eqref{themodel} is similar to Proposition \ref{ThoremII-ii-1}, so we omit it. The proof of this part is complete.

(ii) We only need to show that if \eqref{E12stable} holds, then $E_{12}$ is GAS. The proof is almost the same as Lemma \ref{E12lemma}, so we omit it. The proof of the uniform persistence for system \eqref{themodel} is similar to Proposition \ref{ThoremII-ii-1}, so we omit it. The proof is complete.
\end{proof}

Next, we show the global dynamics of the second row of Table \ref{tableIV-iv}.
\begin{proposition}\label{IV-ivB2}
Let $\phi(t)=\big(u_1(t), u_2(t), v(t)\big)$ be a solution of \eqref{themodel}  and parameters classified in cells  {\rm B2} of Table \ref{tableIV-iv}. If  \eqref{E12unstable} holds, then we have $\lim_{t\to\infty}u_2(t)=0$ and $E_1$ is GAS.
\end{proposition}
\begin{lemma}
If the parameters are classified in the cell {\rm B2} of Table \ref{tableIV-iv} and \eqref{E12unstable} holds, then there exist two positive numbers $\eta_0$, $\xi_0$ such that these inequalities 
\begin{align*}
&\frac{C_2}{C_1}<\eta_0<1,
& \max\left\{\frac{\eta_0 C_1-C_4}{a_{31}}, \frac{\eta_0 C_3-C_2}{a_{32}}\right\}<\xi_0<\eta_0 C_1-C_2,
\end{align*}
are valid.
\end{lemma}
\begin{proof}
Define two functions of $\eta$ on domain $(C_2/C_1, 1]$
\begin{align*}
f(\eta)= \frac{\eta C_1-C_4}{\eta C_1-C_2} \quad\text{ and } \quad g(\eta)=\frac{\eta C_3-C_2}{\eta C_1-C_2},
\end{align*}
and it is easy to see that 
\begin{align*} 
&a_{31}<f(1)<1, \ 1<g(1)<a_{32}, \ f'(\eta)>0, \ g'(\eta)<0, \\& \lim_{\eta\to(\frac{C_2}{C_1})^+}f(\eta)=-\infty, \ \text{ and }\lim_{\eta\to(\frac{C_2}{C_1})^+}g(\eta)=\infty. 
\end{align*}
By the Intermediate Value Theorem of the continuous function, there are two unique positive numbers $\eta_1, \eta_2\in(C_2/C_1, 1)$ such that $f(\eta_1)=a_{31}$ and $g(\eta_2)=a_{32}$, respectively. It is easy to solve
\begin{align*}
 \eta_1=\frac{C_4-a_{31}C_2}{C_1(1-a_{31})},  \qquad\eta_2=\frac{C_2(a_{32}-1)}{a_{32}C_1-C_3}, 
\end{align*}
and verify that the inequality \eqref{E12unstable} is equivalent to
\begin{align*}
(1-a_{31})C_1C_2+a_{32}C_1(C_4-C_2)>C_3(C_4-a_{31}C_2),
\end{align*}
which is also equivalent to $\eta_1>\eta_2$. This implies that $f(\eta_1)=a_{31}$ and $g(\eta_1)<a_{32}$. By the continuity of the functions $f$ and $g$, we take a positive number $\epsilon$ such that $\eta_0=\eta_1-\epsilon$, $0<f(\eta_0)<a_{31}$ and $g(\eta_0)<a_{32}$. More clearly, it means that 
\begin{align*}
\frac{\eta_0 C_1-C_4}{\eta_0 C_1-C_2}<a_{31} \qquad\text{ and } \qquad\frac{\eta_0 C_3-C_2}{\eta_0 C_1-C_2}<a_{32},
\end{align*}
which implies that 
\[ \max\left\{ \frac{\eta_0 C_1-C_4}{a_{31}}, \frac{\eta_0 C_3-C_2}{a_{32}} \right\}<\eta_0 C_1-C_2. \]
Therefore, we can choose $\xi_0$ easily, and the proof is completed.

\end{proof}

\begin{proof}[Proof of Proposition \ref{IV-ivB2}]
Let constants $\eta_0$ and $\xi_0$ be defined on the previous lemma. Consider
\begin{align*}
a_{13}\frac{\dot u_2}{u_2}-&\eta_0 a_{23}\frac{\dot u_1}{u_1}-\xi_0\frac{\dot v}{v}\\
&\le C_2(1-u_2)-C_4u_1-\eta_0C_1(1-u_1)+\eta_0C_3u_2+\xi_0-\xi_0a_{31}u_1-\xi_0a_{32}u_2\\
&\le C_2-\eta_0 C_1+\xi_0+(\eta_0 C_1-C_4-\xi_0 a_{31})u_1+(\eta_0 C_3-C_2-\xi_0 a_{32})u_2\\
&\le C_2-\eta_0 C_1+\xi_0<0,
\end{align*}
which implies $\lim_{t\to\infty}u_2(t)=0$. It is clear that $E_1$ is GAS. We have completed the proof.
\end{proof}


\begin{remark}
In this subsection, we establish some global dynamics of \eqref{themodel} with parameters in case IV-(iv).  For the fourth and fifth rows of Table \ref{tableIV-iv}, we show the GAS of $E_1$, which implies the non-existence of $E^*$. Next, in the third row, it is obvious that $E_{12}$ is a stable node or a saddle with a two-dimensional stable manifold, which is dependent on the validation of \eqref{E12stable} or \eqref{E12unstable}, respectively. Meanwhile, the positive equilibrium $E^*$ exists if \eqref{E12unstable} holds. In particular, for cell D3, we show that $E^*$ indeed exists and is stable. In the first and second rows, we show that $E^*$ exists for cells B1 and D1, although we cannot be sure whether it is stable. In cell B2, by Proposition \ref{IV-ivB2}, if \eqref{E12unstable} holds, then $E_1$ is GAS. Moreover, numerical simulations suggest that $E_1$ is also GAS when \eqref{E12stable} holds, though an analytical proof remains open.

Biologically, the dynamics of some sub-cases of this category IV-(iv) are no longer dominated by negative competition, either. For the third, fourth, and fifth rows, corresponding to the win-win, win-loss, and win-loss states in negative competition, respectively, these cases are compatible with those of the two-dimensional $u_1$-$u_2$ subsystem.
However, in the first row, with the loss-win state under negative competition, we observe bistability. And for the third row, with the win-win state, we observe different dynamics than before. Since the win-loss state, in rows 4 and 5 of Table \ref{tableIV-ii}, implies $E_1$ is GAS; the loss-win state, in row 1 of Table \ref{tableIV-ii}, implies $E_2/E_{2v}$ is GAS; the win-win state, in row 3 of Table \ref{tableIV-ii}, implies $E_{12}$ is GAS; the loss-loss state, in row 2 of Table \ref{tableIV-ii}, implies the bistability $E_1/E_2$ or $E_1/E_{2v}$. All results are consistent.
\end{remark}
\section{Discussions, Biological Interpretations and Numerical Simulations}
In this work, we revisit the two-prey-one-predator model with prey competition \eqref{themodel}. Based on hypothesis {\bf (H2)}, we identified 56 meaningful sub-items, all of which were analytically proven or verified through numerical simulations. We are therefore in a position to summarize and finalize this work. 

In Section 6.1, we compare the results of \cite{Vance1978, Gilpin1979, Hsu1981} with our findings. Moreover, to explore the biological significance of our results, we used three indicators- negative competition, apparent competition, and trade-off competition- to categorize the 56 results into different biological meanings in Section 6.2.

\subsection{Mathematical Discussion}

The present work incorporates both direct and indirect interactions by including prey competition in addition to predation effects. This leads to a significantly richer dynamical structure, involving Negative competition and trade-off mechanisms. When the competition coefficients $a_{12}$ and $a_{21}$ between prey are set to $0$, the dynamics of system \eqref{themodel} is simpler \cite{Yang2025}. More importantly, we provide a complete classification of the global dynamics across the entire parameter space, rather than restricting attention to selected parameter regimes. Our analysis reveals that, despite the increased complexity, the dynamics are largely consistent with those of an associated two-dimensional prey-competition subsystem, except in a small number of exceptional cases. This observation highlights an intrinsic structural relationship between the full system and its lower-dimensional counterpart, and provides a unified perspective on the emergence of different dynamical behaviors.

Hsu, in 1981 \cite{Hsu1981}, considered the model 
\begin{align}\label{Hsumodel}
N_1/d\bar t & = N_1\left[\bar r_1(1-N_1/K_1)-\alpha_{12}N_2-\alpha_{13}P\right],\nonumber \\
N_2/d\bar t & = N_2\left[\bar r_2(1-N_2/K_2)-\alpha_{21}N_1-\alpha_{23}P\right],\\
\dot P/d\bar t & = P\left[\alpha_{31}N_1+\alpha_{32}N_2-D\right],\nonumber
\end{align} 
 with assumptions 
 \begin{align}\label{Hsuassumption}
 K_1>D/\alpha_{31}\quad \text{ and } \quad K_2>D/\alpha_{32}, 
 \end{align}
 and showed the following results.
 \begin{theorem}[\cite{Hsu1981}]\label{Hsuresults}
Let $E_{1p}=(\bar N_1, 0, \bar P_1)$ and $E_{2p}=(0, \bar N_2, \bar P_2)$ be the unique equilibria of the sub-community $(N_1, P)$ and $(N_2, P)$, respectively. Then the following statements are valid.
\begin{enumerate}[{\rm (i)}]
\item  If $\frac{\bar r_1-\alpha_{12}\bar N_2}{\bar r_2(1-\bar N_2/K_2)}>\frac{a_{13}}{\alpha_{23}}$ and $\frac{\bar r_2-\alpha_{21}\bar N_1}{\bar r_1(1-\bar N_1/K_1)}<\frac{\alpha_{23}}{\alpha_{13}}$, then $E_{1p}$ is globally stable {\rm \cite[Theorem 4.3(i)]{Hsu1981}}.
\item If $\frac{\bar r_1-\alpha_{12}\bar N_2}{\bar r_2(1-\bar N_2/K_2)}<\frac{a_{13}}{\alpha_{23}}$ and $\frac{\bar r_2-\alpha_{21}\bar N_1}{\bar r_1(1-\bar N_1/K_1)}>\frac{\alpha_{23}}{\alpha_{13}}$, then $E_{2p}$ is globally stable {\rm \cite[Theorem 4.3(ii)]{Hsu1981}}.
\item If $\frac{\bar r_1-\alpha_{12}\bar N_2}{\bar r_2(1-\bar N_2/K_2)}<\frac{a_{13}}{\alpha_{23}}$ and $\frac{\bar r_2-\alpha_{21}\bar N_1}{\bar r_1(1-\bar N_1/K_1)}<\frac{\alpha_{23}}{\alpha_{13}}$, then {\rm  \cite[Theorem 4.4]{Hsu1981}}
\begin{enumerate}[{\rm (a)}]
\item The positive equilibrium $E_c=(N_{1c}, N_{2c}, P_c)$ of \eqref{Hsumodel} exists uniquely if and only if 
          \[ r_1r_2(1-1/K_1)(1-1/K_2)>(r_1-a_{12})(r_2-a_{21}). \]
\item If $E_c$ exists, then $E_c$ is a saddle point with two-dimensional stable manifold through $E_c$. In this case, if $E_{12}=(N_1^*, N_1^*, 0)$, $K_1>r_2/\alpha_{21}$, and $K_2>r_1/\alpha_{12}$, then $E_{12}$ is a saddle point with one-dimensional stable manifold through $E_{12}$
in the $N_1$-$N_2$ plane, while in the case $K_1<r_2/\alpha_{21}$ and $K_2<r_1/\alpha_{12}$, $E_{12}$ is a saddle point with
two-dimensional stable manifold in the $N_1$-$N_2$ plane.
\item If $E_c$ doesn’t exist, then $E_{12}$ exists with $K_1<r_2/\alpha_{21}$ and $K_2<r_1/\alpha_{12}$, and $E_{12}$ is a saddle point with a two-dimensional stable manifold through it.
\end{enumerate}
\item If $\frac{\bar r_1-\alpha_{12}\bar N_2}{\bar r_2(1-\bar N_2/K_2)}>\frac{a_{13}}{\alpha_{23}}$ and $\frac{\bar r_2-\alpha_{21}\bar N_1}{\bar r_1(1-\bar N_1/K_1)}>\frac{\alpha_{23}}{\alpha_{13}}$, then the positive equilibrium $E_c$ exists and the system \eqref{Hsumodel} is persistent  {\rm \cite[Theorem 4.5]{Hsu1981}}.
\end{enumerate}
\end{theorem}

It is easy to see that systems \eqref{Hsumodel} and \eqref{themodel} are equivalent with the following scaling,
\begin{align*}
 t &= D \bar t, \quad u_1=\frac{N_1}{K_1},  \quad u_2=\frac{N_2}{K_2}, \quad v=\frac{P}{D}, \\
r_1&=\frac{\bar r_1}{D}, \quad r_2=\frac{\bar r_2}{D}, \quad a_{12}=\frac{\alpha_{12}K_2}{D}, \quad a_{21}=\frac{\alpha_{21}K_1}{D}, \\
a_{13}&=\alpha_{13}, \quad a_{23}=\alpha_{23}, \quad a_{31}=\frac{\alpha_{31}K_1}{D}, \quad a_{32}=\frac{\alpha_{32}K_2}{D},
\end{align*}
and assumptions \eqref{Hsuassumption} are equivalent to our assumptions $a_{31}>1$ and $a_{32}>1$, respectively.

Comparing Theorem \ref{Hsuresults} with our results, Hsu only discusses the case where both $a_{31}$ and $a_{32}$ are greater than 1. Moreover, it is easy to see that the results (i)-(iv) of Theorem \ref{Hsuresults} correspond to Proposition \ref{EcaseIV-i} (IV-(i)), the D1 case of Proposition \ref{IV-(ii)-results}(ii), the D1 and D3 cases of Proposition \ref{IV-iii-bistability}-(i) and (iii), and the D1 and D3 cases of Proposition \ref{IV-ivE*}, respectively. However, our work also yields many additional novel results for all positive $a_{31}$ and $a_{32}$, which are investigated analytically or numerically. Additionally, Gilpin \cite{Gilpin1979} presented some complex behavior of \eqref{themodel} for a set of parameters, and we numerically present similar results in the cell IV-(iv)-B1 and show the bifurcation diagram in the next section.

\subsection{Biological Interpretations and Numerical Simulations}

We classify the global dynamics based on the essential inequalities derived from the relationships between $ C_1$, $C_2$, $C_3$, and $ C_4$ with assumption {\bf (H2)}, $C_1>C_2$, which means that the species $u_1$ wins the apparent competition with $u_2$ \cite{Holt1977, Holt2017}. Accordingly, we attempt to use two additional biological factors that interact to interpret the global dynamics: {\it negative} and {\it trade-off} competitions. After integrating all mathematical results in Table \ref{tab:long}, we will use these two indicators to categorize and explain our findings. Meanwhile, we try to answer the Paine conjecture \cite{Paine1966}, which suggests that an ecosystem's species diversity is directly related to its top predators.

Recall the negative competition explained in Remark \ref{negativecompetition}. We say that $u_1$ wins the negative competition with $u_2$ if $C_1 > C_3$, or equivalently, $r_1 > a_{12}$, which means that the $r$-strategy of $u_1$ can overcome the negative effect from the species $u_2$. Similarly, we say that $u_2$ wins the negative competition with $u_1$ if $C_2 > C_4$, or equivalently, $r_2 > a_{21}$. Note that this competition is independent of the predator $v$, so this item is actually two-dimensional competition between prey. A new term, "negative competition," is used instead of "resource competition," which is usually used to describe competition among species at the same trophic level for the same limiting resource. In general, there is only one winner in resource competition. However, in our consideration, there could be two winners in negative competition due to different biological or environmental niches. Hence, there are four generic relationships for $u_1$-$u_2$ in negative competition: win-loss, win-win, loss-loss, and loss-win, corresponding to the values of $C_1$-$C_4$ and named T1-T4 in the first and second columns of Table \ref{tab:long}.

On the other hand, for species $u_1$, there are three alternative choices to make the trade-off strategy: to increase its offspring with the $r$-strategy, to enhance the defense mechanism against $v$, to increase the negative effects for $u_2$, where the parameters $r_1$, $a_{13}$, and $a_{21}$ present the qualities of these strategies. We say that $u_1$ wins the trade-off competition if $C_1=r_1a_{23}>C_4=a_{21}a_{13}$, or equivalently, $\frac{r_1}{a_{13}}>\frac{a_{21}}{a_{23}}$. It can be interpreted that $u_1$ makes a trade-off between getting a higher ratio $\frac{r_1}{a_{13}}$ or a higher ratio $\frac{a_{21}}{a_{23}}$ by increasing its offspring, developing a stronger defense mechanism against $v$, or enhancing its competitive ability against $u_2$. This also leads to a decrease in $v$'s predation on $u_1$, causing $v$ to increase its predation on $u_2$. Similarly, we say that $u_2$ wins the trade-off competition if $C_2=r_2a_{13} > C_3=a_{12}a_{23}$, or equivalently, $\frac{r_2}{a_{23}}>\frac{a_{12}}{a_{13}}$. We present all possibilities in the third column, named Trade-Off, of Table \ref{tab:long}. 

Two points are worth noting. The negative and trade-off competitions are not independent because the parameters depend only on $C_1$-$C_4$ and assumption {\bf (H2)}. For example, in type T4, with loss-win in the negative competition, i.e., $C_1<C_3$ and $C_2>C_4$, and with $C_1>C_2$, we have $C_3>C_1>C_2>C_4$, which implies win-loss in the trade-off competition. The second point is that, under the trade-off strategy of $u_1$, increasing $a_{21}$ to enhance its competitive ability against $u_2$ appears ineffective in winning the trade-off competition. However, $u_1$ can win the negative competition if $a_{21}$ exceeds $r_2$, and we will see that the negative competition is the dominant factor.
 
\begin{figure}
\caption{Selected numerical simulations of \eqref{themodel} with parameters of type T1 on the negative competition and the initial condition $(0.68, 0.5, 0.5)$.}
\label{figT1}
\subfigure[The cell III-ii-A3 with parameters $r_1=0.7$, $r_2=0.9$, $a_{12}=0.5$, $a_{21}=2.9$, $a_{13}=1.5$, $a_{23}=3.9$, $a_{31}=0.85$, $a_{32}=0.25$.]{
\includegraphics[scale=0.32]{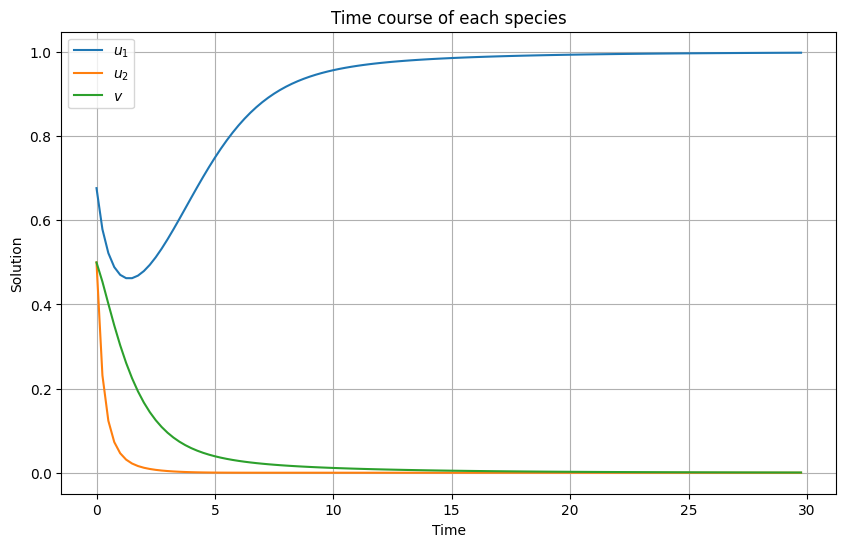}}
\quad
\subfigure[The cell III-ii-C3 with parameters $r_1=0.9$, $r_2=0.1$, $a_{12}=0.7$, $a_{21}=4.5$, $a_{13}=0.3$, $a_{23}=0.3$, $a_{31}=1.5$, $a_{32}=0.7$.]{
\includegraphics[scale=0.32]{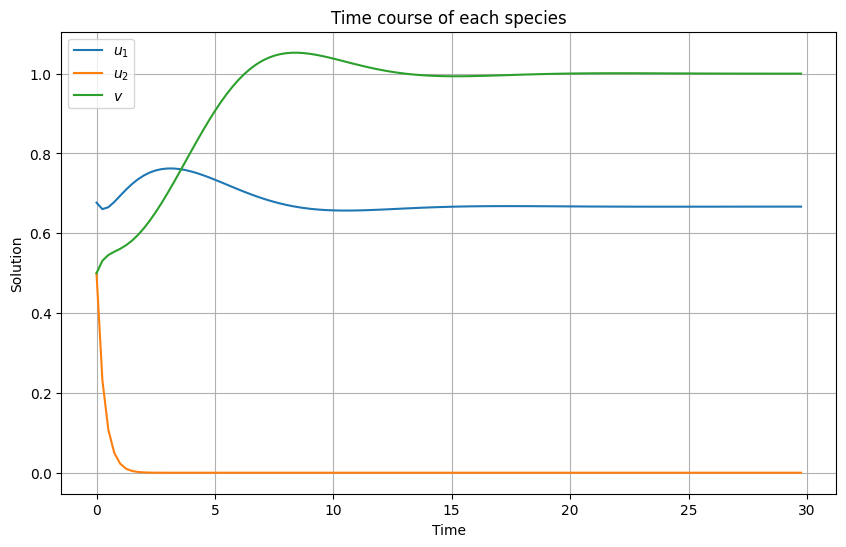}}
\\
\subfigure[The cell IV-iii-A4 with parameters $r_1=0.5$, $r_2=0.1$, $a_{12}=0.4$, $a_{21}=0.3$, $a_{13}=0.1$, $a_{23}=0.1$, $a_{31}=0.55$, $a_{32}=0.1$.]{
\includegraphics[scale=0.32]{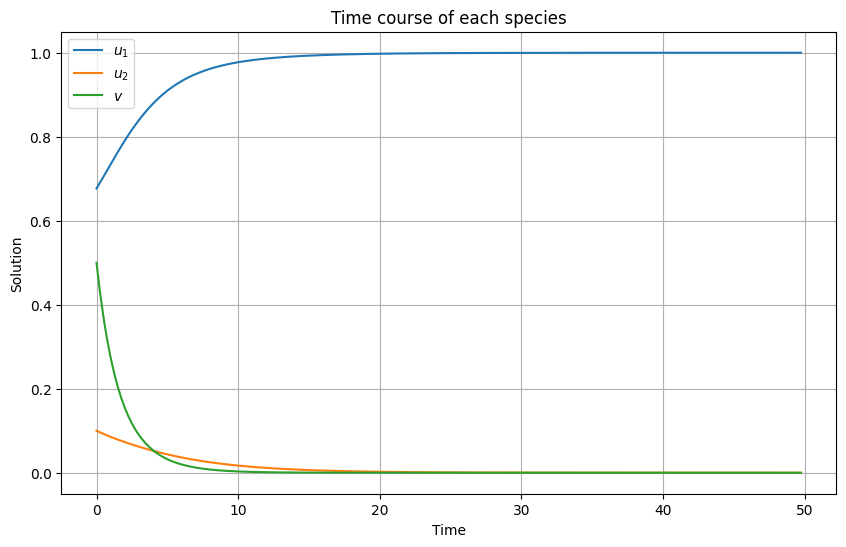}}
\quad
\subfigure[The cell IV-iii-C4 with parameters $r_1=1.5$, $r_2=0.3$, $a_{12}=0.9$, $a_{21}=2.7$, $a_{13}=0.3$, $a_{23}=4.7$, $a_{31}=1.5$, $a_{32}=0.3$.]{
\includegraphics[scale=0.32]{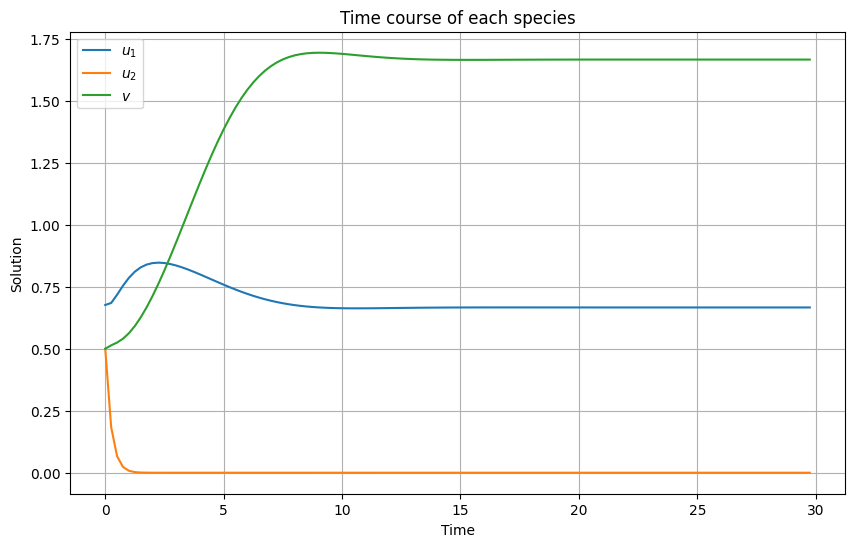}}
\end{figure}

Under the superior apparent competition of $u_1$, we summarize all possibilities of our findings in Table \ref{tab:long}. First, for type T1, the boundary equilibrium $E_1$ or $E_{1v}$ is GAS, corresponding to categories I, II-(i), ..., III-(ii)-C3. All these results are verified analytically, and we conclude that negative competition dominates the global dynamics of \eqref{themodel} under the assumption {\bf (H2)}, regardless of the trade-off competition. That is, all dynamics are consistent with the results of the two-dimensional subsystem \eqref{u1u2system}. The selected numerical simulations are presented in Figure \ref{figT1}.




\begin{table}[h]
\begin{threeparttable}
\caption{Classifications with assumption {\bf (H2)} $C_1>C_2$.} 
    \label{tab:long}
\centering
{\footnotesize
\begin{tabular}{|c|cc|cc|l|l|}
\hline\hline\hline
\multirow{3}{*}{\textbf{}} 
& \multicolumn{2}{|c|}{\textbf{Negative}} 
& \multicolumn{2}{c|}{\textbf{Trade-Off}} 
& \multirow{3}{*}{Dynamics} & \multirow{3}{*}{Classifications}\\
&\multicolumn{2}{c|}{Competition} &\multicolumn{2}{c|}{Competition}&&\\
\cline{2-5}
 & $u_1$ & $u_2$ & $u_1$ & $u_2$ &&\\
\hline
\hline
\hline
\multirow{7}{*}{T1}  & \multirow{7}{*}{Win}  & \multirow{7}{*}{Loss} & 
 \multirow{4}{*}{Win} & \multirow{4}{*}{Loss} & \multirow{3}{*}{$E_1$ GAS} &IV-(ii)-A4 and A5 \\
&&&&& &IV-(iii)-A4 and A5\\
&&&&&&IV-(iv)-A4, A5, B4 and B5 \\
\cline{6-7}
&&&&&$E_{1v}$ GAS&IV-(iii)-C4 and C5\\
\cline{4-7}
&&&\multirow{1}{*}{Win}& \multirow{1}{*}{Win} & $E_1$ GAS & II-(ii)-A1 and B1\\
\cline{4-7}
\cline{4-7}
  &   &  & \multirow{2}{*}{Loss} & \multirow{2}{*}{Loss} &$E_1$ GAS& III-(ii)-A3 \\
  \cline{6-7}
  &   &  &  &  &$E_{1v}$ GAS& III-(ii)-C3 \\
\hline
\hline
\multirow{6}{*}{T2}  & \multirow{6}{*}{Win}  & \multirow{6}{*}{Win}  & 
 \multirow{4}{*}{Win} & \multirow{4}{*}{Loss} &$E_{12}$ GAS& IV-(ii)-A3 and C3\\
\cline{6-7}
&&&&& \red{$E_{1v}$/$E_{12}$}\tnote{1} & \red{IV-(iii)-C3}\\
\cline{6-7}
&&&&& $E_{12}$/$E^*$\tnote{2}& IV-(iv)-A3, B3 and C3\\
\cline{6-7}
&&&&& $E^*$ & IV-(iv)-D3\\
\cline{4-7}
&&&\multirow{3}{*}{Win} & \multirow{3}{*}{Win} &
\multirow{2}{*}{$E_{12}$/$E^*$\tnote{2}}& II-(ii)-A2, A3, B2, \\
&&&&&& B3, C2 and C3\\
\cline{6-7}
&&&&& $E^*$& II-(ii)-D2 and D3\\
\hline 
\hline
\multirow{11}{*}{T3}  & \multirow{11}{*}{Loss} & \multirow{11}{*}{Loss} & \multirow{6}{*}{Win} & \multirow{6}{*}{Loss} &\multirow{2}{*}{Bis $E_1$/$E_2$} & IV-(ii)-A2\\
&&&&&& IV-(iii)-A2\\
\cline{6-7}
&&&&&\multirow{2}{*}{Bis $E_1$/$E_{2v}$} & IV-(ii)-B2\\
&&&&&& IV-(iii)-B2\\
\cline{6-7}
&&&&&Bis $E_{1v}$/$E_2$ & IV-(iii)-C2\\
\cline{6-7}
&&&&&Bis $E_{1v}$/$E_{2v}$& IV-(iii)-D2\\
\cline{6-7}
&&&&& \red{$E_1$ GAS}\tnote{3}& \red{IV-(iv)-B2}\\
\cline{4-7}
  &  &  & \multirow{4}{*}{Loss} & \multirow{4}{*}{Loss} & Bis $E_1$/$E_2$& III-(ii)-A1 and A2\\
  \cline{6-7}
&&&&& Bis $E_1$/$E_{2v}$ &  III-(ii)-B1 and B2\\
\cline{6-7}
&&&&& Bis $E_{1v}$/$E_2$& III-(ii)-C1 and C2\\
\cline{6-7}
&&&&& Bis $E_{1v}$/$E_{2v}$& III-(ii)-D1 and D2\\
\hline
\hline
\multirow{6}{*}{T4} & \multirow{6}{*}{Loss} & \multirow{6}{*}{Win}  & \multirow{6}{*}{Win} & \multirow{6}{*}{Loss} &$E_2$ GAS& IV-(ii)-A1 and C1\\
&&&&& $E_{2v}$ GAS &  IV-(ii)-B1 and D1\\
\cline{6-7}
&&&&& \red{Bis $E_{1v}$/$E_2$} & \red{IV-(iii)-C1}\\
&&&&& \red{Bis $E_{1v}$/$E_{2v}$} & \red{IV-(iii)-D1}\\
\cline{6-7}
&&&&& \red{Complex}\tnote{4} & \red{IV-(iv)-B1}\\
&&&&& \red{$E^*$} & \red{IV-(iv)-D1}\\
\hline
\hline
\end{tabular}
}
\begin{tablenotes}

\item[1] If \eqref{E12stable} holds, then the bistabilty of $E_{1v}/E_{12}$ happens. Otherwise, if \eqref{E12unstable} holds, then $E_{1v}$ is GAS.
\item[2] If \eqref{E12stable} holds, then $E_{12}$ is GAS. Otherwise, if \eqref{E12unstable} holds, then $E^*$ is stable, and system \eqref{themodel} is uniformly persistent.
\item[3] By analytically proof, $E_1$ is GAS if \eqref{E12unstable} holds. However, if \eqref{E12stable} holds, then $E_1$ is also GAS by numerical observations.
\item[4] In this cell, we observe numerically that $E^*$ is GAS, $E^*$ is unstable, and there is a limit cycle, or very complicated behavior.
\end{tablenotes}
\end{threeparttable}
\end{table}

\begin{figure}
\caption{Selected numerical simulations of \eqref{themodel} with parameters of type T2 on the negative competition and the initial condition $(0.68, 0.5, 0.5)$. Parameters of these two simulations are almost the same except for $a_{32}=0.45$ in the left panel and $a_{32}=0.85$ in the right panel, which satisfy \eqref{E12stable} and \eqref{E12unstable}, respectively.}
\label{figT2-1}
\subfigure[The cell II-(ii)-A2 with parameters $r_1=0.3$, $r_2=0.5$, $a_{12}=0.1$, $a_{21}=0.35$, $a_{13}=0.1$, $a_{23}=0.3$, $a_{31}=0.85$, $a_{32}=0.45$.]{
\includegraphics[scale=0.3]{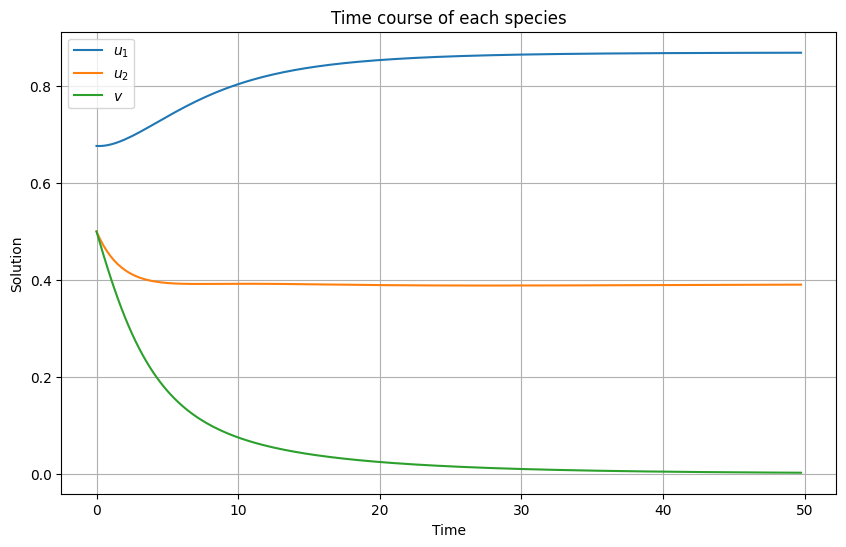}}
\qquad
\subfigure[The cell II-(ii)-A2 with parameters $r_1=0.3$, $r_2=0.5$, $a_{12}=0.1$, $a_{21}=0.35$, $a_{13}=0.1$, $a_{23}=0.3$, $a_{31}=0.85$, $a_{32}=0.85$.]{
\includegraphics[scale=0.3]{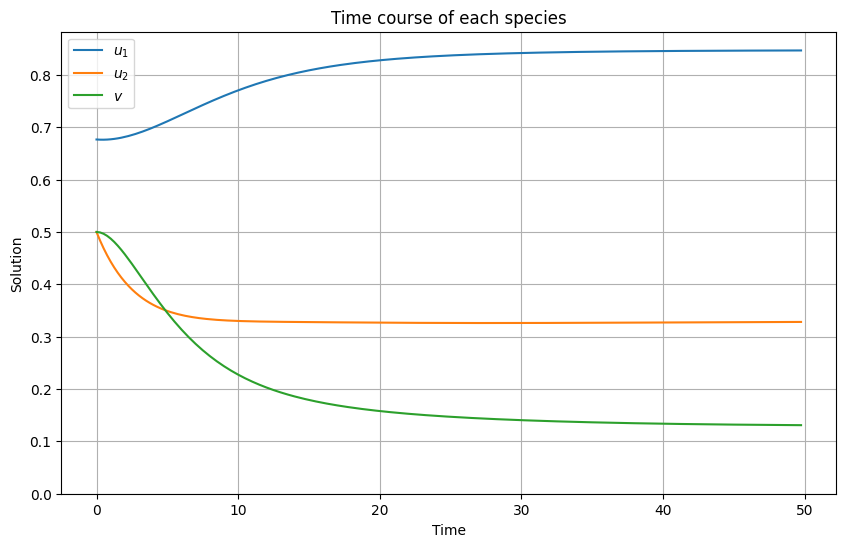}}
%
\end{figure}

For type T2 of the win-win state in the negative competition for system \eqref{themodel}, compared with the results of \eqref{u1u2system}, we expect that $E_{12}$ is stable if it exists and \eqref{E12stable} holds. Alternatively, $E^*$ is stable if it exists and \eqref{E12unstable} holds. These results appear in type T2 of Table \ref{tab:long}, except for case IV-(iii)-C3. We have shown analytically that $E_{12}$ is GAS if it exists and is stable. However, we can only verify the local stability of $E^*$ when it exists. The selected numerical simulations are presented in Figure \ref{figT2-1}. Parameters for two cells, II-ii-A2 and IV-iv-A3, are used to illustrate that $E_{12}$ or $E^*$ is globally asymptotically stable, depending on whether inequality \eqref{E12stable} or \eqref{E12unstable} holds. In panel (a) of Figure \ref{figT2-1}, with the parameters $r_1=0.3$, $r_2=0.5$, $a_{12}=0.1$, $a_{21}=0.35$, $a_{13}=0.1$, $a_{23}=0.3$, $a_{31}=0.85$, $a_{32}=0.45$, we check that $E_{12}$ exists and \eqref{E12stable} holds. By contrast, in panel (b) of Figure \ref{figT2-1}, $E_{12}$ exists and \eqref{E12unstable} holds only when $a_{32}$ is changed from $0.45$ to $0.85$. Similar simulations can be performed for the cells II-(ii)-A2, A3, B2, B3, C3, C3 and IV-(iv)-A3, B3, C3. This exception will be discussed later.

\begin{figure}
\caption{In the left panel, selected numerical simulations of \eqref{themodel}  with parameters of type T3 with initial conditions $(0.8, 0.2, 0.2)$ in red color and $(0.2, 0.8, 0.2)$ in green color, respectively. In the right panel, selected numerical simulations of \eqref{themodel} with parameters of type T3 with initial conditions $(0.9, 0.1, 0.1)$ in red color and $(0.1, 0.9, 0.1)$ in green color, respectively.}
\label{figT3}
\subfigure[The cell IV-ii-A2 with parameters $r_1=0.3$, $r_2=0.3$, $a_{12}=2.9$, $a_{21}=4.5$, $a_{13}=0.1$, $a_{23}=3.7$, $a_{31}=0.9$, $a_{32}=0.7$.]{
\includegraphics[scale=0.3]{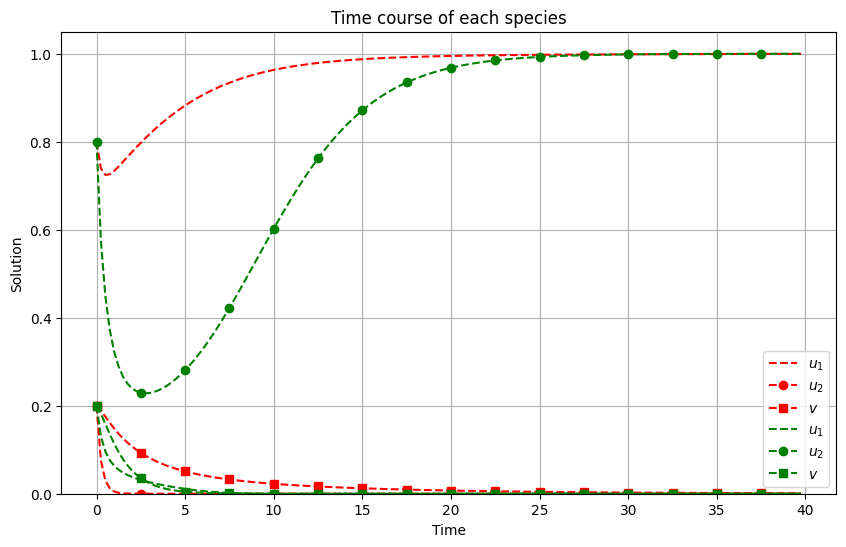}}
\qquad
\subfigure[The cell IV-iii-C2 with parameters $r_1=0.5$, $r_2=1.9$, $a_{12}=3.1$, $a_{21}=3.7$, $a_{13}=0.1$, $a_{23}=4.5$, $a_{31}=1.5$, $a_{32}=0.7$.]{
\includegraphics[scale=0.3]{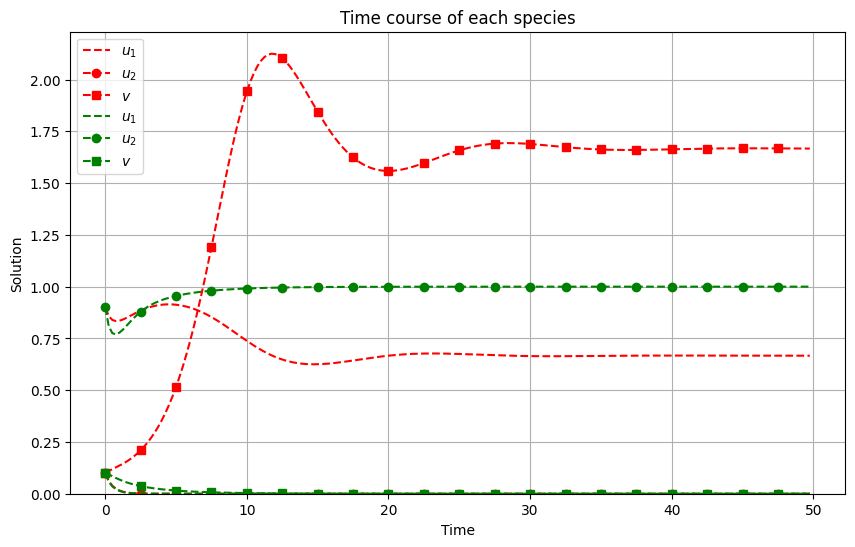}}
\end{figure}

\begin{figure}
\caption{In the left panel, selected numerical simulations of \eqref{themodel}  with parameters of type T3 with initial conditions $(0.9, 0.1, 0.1)$ in red color and $(0.1, 0.9, 0.1)$ in green color, respectively. In the right panel, selected numerical simulations of \eqref{themodel} with parameters of type T3 with initial conditions $(0.9, 0.1, 0.1)$ in red color and $(0.1, 0.9, 0.1)$ in green color, respectively.}
\label{figIII-iibc}
\subfigure[The cell III-ii-B1 with parameters $r_1=0.1$, $r_2=0.1$, $a_{12}=0.3$, $a_{21}=0.3$, $a_{13}=0.7$, $a_{23}=0.9$, $a_{31}=0.9$, $a_{32}=1.7$.]{
\includegraphics[scale=0.3]{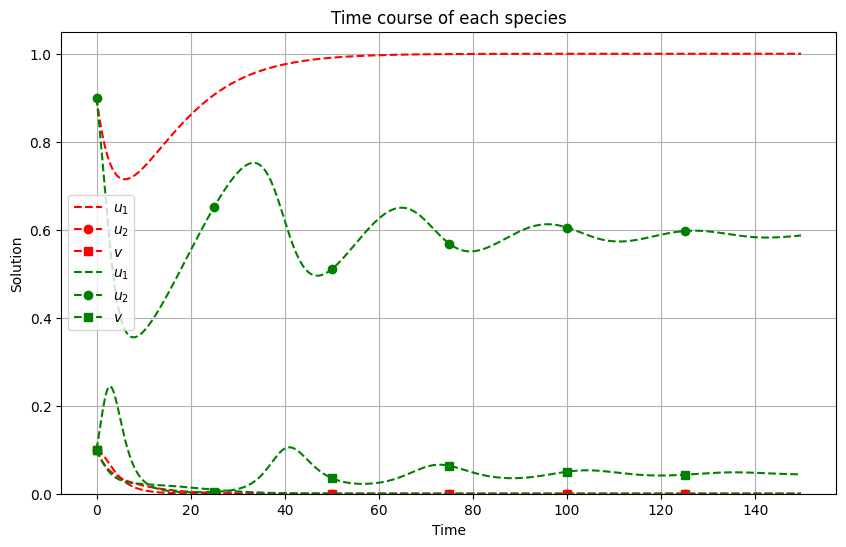}}
\qquad
\subfigure[The cell III-ii-C2 with parameters $r_1=0.1$, $r_2=0.5$, $a_{12}=0.5$, $a_{21}=3.3$, $a_{13}=0.7$, $a_{23}=4.5$, $a_{31}=1.5$, $a_{32}=0.7$.]{
\includegraphics[scale=0.3]{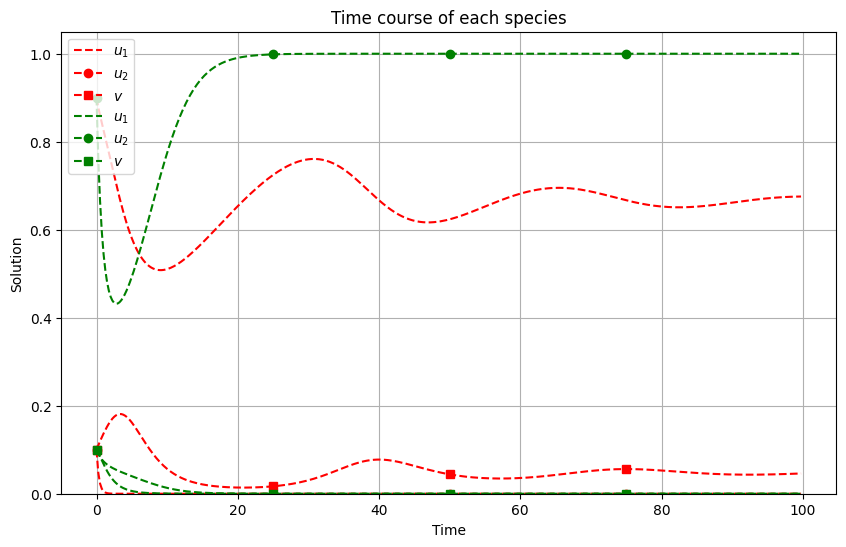}}
\end{figure}

For the T3 loss-loss state of the negative competition in Table \ref{tab:long}, as in \eqref{u1u2system}, we expect bistability. Indeed, there are four types of bistability, $E_1/E_2$, $E_1/E_{2v}$, $E_{1v}/E_2$, and $E_{1v}/E_{2v}$, which depend on the values of $a_{31}$ and $a_{32}$. The selected numerical simulations are presented in Figures \ref{figT3} and \ref{figIII-iibc}. Notably, the only exception is case IV-(iv)-B2. We show analytically that if \eqref{E12unstable} holds, then $E_1$ is GAS, and we find numerically that if \eqref{E12stable} holds, then solutions of \eqref{themodel} also converge to $E^*$. This exception will be discussed later.


\begin{figure}
\caption{Selected numerical simulations of \eqref{themodel} with parameters of type T4 on the negative competition and the initial condition $(0.7, 0.5, 0.5)$.}
\label{figIV-iibc}
\subfigure[The cell IV-ii-B1 with parameters $r_1=0.1$, $r_2=0.3$, $a_{12}=0.3$, $a_{21}=0.1$, $a_{13}=0.3$, $a_{23}=3.5$, $a_{31}=0.6$, $a_{32}=1.3$.]{
\includegraphics[scale=0.3]{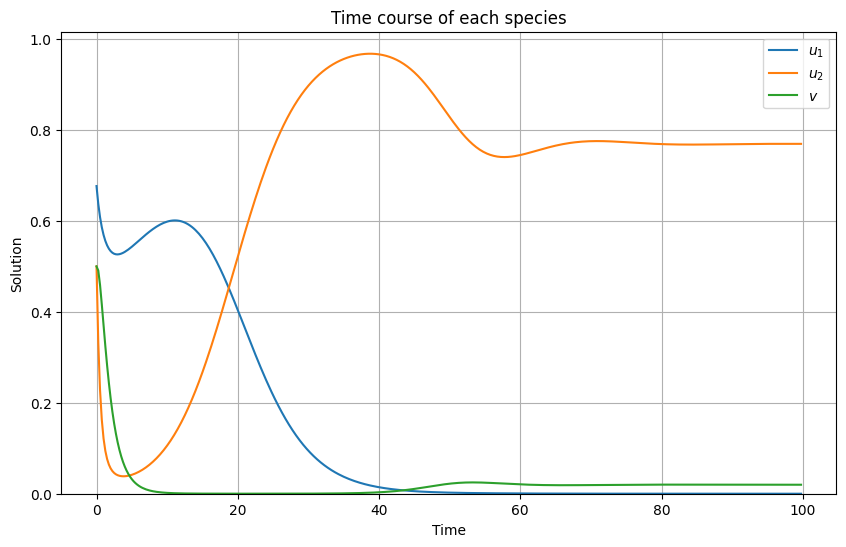}}
\qquad
\subfigure[The cell IV-ii-C1 with parameters $r_1=0.1$, $r_2=0.3$, $a_{12}=0.3$, $a_{21}=0.1$, $a_{13}=0.9$, $a_{23}=3.9$, $a_{31}=1.9$, $a_{32}=0.7$.]{
\includegraphics[scale=0.3]{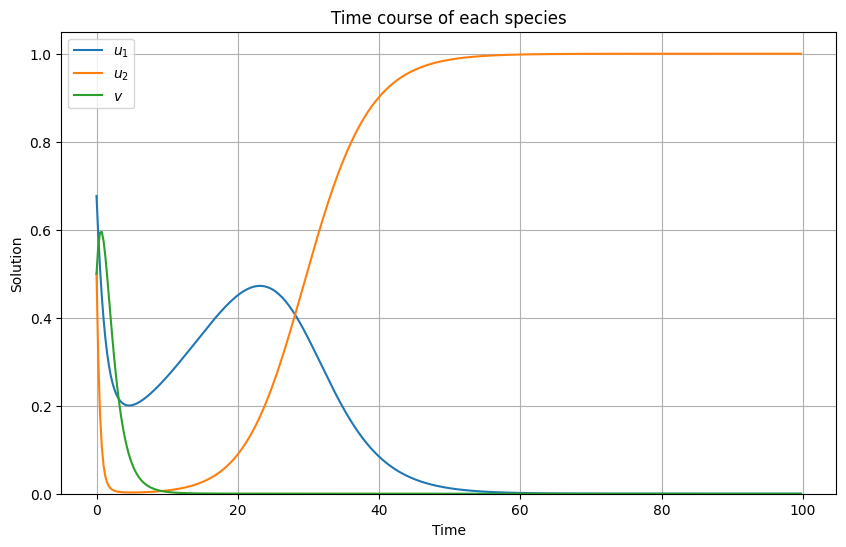}}
\end{figure}

\begin{figure}
\caption{In the left panel, selected numerical simulations of \eqref{themodel}  with parameters of type T4 with initial conditions $(0.9, 0.1, 0.1)$ in red color and $(0.1, 0.9, 0.1)$ in green color, respectively. In the right panel, selected numerical simulations of \eqref{themodel} with parameters of type T3 with initial conditions $(0.9, 0.1, 0.1)$ in red color and $(0.1, 0.9, 0.1)$ in green color, respectively.}
\label{figIV-iiic1d1}
\subfigure[The cell IV-iii-C1 with parameters $r_1=0.1$, $r_2=0.3$, $a_{12}=0.3$, $a_{21}=0.1$, $a_{13}=0.1$, $a_{23}=2.3$, $a_{31}=1.5$, $a_{32}=0.5$.]{
\includegraphics[scale=0.3]{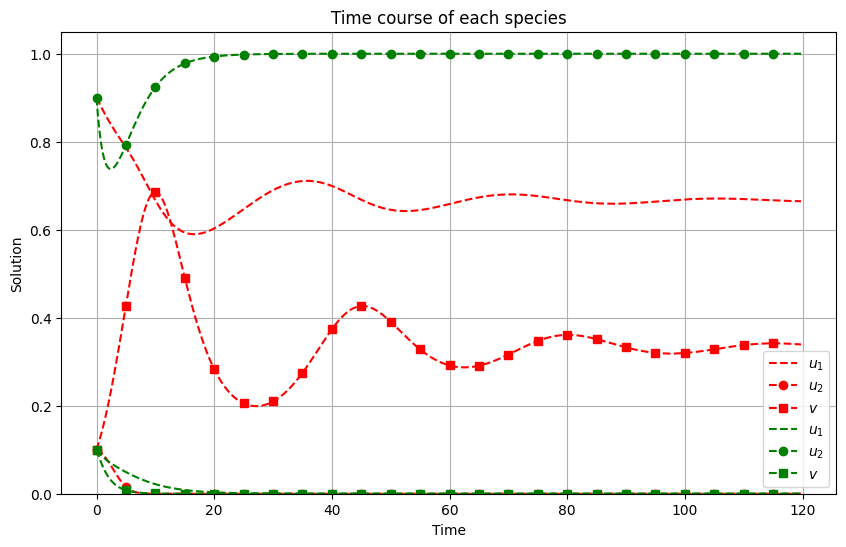}}
\qquad
\subfigure[The cell IV-iii-D1 with parameters $r_1=0.1$, $r_2=0.3$, $a_{12}=0.3$, $a_{21}=0.1$, $a_{13}=0.1$, $a_{23}=0.9$, $a_{31}=1.5$, $a_{32}=2.9$.]{
\includegraphics[scale=0.3]{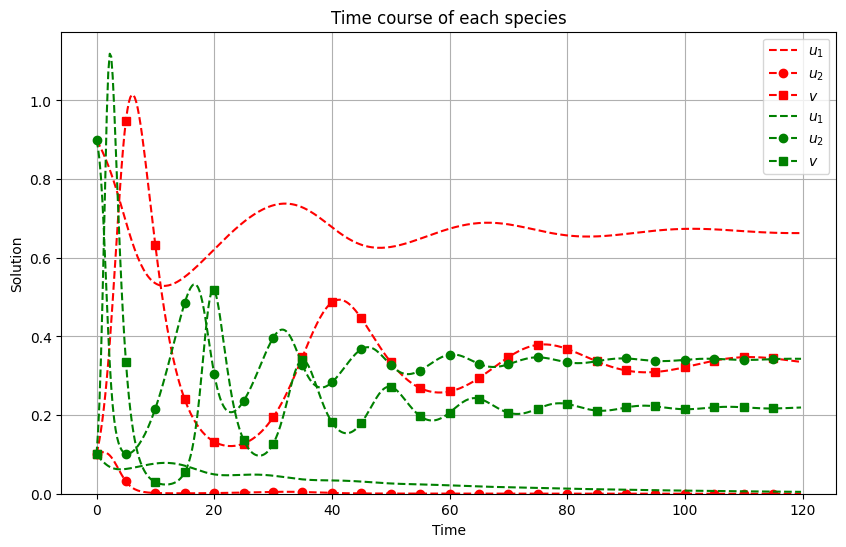}}
\end{figure}

\begin{figure}
\caption{Selected numerical simulations of \eqref{themodel}  for the cell IV-iv-B1 with parameters, $r_1=0.901$, $r_2=0.902$, $a_{12}=1.303$, $a_{21}=0.904$, $a_{13}=0.025$, $a_{23}=0.106$, $a_{31}=0.5$, $a_{32}=4.5$.}
\label{figIV-ivB1}
\subfigure[Time courses of each species for the cell IV-iv-B1.]{
\includegraphics[scale=0.4]{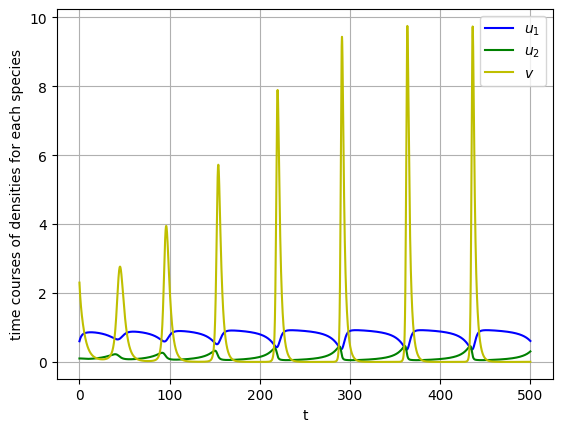}}
\qquad
\subfigure[3D plot of a trajectory for the cell IV-iv-B1.]{
\includegraphics[scale=0.4]{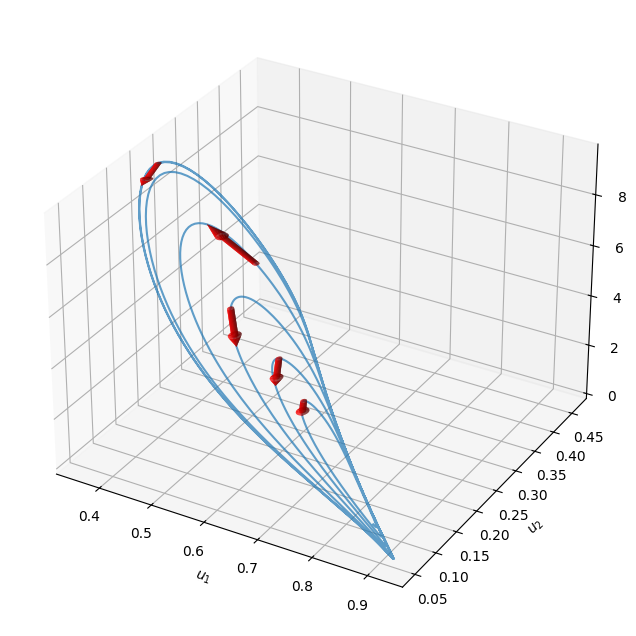}}
\end{figure}

\begin{figure}
\caption{Bifurcation diagram of simulations of \eqref{themodel}  for the cell IV-iv-B1 with parameters, $r_1=0.901$, $r_2=0.908$, $a_{12}=1.303$, $a_{21}=0.904$, $a_{13}=0.02$, $a_{23}=0.206$, $a_{31}=0.5$, using bifurcation parameter $a_{32}$ from 3.0 to 8 with step size 0.002.}
\label{figIV-ivB1-bif}
\begin{center}
\includegraphics[scale=0.5]{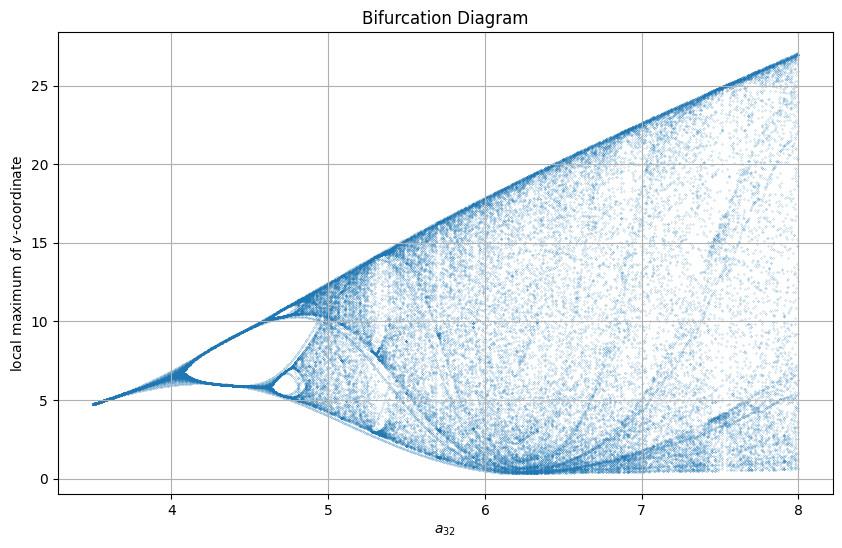}
\end{center}
\end{figure}

For the type T4 loss-win state in the negative competition, compared with the results of \eqref{u1u2system}, we expect that $E_2$ or $E_{2v}$ is GAS, which depends on whether $a_{32}<1$ or $a_{32}>1$. You can see that, in Table \ref{tab:long}, it is true only for the second row of case IV-(ii).  Numerical simulations of cells IV-ii-B1 and IV-ii-C1 are presented in Figure \ref{figIV-iibc}.
There are four exceptions in type T4. Numerically, we observe the bistability $E_{1v}/E_2$ or $E_{1v}/E_{2v}$ in the cells IV-(iii)-C1 and IV-(iii)-D1, respectively. Please refer to Figure \ref{figIV-iiic1d1}.
Moreover, we show the existence and stability of $E^*$ for case IV-(iv)-D1, and we conjecture that coexistence is globally asymptotically stable. Numerical simulations indicate that periodic solutions or chaos occur in cell IV-(iv)-B1. Selected numerical simulations of a periodic solution are presented in Figure \ref{figIV-ivB1}. In Figure \ref{figIV-ivB1-bif}, we plot a bifurcation diagram for parameters in cell IV-(iv)-B1, varying $a_{32}$ from 3 to 8 in steps of 0.002. For each fixed parameter set, the trajectory of \eqref{themodel} is simulated for 8000 time units, and the local maxima of the $v$-coordinate are projected to form Figure \ref{figIV-ivB1-bif}. Based on numerical observations, the period-doubling route to chaos appears likely. However, a detailed discussion of this issue is outside the scope of this work.

Finally, we discuss some exceptions in Table \ref{tab:long}, highlighted in red. Here, the exceptions indicate that the dynamics of \eqref{themodel} differ from those of \eqref{u1u2system} without the top predator $v$ involved and sharing the same parameters. We try to answer the following question related to the Paine conjecture: how does the top predator interact with two competing prey species to affect the diversity of the whole system? 

First, we observe that the exception occurs only in categories IV-(iii) and IV-(iv), and those are all the win-loss states of the trade-off competition. To clarify the factors underlying the different dynamics, we list all competitions of three species: trade-off competition, apparent competition, and extensive negative competition. We say that $u_1$ wins the {\it extensive negative competition} which generalizes the negative competition of $u_1$-$u_2$ submodule if 
\begin{align}\label{u1venc}
r_1-a_{12}\bar u_2-a_{13}\bar v_2>0.
\end{align}
This means that species $u_1$ can withstand the negative effects from the balance states of the $u_2$-$v$ submodule. And it is clear that \eqref{u1venc} is the converse of \eqref{E2vstable} and is also equivalent to $\frac{C_1-C_2}{C_3-C_2}>\frac{1}{a_{32}}$. Hence, if $u_1$ wins the extensive negative competition, then $E_{2v}$ is unstable. Similarly, we say that $u_2$ wins the extensive negative competition if $r_2-a_{21}\bar u_1-a_{23}\bar v_1>0$ which is the converse of \eqref{E1vstable} and is also equivalent to $\frac{C_1-C_2}{C_1-C_4}<\frac{1}{a_{31}}$. Moreover, if $u_2$ wins the extensive negative competition, then $E_{1v}$ is unstable


\begin{enumerate}[(i)]
\item IV-(iii)-C3: In this case, this is the only exception for the win-win state of Type 2. For system \eqref{u1u2system}, the equilibrium $E_{12}$ is GAS on the $u_1$-$u_2$ subcommunity. And $v$ can invade if \eqref{E12unstable} holds, since $v$ can overcome its mortality by consuming $u_1$ and $u_2$ simultaneously. On the other hand, if \eqref{E12stable} holds, then $v$ cannot invade. So we expect that for type T2, the dynamics of all cases are either  $E_{12}$ is GAS if \eqref{E12stable} holds or $E^*$ is GAS if \eqref{E12unstable} holds, if the negative competition dominates the dynamics of model \eqref{themodel}. However, for case IV-(iii), $E^*$ is unstable if it exists by Lemma \ref{IV-iiiE*unstable} since \eqref{RH-2-1} does not hold. So it is interesting what happens for cell IV-(iii)-C3 when \eqref{E12unstable} holds.

Mathematically, for cell IV-(iii)-C3, we show that bistability $E_{1v}/E_{12}$ happens if \eqref{E12stable} holds and $E_{1v}$ is GAS if  \eqref{E12unstable} holds. This implies that $E^*$ does not exist and the role of $E^*$ is replaced by $E_{1v}$. Scenario-wise, if $v$ invades the $u_1$-$u_2$ only competing sub-community and \eqref{E12unstable} holds, then $v$ can survive. Then $u_2$ cannot persist because it not only loses the apparent competition but also makes the wrong trade-off strategy because of $C_2=r_2a_{13}<a_{12}a_{23}=C_3$. This means that $u_2$ makes an effort to resist the negative effect from $u_1$ and to escape the predation of $v$. This is the first case to show that the top predator $v$ affects the dynamics of competing prey after invading.


\item IV-(iv)-B2 :  In this case, this is the only exception for the loss-loss state of Type 3. For submodule \eqref{u1u2system}, the bistability phenomenon was established. Therefore, for system \eqref{themodel} with the top predator involved, we expected it to also be bistable in the form of $u_1$/$u_2$. Because $a_{31}<1$ implies that $E_{1v}$ does not exist, and $E_{2v}$ is saddle with the classification conditions $\frac{1}{a_{32}}<\frac{C_1-C_2}{C_3-C_2} <1$. However, we show that if \eqref{E12unstable} holds, then $\lim_{t\to\infty}u_2(t)=0$ and $E_1$ is GAS. Numerical observations also suggest that $E_1$ is GAS even when \eqref{E12stable} holds.

For cell IV-(iv)-B2, the key lies in the interplay between the predator preference and the trade-off competition. System \eqref{themodel} can be viewed as the combination of the prey-prey competition subsystem \eqref{u1u2system} and the predator-prey interactions. In the loss-loss state of negative competition, both prey species suffer from direct competition, and without the predator, the system exhibits bistability: the outcome depends on initial conditions, with either $E_1$ or $E_2$ being stable. 

However, when the predator is introduced with $a_{31} < 1 < a_{32}$, the predator preferentially consumes $u_2$ (since $a_{32} > 1$ means $u_2$ can support a positive predator population in the $u_2$-$v$ subsystem). This preferential predation on $u_2$ reduces its population below the threshold needed to resist competitive exclusion by $u_1$. On the other hand, since $a_{31} < 1$, the predator cannot invade the $u_1$ one-dimensional subsystem, meaning $u_1$ can persist even in the absence of $u_2$. Moreover, like the previous case, $u_2$ makes the wrong trade-off strategy. Hence, whether or not the condition \eqref{E12unstable} holds, due to the asymmetric predation preferences and the trade-off structure ($C_2< C_4$ disfavors $u_2$ in the trade-off and $a_{32} > 1$ hurts $u_2$ through increased predation), the net effect is that $u_2$ is driven to extinction while $u_1$ persists at its carrying capacity. This explains why $E_1$ is GAS rather than bistability: the predator effectively mediates the competitive exclusion of $u_2$ by $u_1$ through preferential predation, overriding the bistable dynamics expected from the two-dimensional subsystem. This is the second case to show that the top predator $v$ affects the dynamics of competing prey after invading.

\item IV-(iii)-C1/D1: In this case, this is an exception for the loss-win state of Type 4. Similarly, for system \eqref{u1u2system}, the equilibrium $E_2$ is GAS on the $u_1$-$u_2$ subcommunity since $u_1$ loses and $u_2$ wins in the negative competition (i.e., $C_1 < C_3$ and $C_2 > C_4$). We expect that $E_2$ (or $E_{2v}$ if $a_{32}>1$) is GAS for system \eqref{themodel} as well. However, for cells C1 and D1 of case IV-(iii), we show bistability: $E_{1v}/E_2$ for C1 and $E_{1v}/E_{2v}$ for D1. This means that $u_2$ does not completely exclude $u_1$ as expected; instead, the predator $v$ plays a crucial role in maintaining $u_1$ through predation on $u_2$. We think that the reasons are threefold. One is that $u_1$ wins the apparent competition, and the other is that $u_1$ wins the trade-off competition. The final, and most important, reason is that $u_1$ and $u_2$ both win extensive negative competition, with $\frac{C_1-C_2}{C_1-C_4} > \frac{1}{a_{31}}$ and $\frac{C_1-C_2}{C_3-C_2} < \frac{1}{a_{32}}$, respectively.

For cells C1 and D1 of case IV-(iii), the advantages of $u_1$ and $u_2$ in the extensive negative competition mean that species $u_1$ and $u_2$ can be sustained despite the negative effect of the top predator $v$ and competing species. In cell C1, with $a_{32} < 1 < a_{31}$, this creates a scenario where the outcome depends on initial conditions: if the predator population is sufficiently large (the initial is close to $E_{1v}$), it can suppress $u_2$ enough to allow $u_1$ to coexist with the predator; if the predator population is small (the initial is close to $E_{2}$), $u_2$ can exclude $u_1$ as in the two-dimensional subsystem.

In cell D1, we have $1 < a_{31}$ and $1 < a_{32}$, meaning the predator can invade both single-prey subsystems. The bistability arises because the predator preferentially exploits $u_2$ (which is the superior competitor in the absence of the predator), thereby relieving the competitive pressure on $u_1$ and allowing it to persist. 

\item IV-(iv)-B1/D1: In this case, this is another exception for the loss-win state of Type 4. For cells B1 and D1 of case IV-(iv), we observe different dynamics: in cell B1, numerical simulations suggest the existence of periodic solutions or a stable positive equilibrium $E^*$, indicating complex dynamics shown in Figure \ref{figIV-ivB1-bif}; in cell D1, we show that $E^*$ exists and the system \eqref{themodel} is uniformly persistent. It is also stable by numerical observation. These results are quite different from the expected simple convergence to $E_2$ or $E_{2v}$.

For cells D1 and B1 of case IV-(iv), the conditions exclude the stability of all boundary equilibria. Therefore, we show that the system \eqref{themodel} is uniformly persistent. This is the only analytical result that we can show.

\end{enumerate}

\subsection{Conclusion} In this work, we have completed a comprehensive global analysis of the two-prey-one-predator Lotka-Volterra model with direct prey competition. Through a complete classification of the parameter space into 56 subitems across four generic categories, we have identified four essential factors governing the system's dynamics: (1) \textit{apparent competition}, which determines the superior prey through the ratio of birth rate to consumption rate (assumption (H2)); (2) \textit{negative competition}, representing the direct resource competition between prey species; (3) \textit{trade-off competition}, capturing the strategic choices prey make between reproduction, defense against predators, and competitive ability; and (4) \textit{extensive negative competition}, which extends the notion of negative competition to the full three-dimensional system by incorporating the predator's effect. Our main result demonstrates that, despite the complexity arising from the interplay of these four mechanisms, the global dynamics of the three-dimensional system are largely consistent with those of the associated two-dimensional prey-competition subsystem, with only six notable exceptions. These exceptions, discussed in detail in Section 6.2, arise from the intricate interaction between predator preference and the competitive structures, leading to dynamical phenomena such as bistability, periodic solutions, and chaos that differ from the predictions of the two-dimensional subsystem. Notably, these exceptional cases provide mathematical support for the Paine conjecture, illustrating how top predators can promote species diversity by mediating competition between prey species and preventing competitive exclusion.

\bibliographystyle{plain}

\end{document}